\documentclass[12pt]{amsart}
\usepackage{mathrsfs}
\usepackage{cases}
\usepackage{epic}
\usepackage{amsfonts}
\usepackage{graphicx}
\usepackage{amsmath}
\usepackage{amssymb, upgreek}
\usepackage{bm}
\usepackage{latexsym,todonotes}
\usepackage{pdflscape}
\usepackage[all]{xypic}
\usepackage[all]{xy}
\usepackage{color}
\usepackage{colordvi}
\usepackage{multicol}
\usepackage[normalem]{ulem}
\input diagxy
\usepackage[linktocpage=true]{hyperref}
\hypersetup{colorlinks,linkcolor=blue,urlcolor=cyan,citecolor=blue}

\begin{document}
\input xy
\xyoption{all}

\allowdisplaybreaks
%Theorem for the introduciton
\newtheorem{innercustomthm}{{\bf Theorem}}
\newenvironment{customthm}[1]
  {\renewcommand\theinnercustomthm{#1}\innercustomthm}
  {\endinnercustomthm}

  \newtheorem{innercustomcor}{{\bf Corollary}}
\newenvironment{customcor}[1]
  {\renewcommand\theinnercustomcor{#1}\innercustomcor}
  {\endinnercustomthm}

  \newtheorem{innercustomprop}{{\bf Proposition}}
\newenvironment{customprop}[1]
  {\renewcommand\theinnercustomprop{#1}\innercustomprop}
  {\endinnercustomthm}

\newcommand{\ev}{\bar{0}}
\newcommand{\odd}{\bar{1}}

\newcommand{\sub}{\operatorname{sub}\nolimits}

\newcommand{\fac}{\operatorname{fac}\nolimits}

\newcommand{\iadd}{\operatorname{iadd}\nolimits}

\renewcommand{\mod}{\operatorname{mod}\nolimits}
\newcommand{\proj}{\operatorname{proj}\nolimits}
\newcommand{\inj}{\operatorname{inj}\nolimits}
\newcommand{\rad}{\operatorname{rad}\nolimits}
\newcommand{\Span}{\operatorname{Span}\nolimits}
\newcommand{\soc}{\operatorname{soc}\nolimits}
\newcommand{\ind}{\operatorname{inj.dim}\nolimits}
\newcommand{\Ginj}{\operatorname{Ginj}\nolimits}
\newcommand{\res}{\operatorname{res}\nolimits}
\newcommand{\np}{\operatorname{np}\nolimits}
\newcommand{\Fac}{\operatorname{Fac}\nolimits}
\newcommand{\Aut}{\operatorname{Aut}\nolimits}
\newcommand{\DTr}{\operatorname{DTr}\nolimits}
\newcommand{\TrD}{\operatorname{TrD}\nolimits}
\newcommand{\Gr}{\operatorname{Gr}\nolimits}
\newcommand{\oH}{\ov{\BH}}
\newcommand{\FGS}{\operatorname{FGS}\nolimits}
\newcommand{\Mod}{\operatorname{Mod}\nolimits}
\newcommand{\R}{\operatorname{R}\nolimits}
\newcommand{\End}{\operatorname{End}\nolimits}
\newcommand{\lf}{\operatorname{l.f.}\nolimits}
\newcommand{\Iso}{\operatorname{Iso}\nolimits}
\newcommand{\aut}{\operatorname{Aut}\nolimits}
\newcommand{\Ui}{{\mathbf U}^\imath}
\newcommand{\UU}{{\mathbf U}\otimes {\mathbf U}}
\newcommand{\UUi}{(\UU)^\imath}
\newcommand{\tUU}{{\tU}\otimes {\tU}}
\newcommand{\tUUi}{(\tUU)^\imath}
\newcommand{\tUi}{\widetilde{{\mathbf U}}^\imath}
\newcommand{\sqq}{{\bf v}}
\newcommand{\sqvs}{\sqrt{\vs}}
\newcommand{\dbl}{\operatorname{dbl}\nolimits}
\newcommand{\swa}{\operatorname{swap}\nolimits}
\newcommand{\Gp}{\operatorname{Gp}\nolimits}

\newcommand{\tCMH}{{\cc\widetilde{\ch}(\K Q,\btau)}}

\newcommand{\tMHk}{{\widetilde{\ch}(\K Q,\btau)}}
\newcommand{\tCMHg}{\cc\widetilde{\ch}(Q,\btau)}

\newcommand{\iLa}{\Lambda^{\imath}}
\newcommand{\U}{{\mathbf U}}
\newcommand{\tU}{\widetilde{\mathbf U}}
\newcommand{\bvs}{\mathbf{\varsigma}}

\newcommand{\vs}{\varsigma}
\newcommand{\ov}{\overline}
\newcommand{\tk}{\widetilde{k}}
\newcommand{\tK}{\widetilde{K}}

\newcommand{\tH}{\widetilde{\ch}}

\newcommand{\utM}{\operatorname{\cs\cd\ch}\nolimits}
\newcommand{\tM}{\operatorname{\cs\cd\widetilde{\ch}}\nolimits}
\newcommand{\rM}{\operatorname{\cs\cd\ch_{\rm{red}}}\nolimits}
\newcommand{\utMH}{\operatorname{\cs\cd\ch(\iLa)}\nolimits}
\newcommand{\tMH}{{{\widetilde{\ch}(\K Q,\btau)}}}
\newcommand{\rMH}{\widetilde{\ch}_{\rm{red}}(kQ,\btau)}
\newcommand{\utMHg}{\operatorname{\ch(Q,\btau)}\nolimits}
\newcommand{\tMHg}{\operatorname{\widetilde{\ch}(Q,\btau)}\nolimits}
\newcommand{\rMHg}{\operatorname{\ch_{\rm{red}}(Q,\btau)}\nolimits}

\newcommand{\rMHd}{\widetilde{\ch}_{\rm{red}}(kQ,\btau)_{\bvsd}}
\newcommand {\rMHdi} {\widetilde{\ch}_{\rm{red}}(k\bs_iQ,\btau)_{\bvsd}}
%\newcommand{\tMHd}{\operatorname{\cs\cd\widetilde{\ch}(\iLa)_{\bvsd}}\nolimits}

%LU
%\newcommand{\utMHl}{\cs\cd\ch({\bs}_\ell\iLa)_{\diamond}}
\newcommand{\tMHl}{\widetilde{\ch}(\K \bs_\ell Q,\btau) }
\newcommand{\rMHl}{\cs\cd\ch_{\rm{red}}({\bs}_\ell\iLa)_{\bvsd}}
\newcommand{\tf}{\widetilde{f}}
\newcommand{\tMHi}{ \widetilde{\ch}(\K \bs_i Q,\btau)}
\newcommand{\rMHi}{\cs\cd\ch_{\rm{red}}({\bs}_i\iLa)_{\bvsd}}
\newcommand{\tMHgl}{\widetilde{\ch}({\bs}_\ell Q,\btau)}
\newcommand{\tMHgi}{\widetilde{\ch}({\bs}_i Q,\btau)}

\newcommand{\utGpg}{\operatorname{\ch^{\rm Gp}(Q,\btau)}\nolimits}
\newcommand{\tGpg}{\operatorname{\widetilde{\ch}^{\rm Gp}(Q,\btau)}\nolimits}
\newcommand{\rGpg}{\operatorname{\ch_{red}^{\rm Gp}(Q,\btau)}\nolimits}

\newcommand{\colim}{\operatorname{colim}\nolimits}
\newcommand{\gldim}{\operatorname{gl.dim}\nolimits}
\newcommand{\cone}{\operatorname{cone}\nolimits}
\newcommand{\rep}{\operatorname{rep}\nolimits}
\newcommand{\Ext}{\operatorname{Ext}\nolimits}
\newcommand{\Tor}{\operatorname{Tor}\nolimits}
\newcommand{\Hom}{\operatorname{Hom}\nolimits}
\newcommand{\Top}{\operatorname{top}\nolimits}
\newcommand{\Coker}{\operatorname{Coker}\nolimits}
\newcommand{\thick}{\operatorname{thick}\nolimits}
\newcommand{\rank}{\operatorname{rank}\nolimits}
\newcommand{\Gproj}{\operatorname{Gproj}\nolimits}
\newcommand{\Len}{\operatorname{Length}\nolimits}
\newcommand{\RHom}{\operatorname{RHom}\nolimits}
\renewcommand{\deg}{\operatorname{deg}\nolimits}
\renewcommand{\Im}{\operatorname{Im}\nolimits}
\newcommand{\Ker}{\operatorname{Ker}\nolimits}
\newcommand{\Coh}{\operatorname{Coh}\nolimits}
\newcommand{\Id}{\operatorname{Id}\nolimits}
\newcommand{\Qcoh}{\operatorname{Qch}\nolimits}
\newcommand{\CM}{\operatorname{CM}\nolimits}
\newcommand{\sgn}{\operatorname{sgn}\nolimits}
\newcommand{\Gdim}{\operatorname{G.dim}\nolimits}
\newcommand{\fpr}{\operatorname{\mathcal{P}^{\leq1}}\nolimits}
\newcommand{\ff}{B}

\newcommand{\For}{\operatorname{{\bf F}or}\nolimits}
\newcommand{\coker}{\operatorname{Coker}\nolimits}
\renewcommand{\dim}{\operatorname{dim}\nolimits}
\newcommand{\rankv}{\operatorname{\underline{rank}}\nolimits}
\newcommand{\dimv}{{\operatorname{\underline{dim}}\nolimits}}
\newcommand{\diag}{{\operatorname{diag}\nolimits}}
\newcommand{\qbinom}[2]{\begin{bmatrix} #1\\#2 \end{bmatrix} }

\renewcommand{\Vec}{{\operatorname{Vec}\nolimits}}
\newcommand{\pd}{\operatorname{proj.dim}\nolimits}
\newcommand{\gr}{\operatorname{gr}\nolimits}
\newcommand{\id}{\operatorname{Id}\nolimits}
\newcommand{\Res}{\operatorname{Res}\nolimits}

\newcommand{\pdim}{\operatorname{proj.dim}\nolimits}
\newcommand{\idim}{\operatorname{inj.dim}\nolimits}
\newcommand{\Gd}{\operatorname{G.dim}\nolimits}
\newcommand{\Ind}{\operatorname{Ind}\nolimits}
\newcommand{\add}{\operatorname{add}\nolimits}
\newcommand{\pr}{\operatorname{pr}\nolimits}
\newcommand{\oR}{\operatorname{R}\nolimits}
\newcommand{\oL}{\operatorname{L}\nolimits}
\newcommand{\ext}{{ \mathfrak{Ext}}}
\newcommand{\Perf}{{\mathfrak Perf}}
\def\scrP{\mathscr{P}}
\newcommand{\bk}{{\mathbb K}}
\newcommand{\cc}{{\mathcal C}}
\newcommand{\gc}{{\mathcal GC}}
\newcommand{\dg}{{\rm dg}}
\newcommand{\ce}{{\mathcal E}}
\newcommand{\cs}{{\mathcal S}}
\newcommand{\cl}{{\mathcal L}}
\newcommand{\cf}{{\mathcal F}}
\newcommand{\cx}{{\mathcal X}}
\newcommand{\cy}{{\mathcal Y}}
\newcommand{\ct}{{\mathcal T}}
\newcommand{\cu}{{\mathcal U}}
\newcommand{\cv}{{\mathcal V}}
\newcommand{\cn}{{\mathcal N}}
\newcommand{\mcr}{{\mathcal R}}
\newcommand{\ch}{{\mathcal H}}
\newcommand{\ca}{{\mathcal A}}
\newcommand{\cb}{{\mathcal B}}
\newcommand{\ci}{{\I}_{\btau}}
\newcommand{\cj}{{\mathcal J}}
\newcommand{\cm}{{\mathcal M}}
\newcommand{\cp}{{\mathcal P}}
\newcommand{\cg}{{\mathcal G}}
\newcommand{\cw}{{\mathcal W}}
\newcommand{\co}{{\mathcal O}}
\newcommand{\cq}{{Q^{\rm dbl}}}
\newcommand{\cd}{{\mathcal D}}
\newcommand{\cz}{{\mathcal Z}}
\newcommand{\ck}{\widetilde{\mathcal K}}
\newcommand{\calr}{{\mathcal R}}
\newcommand{\La}{\Lambda}
\newcommand{\LL}{\texttt{L}}
\newcommand{\RR}{\texttt{R}}
\newcommand{\ol}{\overline}
\newcommand{\st}{[1]}
\newcommand{\ow}{\widetilde}
\renewcommand{\P}{\mathbf{P}}
\newcommand{\pic}{\operatorname{Pic}\nolimits}
\newcommand{\Spec}{\operatorname{Spec}\nolimits}

\newtheorem{theorem}{Theorem}[section]
\newtheorem{acknowledgement}[theorem]{Acknowledgement}
\newtheorem{algorithm}[theorem]{Algorithm}
\newtheorem{axiom}[theorem]{Axiom}
\newtheorem{case}[theorem]{Case}
\newtheorem{claim}[theorem]{Claim}
\newtheorem{conclusion}[theorem]{Conclusion}
\newtheorem{condition}[theorem]{Condition}
\newtheorem{conjecture}[theorem]{Conjecture}
\newtheorem{construction}[theorem]{Construction}
\newtheorem{corollary}[theorem]{Corollary}
\newtheorem{criterion}[theorem]{Criterion}
\newtheorem{definition}[theorem]{Definition}
\newtheorem{example}[theorem]{Example}
\newtheorem{exercise}[theorem]{Exercise}
\newtheorem{lemma}[theorem]{Lemma}
\newtheorem{notation}[theorem]{Notation}
\newtheorem{problem}[theorem]{Problem}
\newtheorem{proposition}[theorem]{Proposition}
\newtheorem{solution}[theorem]{Solution}
\newtheorem{summary}[theorem]{Summary}
\numberwithin{equation}{section}

\theoremstyle{remark}
\newtheorem{remark}[theorem]{Remark}
\newcommand{\Pd}{\pi_*}
\def \bvs{{\boldsymbol{\varsigma}}}
\def \bvsd{{\boldsymbol{\varsigma}_{\diamond}}}
\def \btau{{{\tau}}}

\def \bp{{\mathbf p}}
\def \bq{{\bm q}}
\def \bv{{v}}
\def \bs{{\bm s}}
\newcommand{\bsigma}{\bm{\sigma}}

\newcommand{\bfv}{\mathbf{v}}
\def \brW{{\rm Br}(W^{\btau})}

\def \bfK{{\mathbf K}}
\def \bA{{\mathbf A}}
\def \ba{{\mathbf a}}
\def \bb{{\mathbf b}}
\def \bL{{\mathbf L}}
\def \bF{{\mathbf F}}
\def \bS{{\mathbf S}}
\def \bC{{\mathbf C}}
\def \bU{{\mathbf U}}
\def \bc{{\mathbf c}}
\def \fpi{\mathfrak{P}^\imath}
\def \Ni{N_\imath}
\def \fp{\mathfrak{P}}
\def \fg{\mathfrak{g}}
\def \fb{\mathfrak{b}}
\def \fk{\fg^\theta}  %\mathfrak{k}}
\def \p{p}
\def \fn{\mathfrak{n}}
\def \fh{\mathfrak{h}}
\def \fu{\mathfrak{u}}
\def \fv{\mathfrak{v}}
\def \fa{\mathfrak{a}}
\def \Z{{\Bbb Z}}
\def \F{{\Bbb F}}
\def \D{{\Bbb D}}
\def \C{{\Bbb C}}
\def \N{{\Bbb N}}
\def \Q{{\Bbb Q}}
\def \G{{\Bbb G}}
\def \P{{\Bbb P}}
\def \K{{k}}
\def \E{{\Bbb K}}
\def \A{{\Bbb A}}
\def \L{{\Bbb L}}
\def \I{{\Bbb I}}
\def \BH{{\Bbb H}}
\def \T{{\mathcal T}}
\def \tT{\widetilde{\mathcal T}}
\def \tTL{\tT(\iLa)}
\newcommand{\TT}{\bold{T}}

\newcommand {\lu}[1]{\textcolor{red}{$\clubsuit$: #1}}

\newcommand{\browntext}[1]{\textcolor{brown}{#1}}
\newcommand{\greentext}[1]{\textcolor{green}{#1}}
\newcommand{\redtext}[1]{\textcolor{red}{#1}}
\newcommand{\bluetext}[1]{\textcolor{blue}{#1}}
\newcommand{\brown}[1]{\browntext{ #1}}
\newcommand{\green}[1]{\greentext{ #1}}
\newcommand{\red}[1]{\redtext{ #1}}
\newcommand{\blue}[1]{\bluetext{ #1}}

%todo
\newcommand{\wtodo}{\todo[inline,color=orange!20, caption={}]}
\newcommand{\lutodo}{\todo[inline,color=green!20, caption={}]}

%%%%%
\title[Braid group symmetries on quasi-split $\imath$quantum groups]{Braid group symmetries on quasi-split $\imath$quantum groups  via $\imath$Hall algebras}

\author[Ming Lu]{Ming Lu}
\address{Department of Mathematics, Sichuan University, Chengdu 610064, P.R.China}
\email{luming@scu.edu.cn}

\author[Weiqiang Wang]{Weiqiang Wang}
\address{Department of Mathematics, University of Virginia, Charlottesville, VA 22904, USA}
\email{ww9c@virginia.edu}


\begin{abstract}
We establish automorphisms with closed formulas on quasi-split $\imath$quantum groups of symmetric Kac-Moody type associated to restricted Weyl groups. The proofs are carried out in the framework of $\imath$Hall algebras and reflection functors, thanks to the $\imath$Hall algebra realization of $\imath$quantum groups in our previous work. Several quantum binomial identities arising along the way are established.
\end{abstract}

\maketitle
\setcounter{tocdepth}{1}
 \tableofcontents

%%%%%%%
%%%%%%%
\section{Introduction}

\subsection{Background}

Among the fundamental structures of Drinfeld-Jimbo quantum groups is the existence of braid group symmetries \cite{Lus90, Lus90b}; also see \cite{KR90, LS90} for different formulations. The formulas for the actions of these automorphisms are intimately related to Lusztig's higher order Serre relations \cite{Lus93}. Reflection functors on Hall algebras can be used to construct braid group symmetries for quantum groups; see \cite{Rin96, X97, SV99, XY01}. Braid group actions have played a fundamental role in the constructions of PBW bases and canonical bases, and they also have applications in geometric representation theory and categorification.

A Satake diagram
\[
(\I =\I_\circ \cup \I_\bullet, \; \tau)
\]
 consists of a bicolored partition of the Dynkin diagram $\I =\I_\circ \cup \I_\bullet$ and a (possibly trivial) diagram involution $\tau$ subject to some compatibility conditions. Associated to a Satake diagram, a quantum symmetric pair $(\U, \Ui)$ \cite{Let99, Let02} consists of a Drinfeld-Jimbo quantum group $\U$ and its coideal subalgebra $\Ui$; we shall refer to $\Ui$ as an $\imath$quantum group and further call $\Ui$ quasi-split if $\I_\bullet =\emptyset$. On the other hand, a universal $\imath$quantum group $\tUi$ \cite{LW19} is a coideal subalgebra of the Drinfeld double quantum group $\tU$, and the $\imath$quantum group $\Ui$ with parameters \`a la Letzter is recovered by a central reduction of $\tUi$.

For (mostly) quasi-split $\imath$quantum groups of finite type with distinguished parameters, Kolb and Pellegrini \cite{KP11} constructed automorphisms $\TT_i$ of $\Ui$ for $i \in \I_\circ$ and show they satisfy the braid group relations associated to the restricted Weyl group of the symmetric pair; the formulas and the proofs therein relied essentially on computer computations. In type AI, the braid group action of $\Ui$ was noted earlier independently in \cite{Ch07, MR08}; see \cite{Dob20} for a recent progress. It was shown in \cite{BW18b} that Lusztig's braid group action $T_i$, for $i \in \I_\bullet$, preserves $\Ui$ of arbitrary Kac-Moody type. A natural and challenging question since then has been to establish these symmetries with closed formulas for their actions on generators of $\Ui$ in a conceptual way and in great generalities (such as Kac-Moody type and/or beyond quasi-split type).

The $\imath$Program \cite{BW18} aims at generalizing fundamental (algebraic, geometric, and categorical) constructions for quantum groups to $\imath$quantum groups. In case of (quasi-split) quantum symmetric pairs of diagonal type, we recover constructions for quantum groups.
In the framework of semi-derived Hall algebras \cite{Gor18, LP21, Lu19} (generalizing \cite{Rin90, Gr95, Br13} and \cite{T06, XX08}), the authors have developed an $\imath$Hall algebra realization for the quasi-split universal $\imath$quantum groups $\tUi$, first for finite type in \cite{LW19} and then for symmetric Kac-Moody type in \cite{LW20a}. In \cite{LW21}, we have used reflection functors to provide a conceptual realization of the braid group symmetries on a class of quasi-split universal $\imath$quantum groups $\tUi$ of {\em finite type}, and subsequently on the corresponding $\Ui$, which agree with \cite{KP11} for distinguished parameters.

The reflection functors are formulated in the Hall basis and have the advantage that the resulting maps on $\tUi$ are automatically algebra automorphisms; however, the reflection functor approach provides very little clue on the closed formulas for the corresponding automorphisms on ($\imath$-)quantum groups in terms of Chevalley generators. Therefore, a complementary approach is needed.

In \cite{CLW18} joint with X.~Chen, the authors obtained a Serre presentation for quasi-split $\imath$quantum groups $\Ui$ of arbitrary Kac-Moody type (cf. \cite{Ko14}), building on earlier work of Letzter \cite{Let02} and others. Similarly, a universal $\imath$quantum group $\tUi$ admits a Serre presentation (see Proposition~\ref{prop:Serre}) with Chevalley generators $B_i$ and $\tk_i$, for $i\in \I$ \cite{LW20a}:
\begin{align}
   \label{eq:Che}
\tUi = \big\langle B_i,\tk_i \; (i\in \I) \mid \text{relations } \eqref{relation1}\text{--}\eqref{relation6} \big\rangle.
\end{align}
Chen and the authors established in \cite{CLW21} the Serre-Lusztig (i.e., higher order Serre) relations between $B_i, B_j$ in $\tUi$ for $i=\tau i$, generalizing Lusztig's  construction in \cite[Chapter 7]{Lus93}. Based on the expectation that the close relations between braid group actions and Serre-Lusztig relations in quantum groups carry over to the setting of $\imath$quantum groups, we made a conjecture in \cite[Conjecture~ 6.5]{CLW21} on closed formulas for automorphisms $\TT'_{i,e}, \TT''_{i,e}$, for $i =\tau i$ and $e\in \{\pm 1\}$.  Recently, we established in \cite{CLW21c} the Serre-Lusztig relations between $B_i, B_j$ in $\tUi$ for $i \neq \tau i$; in addition, we made a conjecture in \cite[Conjecture~3.7]{CLW21c} on closed formulas for automorphisms $\TT'_{i,e}, \TT''_{i,e}$, for $i \neq \tau i$, though the relations with Serre-Lusztig relations were not as direct as in the earlier cases.

\subsection{The goal}

Recall Lusztig has constructed 4 variants of braid group symmetries denoted by $T_{i,e}'$ and $T_{i,e}''$, for $i\in \I$ and $e \in \{\pm 1\}$, on the quantum group $\U$ \cite[Chapter~ 37]{Lus93}. One can further transform from one variant to another by twisting with various well-known involutions and anti-involutions.

In this paper, we work with quasi-split $\imath$quantum groups $\tUi$ and $\Ui$ of arbitrary symmetric Kac-Moody type, where a mild condition that the Cartan integers $c_{j, \tau j}$ for all $j \in \I$ are even is further imposed due to a use of the $\imath$Hall algebra technique. We fix a set $\I_\tau$ of representatives of $\tau$-orbits on $\I$, and let $\ov{\I}_\btau$ denote the subset of $i\in \I_\tau$ whose $\tau$-orbit is of finite type.
By definition, the restricted Weyl group $W^\tau$ is the $\tau$-fixed point subgroup of $W$. According to \cite[Appendix]{Lus03}, $W^\tau$ is a Coxeter group generated by ${\bf s}_i$ defined in \eqref{def:si}, for $i \in \ov{\I}_\btau$.

We shall establish 4 versions of automorphisms on $\tUi$, 
\[
\TT_{i,e}',\; \TT_{i,e}'', 
\quad  \text{ for   }i\in \ov{\I}_\btau, \; e \in \{\pm 1\}
\]
with closed formulas on the Chevalley generators in \eqref{eq:Che}; the  $\TT'_{i,e}$ and $\TT''_{i,e}$ are related to each other via a bar involution $\psi_\imath$ and an anti-involution $\sigma_\imath$ on $\tUi$. Our results have confirmed substantially in the setting of quasi-split $\imath$quantum groups \cite[Conjecture~ 6.5]{CLW21} in case $i=\tau i$ and  \cite[Conjecture~3.7]{CLW21c} in case $i\not =\tau i$.

\subsection{The main results}

The formulas for the automorphisms $\TT''_{i,e}$ are given in terms of the $\imath$divided powers $B_{i,\ov{p}}^{(r)}$ \eqref{eq:iDPodd}--\eqref{eq:iDPev} (for $i =\tau i$) and standard divided powers \eqref{eq:DP} $B_{i}^{(r)}$ (for $i\neq \tau i$). The $\imath$divided powers arose from the theory of $\imath$canonical basis and they depend on a parity $\ov p \in \Z_2$ (cf. \cite{BW18, BeW18, CLW18}).

\begin{customthm} {\bf A}
[Theorems~\ref{thm:BG}, \ref{thm:BG1}]
   \label{thm:A}
For $i\in \bar{\I}_\btau$ and $e \in \{\pm 1\}$, there are automorphisms $\TT''_{i,e}$ on $\tUi$ such that
\begin{enumerate}
\item
$\underline{(i=\btau i)}:$ \;
$\TT''_{i,e}(\tk_j)= (-v^{1+e} \tk_i)^{-c_{ij}}\tk_j$, and
\begin{align*}
\TT''_{i,e}(B_i) &= (-v^{1+e}\tk_{i})^{-1}B_i,
\\
\TT''_{i,e}(B_j) &= \sum_{r+s=-c_{ij}}(-1)^{r}v^{er}B_{i,\ov{p}}^{(r)}B_jB_{i,\ov{c_{ij}}+\ov{p}}^{(s)}\\
&+\sum_{u\geq1}\sum_{\stackrel{r+s+2u=-c_{ij}}{
\ov{r}=\ov{p}}}(-1)^{r+u}v^{er+eu}B_{i,\ov{p}}^{(r)}B_jB_{i,\ov{c_{ij}}+\ov{p}}^{(s)}(v\tk_i)^u, \quad \text{ for }j\neq i;
\end{align*}
\item
$\underline{(i\neq \btau i)}:$ \;
$\TT''_{i,e}(\tk_j)= \tk_i^{-c_{ij}} \tk_{\btau i}^{-c_{\tau i,j}} \tk_j$,
\begin{align*}
\TT''_{i,1}(B_j) &=
\begin{cases}  -\tk_{i}^{-1}B_{\btau i},  & \text{ if }j=i \\
-B_i\tk_{\btau i}^{-1}  ,  &\text{ if }j=\btau i,\end{cases}
\qquad\qquad
\TT''_{i,-1}(B_j)= \begin{cases}  -\tk_{\btau i}^{-1}B_{\btau i},  & \text{ if }j=i \\
-B_i\tk_{i}^{-1},  &\text{ if }j=\btau i,\end{cases}
\end{align*}
and for $j\neq i,\btau i$,
\begin{align*}
\TT''_{i,1}(B_j)
&= \sum^{-\max(c_{ij},c_{\tau i,j})}_{u=0} \; \sum^{-c_{ i,j}-u}_{r=0} \; \sum_{s=0}^{-c_{\tau i,j}-u} (-1)^{r+s} v^{ r-s+(-c_{ij}-r-s-u)u } \\
&\qquad\qquad\qquad\qquad\qquad\qquad
 \times B_i^{(r)} B_{\tau i}^{(-c_{\tau i,j}-u-s)} B_j B_{\tau i}^{(s)} B_i^{(-c_{ij}-r-u)}\tk_{\tau i}^u,
\\
\TT''_{i,-1}(B_j)
= &\sum^{-\max(c_{ij},c_{\tau i,j})}_{u=0} \; \sum^{-c_{ i,j}-u}_{r=0} \; \sum_{s=0}^{-c_{\tau i,j}-u} (-1)^{r+s} v^{- (r-s+(-c_{ij}-r-s-u)u )} \\
&\qquad\qquad\qquad\qquad\qquad\qquad
 \times B_i^{(r)} B_{\tau i}^{(-c_{\tau i,j}-u-s)} B_j B_{\tau i}^{(s)} B_i^{(-c_{ij}-r-u)}\tk_{i}^u.
\end{align*}
\end{enumerate}
\end{customthm}

The closed formulas for the other automorphisms $\TT'_{i,e}$ of $\tUi$ can be found in Theorem~\ref{thm:BG2}.

\begin{customthm}{\bf B} [Theorem~\ref{thm:inverse}]
  \label{thm:B}
We have $\TT_{i,e}' =(\TT_{i,-e}'')^{-1}$, for any $i\in\bar{\I}_\btau$ and $e=\{\pm1\}$.
\end{customthm}

Note that the leading terms (i.e., the $u=0$ summands) in the formulas for the symmetries $\TT'_{i,e}, \TT''_{i,e}$ acting on $B_j$ (for $j\neq i, j\neq \tau i$) in $\tUi$ are precisely the formulas for the symmetries $T'_{i,e}, T''_{i,e}$ on $F_j$ (and also $E_j$ in case $i \neq \tau i$) in $\U$; see \cite[37.1.3]{Lus93}.

Theorem~\ref{thm:A} in the special case of universal affine $\imath$quantum group of split rank one (i.e, universal q-Onsager algebra) has already found applications in \cite{LW20b}; compare with the braid group action on the q-Onsager algebra in \cite{BK20}. For $\imath$quivers of diagonal type, the $\imath$quantum group is the Drinfeld double $\tU$, and 
we reformulate Theorem~\ref{thm:A} as Propositions~\ref{prop:BG1U}--\ref{prop:BG2U} on the automorphisms on $\tU$; they descend to automorphisms on $\U$ upon the identification $\tK_i'=\tK_i^{-1}$, providing a new $\imath$Hall algebra approach for some main formulas in  \cite[Chapter 37]{Lus93}; compare \cite{SV99, XY01}.

The automorphisms $\TT'_{i,e}, \TT''_{i,e}$ on $\tUi$ descend to automorphisms on an $\imath$quantum group $\Ui$ associated with distinguished parameters.

It remains a fundamental problem to formulate and establish conceptually in full generality the braid group actions on $\imath$quantum groups of arbitrary (not necessarily quasi-split) finite type and on their modules. This shall require completely different new ideas.

In the process of proving Theorem~\ref{thm:A}, we are led to the following $v$-binomial identities which are of independent interest; for an additional novel $v$-binomial identity see Proposition~\ref{prop:B4}. We refer to \eqref{eq:binom} for notation $[2k]^{!!}$.

\begin{customthm} {\bf C}
[Propositions~\ref{prop:DD}--\ref{prop:DD=CC}, Theorem~\ref{thm:DDCC}]
  \label{thm:C}
The following identities hold, for $d\ge 1$:
\begin{align*}
\sum_{\stackrel{k,m,n \in \N}{k+m+n =d} }  \varepsilon^n
 \frac{v^{\binom{n+1}{2} -2km  +2k}}{[n]^! [2k]^{!!} [2m]^{!!}}
&=
\sum_{\stackrel{t,k,m,n \in \N}{t+k+m+n =d} }  \varepsilon^n
 \frac{v^{t^2 -2dt +t +2nt +\binom{n+1}{2} -2km -2m}}{[n]^! [2k]^{!!} [2m]^{!!}} (v -v^{-1})^{t}
\\
&=
\begin{cases}
0, & \text{ if }  \varepsilon =-1,
\\
\frac{2v^d (v +v^{-1}) (v^2 +v^{-2}) \cdots (v^{d-1} +v^{1-d})}{[d]!},  & \text{ if }  \varepsilon =1.
\end{cases}
\end{align*}
\end{customthm}
The $v$-binomial identities in Theorem~\ref{thm:C} are closed related to some $v$-binomial identities arising from the $\imath$Hall algebra realization of $\imath$quantum groups of Kac-Moody type in \cite{LW20a}.

\subsection{Our approach}

Lusztig formulated and established the braid group actions on both the algebra $\U$ and on its integrable highest weight modules, with the actions on the module level established first; cf. \cite[\S5.2, Chapter 37]{Lus93}. In the $\imath$quantum group setting, we do not have the braid group action on the module level available for now, and so we cannot follow Lusztig's approach. Instead, we shall resort to $\imath$Hall algebras and reflection functors, and accordingly, most constructions and computations in this paper will be carried out in the setting of $\imath$Hall algebras.

Our strategy is to verify the closed formulas for the action of reflection functors $\Gamma_i$, for $i \in \ov{\I}_\btau$. Via the algebra isomorphism from $\tUi$ to the composition $\imath$Hall algebra (see Theorem~\ref{thm:Ui=iHall} and Corollary~\ref{cor:isomgeneric}), we transfer the isomorphism $\Gamma_i$ of $\imath$Hall algebras to an automorphism $\TT''_{i,1}$ of the $\imath$quantum group $\tUi$. The detailed study of reflection functors can be of interest in its own. 

The formulas for $\TT''_{i,e}$ and $\TT'_{i,e}$ are very different and so are the proofs, depending on whether or not $i =\tau i$.

Let us explain in some detail the proof for the formula for $\TT''_{i,1}(B_j)$ in Theorem~\ref{thm:A}(1) for $i = \tau i$. We shall aim at establishing the closed formulas for the reflection functor $\Gamma_i$ in the $\imath$Hall algebra corresponding to $\TT''_{i,1}(B_j)$; see Theorem~\ref{thm:Braidsplit}. To that end, we shall
compute the $\imath$Hall algebra version of the RHS of the formula for $\TT''_{i,1}(B_j)$ in Theorem~\ref{thm:A} in the $\imath$Hall basis, starting with an expansion formula of the $\imath$divided powers in the $\imath$Hall basis established in \cite{LW20a}. This (mostly) homological algebra computation in a rank 2 $\imath$Hall algebra is challenging and tedious. By a comparison with the definition of the reflection functor $\Gamma_i$, we then need to show that all except the leading summand in the $\imath$Hall basis vanishes. After several steps of combinatorial reduction, we reduce the desired vanishing properties to a combinatorial identity, which is then derived from the $v$-binomial identities in Theorem~\ref{thm:C}.

While the strategy for the proof of Theorem~\ref{thm:A}(2) in case $i\neq \tau i$ is similar, the details are all different and separated from the case $i=\tau i$; we shall prove its $\imath$Hall algebra counterpart (see Theorem~ \ref{thm:braidqsplit}), and this requires new homological algebra computations and combinatorial reductions in  Section~\ref{sec:formulaRF2}. This explains in part the length of the paper. 

Recall the counterpart of Theorem~\ref{thm:B} holds for the corresponding braid group symmetries on quantum groups, and the proof in \cite[Chapter 37]{Lus93} relies essentially on the braid group actions on the module level. We present a conceptual simple proof of Theorem~\ref{thm:B} in the framework of $\imath$Hall algebras.

\subsection{The organization}

This paper is organized as follows.
In Section~\ref{sec:prelim}, we review the basics on $\imath$Hall algebras, modulated graphs for $\imath$quivers, and reflection functors.

In Section~\ref{sec:formulaRF}, we establish the reflection functor counterpart of Theorem~\ref{thm:A} for $i=\tau i$, by reducing it to a $v$-binomial identity. This novel identity is then derived from Theorem~\ref{thm:C} which is established in Section~\ref{sec:identities}.
%To be removed after arxiv v1
We also relegate to Appendix~\ref{app:braidsplitodd} the detail for a proof of a second half of $\imath$Serre relations (which is very similar to the one of the first half of $\imath$Serre relations given in this section).

In Section~\ref{sec:formulaRF2}, we establish the reflection functor counterpart of Theorem~\ref{thm:A} for $i \neq \tau i$. The long proof of Proposition~ \ref{prop:build-block} is given in Appendix~\ref{app:prop}.

In Section~\ref{sec:braid}, we review $\imath$quantum groups and their $\imath$Hall algebra realizations. Then we reformulate in terms of $\tUi$ the formulas for the reflection functors on $\imath$Hall algebras obtained in the previous sections. We also prove Theorem~\ref{thm:B}.

\vspace{2mm}
{\bf Acknowledgement.}
We thank Xinhong Chen for her collaboration in related projects and  the joint formulation of the conjectural formulas proved in this paper. ML is partially supported by the National Natural Science Foundation of China (No. 12171333). WW is partially supported by the NSF grant DMS-2001351.

%%%%%%%
%%%%%%%
\section{$\imath$Quivers and $\imath$Hall algebras}
\label{sec:prelim}

In this section, we review and set up notations for $\imath$quivers, $\imath$Hall algebras, and the reflection functors on $\imath$quiver algebras, following \cite{LW19, LW20a, LW21} (also cf. \cite{LP21, Lu19}).

\subsection{Notations}
\label{subsec: notat}
For an additive category $\ca$ and $M\in \ca$, we denote

$\triangleright$ $\add M$ -- subcategory of $\ca$ whose objects are the direct summands of finite direct sums of copies of $M$,

$\triangleright$ $\Fac M$-- the full subcategory of $\ca$ of epimorphic images of objects in $\add M$,

$\triangleright$ $\Iso(\ca)$ -- set of the isoclasses of objects in $\ca$,

$\triangleright$ $[M]$ -- the isoclass of $M$.
\medskip

For an exact category $\ca$ and $M\in\ca$, we denote

$\triangleright$ $K_0(\ca)$ -- the Grothendieck group of $\ca$,

$\triangleright$ $\widehat{M}$ -- the class of $M$ in $K_0(\ca)$.

\medskip

Let $\K$ be a field.
For a quiver algebra $A=\K Q/I$ (not necessarily finite-dimensional), we always identify left $A$-modules with representations of $Q$ satisfying relations in $I$. A representation $V=(V_i,V(\alpha))_{i\in Q_0,\alpha\in Q_1}$ of $A$ is called {\em nilpotent} if for each oriented cycle $\alpha_m\cdots\alpha_1$ at a vertex $i$, the $\K$-linear map $V(\alpha_m)\cdots V(\alpha_1):V_i\rightarrow V_i$ is nilpotent. We denote

$\triangleright$ $\mod(A)$ -- category of finite-dimensional nilpotent $A$-modules,

$\triangleright$ ${\rm proj.dim}_AM$ -- projective dimension of an $A$-module $M$,

$\triangleright$ ${\rm inj.dim}_AM$ -- injective dimension of $M$,

$\triangleright$ $D=\Hom_\K(-,\K)$ -- the standard duality.

\subsection{The $\imath$quiver algebras}

We recall the $\imath$quiver algebras from \cite[\S2]{LW20a}; see also \cite[\S2]{LW19}.

Let $\K$ be a field.
Let $Q=(Q_0,Q_1)$ be a quiver (not necessarily acyclic). Let $n_{ij}$ be the number of edges connecting vertices $i$ and $j$. Throughout the paper, we shall identify
\[
\I=Q_0.
\]
An {\em involution} of $Q$ is defined to be an automorphism $\btau$ of the quiver $Q$ such that $\btau^2=\Id$. In particular, we allow the {\em trivial} involution $\Id:Q\rightarrow Q$. An involution $\btau$ of $Q$ induces an involution of the path algebra $\K Q$, again denoted by $\btau$.
A quiver together with a specified involution $\btau$, $(Q, \btau)$, will be called an {\em $\imath$quiver}.

Let $R_1$ denote the truncated polynomial algebra $\K[\varepsilon]/(\varepsilon^2)$.
Let $R_2$ denote the radical square zero of the path algebra of $\xymatrix{1 \ar@<0.5ex>[r]^{\varepsilon} & 1' \ar@<0.5ex>[l]^{\varepsilon'}}$, i.e., $\varepsilon' \varepsilon =0 =\varepsilon\varepsilon '$. Define a $\K$-algebra
\begin{equation}
  \label{eq:La}
\Lambda=\K Q\otimes_\K R_2.
\end{equation}

Associated to the quiver $Q$, the {\em double framed quiver} $Q^\sharp$
is the quiver such that
\begin{itemize}
\item the vertex set of $Q^{\sharp}$ consists of 2 copies of the vertex set $Q_0$, $\{i,i'\mid i\in Q_0\}$;
\item the arrow set of $Q^{\sharp}$ is
\[
\{\alpha: i\rightarrow j,\alpha': i'\rightarrow j'\mid(\alpha:i\rightarrow j)\in Q_1\}\cup\{ \varepsilon_i: i\rightarrow i' ,\varepsilon'_i: i'\rightarrow i\mid i\in Q_0 \}.
\]
\end{itemize}
Note $Q^\sharp$ admits a natural involution, $\swa$.
The involution $\btau$ of a quiver $Q$ induces an involution ${\btau}^{\sharp}$ of $Q^{\sharp}$ which is basically the composition of $\swa$ and $\btau$ (on the two copies of subquivers $Q$ and $Q'$ of $Q^\sharp$). The algebra $\La$ can be realized in terms of the quiver $Q^{\sharp}$ and a certain ideal $I^{\sharp}$
so that $\Lambda\cong \K Q^{\sharp} \big/ I^{\sharp}$.

By definition, ${\btau}^{\sharp}$ on $Q^\sharp$ preserves $I^\sharp$ and hence induces an involution ${\btau}^{\sharp}$ on the algebra $\Lambda$. The {\rm $\imath$quiver algebra} of $(Q, \btau)$ is the fixed point subalgebra of $\Lambda$ under ${\btau}^{\sharp}$,
\begin{equation}
   \label{eq:iLa}
\iLa
= \{x\in \Lambda\mid {\btau}^{\sharp}(x) =x\}.
\end{equation}
The algebra $\iLa$ can be described in terms of a certain quiver $\ov Q$ and its ideal $\ov{I}$ so that $\iLa \cong \K \ov{Q} / \ov{I}$; see \cite[Proposition 2.6]{LW19}.
We recall $\ov{Q}$ and $\ov{I}$ as follows:
\begin{itemize}
\item[(i)] $\ov{Q}$ is constructed from $Q$ by adding a loop $\varepsilon_i$ at the vertex $i\in Q_0$ if $i=\btau i$, and adding an arrow $\varepsilon_i: i\rightarrow \btau i$ for each $i\in Q_0$ if $i\neq \btau i$;
\item[(ii)] $\ov{I}$ is generated by
\begin{itemize}
\item[(1)] (Nilpotent relations) $\varepsilon_{i}\varepsilon_{\btau i}$ for any $i\in\I$;
\item[(2)] (Commutative relations) $\varepsilon_i\alpha-\btau(\alpha)\varepsilon_j$ for any arrow $\alpha:j\rightarrow i$ in $Q_1$.
\end{itemize}
\end{itemize}
Moreover, it follows by \cite[Proposition 2.2]{LW20a} that $\Lambda^{\imath}$ is a $1$-Gorenstein algebra.

By \cite[Corollary 2.12]{LW19}, $\K Q$ is naturally a subalgebra and also a quotient algebra of $\Lambda^\imath$.
Viewing $\K Q$ as a subalgebra of $\Lambda^{\imath}$, we have a restriction functor
\[
\res: \mod (\Lambda^{\imath})\longrightarrow \mod (\K Q).
\]
Viewing $\K Q$ as a quotient algebra of $\Lambda^{\imath}$, we obtain a pullback functor
\begin{equation}\label{eqn:rigt adjoint}
\iota:\mod(\K Q)\longrightarrow\mod(\Lambda^{\imath}).
\end{equation}
Hence a simple module $S_i (i\in Q_0)$ of $k Q$ is naturally a simple $\iLa$-module.

For each $i\in Q_0$, define a $k$-algebra (which can be viewed as a subalgebra of $\iLa$)
\begin{align}\label{dfn:Hi}
\BH _i:=\left\{ \begin{array}{cc}  \K[\varepsilon_i]/(\varepsilon_i^2) & \text{ if }i=\btau i,
 \\
\K(\xymatrix{i \ar@<0.5ex>[r]^{\varepsilon_i} & \btau i \ar@<0.5ex>[l]^{\varepsilon_{\btau i}}})/( \varepsilon_i\varepsilon_{\btau i},\varepsilon_{\btau i}\varepsilon_i)  &\text{ if } \btau i \neq i .\end{array}\right.
\end{align}
Note that $\BH _i=\BH _{\btau i}$ for any $i\in Q_0$.
Choose one representative for each $\btau$-orbit on $\I$, and let
\begin{align}
\label{eq:ci}
\ci = \{ \text{the chosen representatives of $\btau$-orbits in $\I$} \}.
\end{align}

Define the following subalgebra of $\Lambda^{\imath}$:
\begin{equation}  \label{eq:H}
\BH =\bigoplus_{i\in \ci }\BH _i.
\end{equation}
Note that $\BH $ is a radical square zero selfinjective algebra. Denote by
\begin{align}
\res_\BH :\mod(\iLa)\longrightarrow \mod(\BH )
\end{align}
the natural restriction functor.
On the other hand, as $\BH $ is a quotient algebra of $\iLa$ (cf. \cite[proof of Proposition~ 2.15]{LW19}), every $\BH $-module can be viewed as a $\iLa$-module.

Recall the algebra $\BH _i$ for $i \in \ci$ from \eqref{dfn:Hi}. For $i\in Q_0 =\I$, define the indecomposable module over $\BH _i$ (if $i\in \ci$) or over $\BH_{\btau i}$ (if $i\not \in \ci$)
\begin{align}
  \label{eq:E}
\E_i =\begin{cases}
k[\varepsilon_i]/(\varepsilon_i^2), & \text{ if }i=\btau i;
\\
\xymatrix{\K\ar@<0.5ex>[r]^1 & \K\ar@<0.5ex>[l]^0} \text{ on the quiver } \xymatrix{i\ar@<0.5ex>[r]^{\varepsilon_i} & \btau i\ar@<0.5ex>[l]^{\varepsilon_{\btau i}} }, & \text{ if } i\neq \btau i.
\end{cases}
\end{align}
Then $\E_i$, for $i\in Q_0$, can be viewed as a $\iLa$-module and will be called a {\em generalized simple} $\iLa$-module.

Let $\cp^{<\infty}(\Lambda^\imath)$ be the subcategory of $\mod(\Lambda^\imath)$ formed by modules of finite projective dimensions.
Let $\cp^{\leq d}(\Lambda^\imath)$ be the subcategory of $\mod(\Lambda^\imath)$ which consists of $\Lambda^\imath$-modules of projective dimension less than or equal to $d$, for $d\in\N$.
Then $\cp^{<\infty}(\Lambda^\imath)=\cp^{\leq1}(\Lambda^\imath)$, and $\E_i\in\cp^{<\infty}(\Lambda^\imath)$ for any $i\in \I$; see \cite[Lemma 2.3]{LW20a}.

Following \cite{LW19}, we can define the Euler forms $\langle K,M\rangle =\langle K,M\rangle_{\Lambda^\imath}$ and $\langle M,K\rangle =\langle M,K\rangle_{\Lambda^\imath}$ for any $K\in\cp^{\leq1}(\Lambda^\imath)$, $M\in\mod(\Lambda^\imath)$. These forms descend to bilinear Euler forms on the Grothendieck groups:
\begin{eqnarray*}
\langle\cdot,\cdot\rangle: K_0(\cp^{\leq 1}(\Lambda^\imath))\times K_0(\mod(\Lambda^\imath))\longrightarrow \Z,
\\
\langle\cdot,\cdot\rangle: K_0(\mod(\Lambda^\imath))\times K_0(\cp^{\leq 1}(\Lambda^\imath))\longrightarrow \Z,
\end{eqnarray*}
such that
\begin{equation}
\langle \widehat{K},\widehat{M}\rangle=\langle K,M\rangle,\qquad\langle \widehat{M},\widehat{K}\rangle=\langle M,K\rangle,\qquad\forall\, K\in\cp^{\leq1}(\Lambda^\imath), M\in\mod(\Lambda^\imath).
\end{equation}

Denote by $\langle\cdot,\cdot\rangle_Q$ the Euler form of $\K Q$. Denote by $S_i$ the simple $\K Q$-module (respectively, $\Lambda^{\imath}$-module) corresponding to vertex $i\in Q_0$ (respectively, $i\in\ov{Q}_0$).
These 2 Euler forms are related via the restriction functor $\res:\mod(\Lambda^\imath)\rightarrow \mod(\K Q)$ as follows.

\begin{lemma}
[\text{\cite[Lemma 3.1]{LW20a}}]
   \label{lemma compatible of Euler form}
We have
\begin{enumerate}
\item $\langle K,M\rangle =\langle \res_{\BH}(K),M\rangle$, $\langle M,K\rangle =\langle M,\res_{\BH}(K)\rangle$, $\forall M\in \mod(\Lambda^\imath), K\in\cp^{\leq1}(\Lambda^\imath)$;
\item
$\langle \E_i, M\rangle = \langle S_i,\res (M) \rangle_Q$, $\langle M,\E_i\rangle =\langle \res(M), S_{\btau i} \rangle_Q$, $\forall i\in Q_0$, $M\in\mod(\Lambda^{\imath})$;
\item $\langle M,N\rangle=\frac{1}{2}\langle \res(M),\res(N)\rangle_Q$, $\forall M,N\in\cp^{\leq 1}(\Lambda^\imath)$.
\end{enumerate}
\end{lemma}

An oriented cycle $c$ of $Q$ is called \emph{minimal} if
$c$ does not contain any proper oriented cycle $c'$. For any minimal cycle of length $m$, we call it an \emph{$m$-cycle} for short. In particular, $1$-cycles are called {\em loops}.
An $\imath$quiver $(Q,\btau)$ is called {\em virtually acyclic} if it satisfies the conditions (A1)--(A2):
\begin{itemize}
\item[(A1)] $Q$ does not have any minimal $m$-cycles for $m\neq2$;
\item[(A2)] For any $i,j\in Q_0$, if $\btau i\neq j$, then the full minimal subquiver of $Q$ containing $i,j$ is acyclic.
\end{itemize}

In the remaining of the paper, {\em we always assume that $\imath$quivers $(Q,\btau)$ are virtually acyclic.}

\subsection{Modulated graphs for $\imath$quivers}
  \label{subsec:modulated}
In this subsection, we study representations of modulated graphs associated to virtually acyclic $\imath$quivers, slightly generalizing \cite[\S2.5]{LW19} to the generality of \cite{LW20a}. We then formulate the BGP type reflection functors for $\imath$quivers (not necessarily acyclic). See  \cite{DR74, Li12, GLS17} for more details about representations of modulated graphs.

\subsubsection{}

Let $(Q, \btau)$ be a virtually acyclic $\imath$quiver, and $\Lambda^\imath=\K\ov{Q}/\ov{I}$ with $(\ov{Q}, \ov{I})$ being defined in \cite[Proposition 2.6]{LW20a}.
For each $i\in Q_0$, define a $k$-algebra
\begin{align}\label{dfn:bHi}
{\ov{\BH}} _i:=\left\{ \begin{array}{ll}  \K[\varepsilon_i]/(\varepsilon_i^2) & \text{ if }i=\btau i,
 \\
\K(\xymatrix{ i  \ar@/^0.5pc/@<0.8ex>[r]^{\varepsilon_i} \ar@<0.75ex>[r]|-{r_i} & \btau i\ar@<0.75ex>[l]|-{r_i}  \ar@/^0.5pc/@<0.8ex>[l]^{\varepsilon_{\btau i}} })/( \varepsilon_i\varepsilon_{\btau i},\varepsilon_{\btau i}\varepsilon_i, \alpha_j\varepsilon_{\btau i}-\varepsilon_i\beta_{j},\beta_j\varepsilon_i-\varepsilon_{\btau i}\alpha_j\mid 1\leq j\leq r)  &\text{ if } \btau i \neq i ,\end{array}\right.
\end{align}
with $n_{i,\btau i}=2r_i$. 
Here and below, $\xymatrix{1\ar[r]|-m&2}$ means there are $m$ arrows from $1$ to $2$. Define the following subalgebra of $\Lambda^{\imath}$:
\begin{equation}  \label{eq:H}
	\ov{\BH} =\bigoplus_{i\in \ci }\ov{\BH} _i.
\end{equation}

Define
\begin{align}
\label{eqn:orientation}
\Omega:=\Omega(Q) =\{(i,j) \in Q_0\times Q_0\mid  \exists (\alpha:i\rightarrow j)\in Q_1, \btau i\neq j\}.
\end{align}
Then $\Omega$ represents the orientation of $Q$. Since $(Q,\btau)$ is virtually acyclic, if $(i,j)\in\Omega$, then $(j,i)\notin \Omega$.
We also use $\Omega(i,-)$ to denote the subset $\{j\in Q_0 \mid \exists (\alpha:i\rightarrow j)\in Q_1\}$, and $\Omega(-,i)$ is defined similarly.

For any $(i,j)\in \Omega$, we define
\begin{equation}  \label{eq:jHi}
{}{}_j{\ov{\BH}}_i := {\ov{\BH}} _j\Span_k\{\alpha,\btau\alpha\mid(\alpha:i\rightarrow j)\in Q_1\text{ or } (\alpha:i\rightarrow \btau j)\in Q_1\}{\ov{\BH}} _i.
\end{equation}
Note that ${}{}_j{\ov{\BH}}_i ={}_{\btau j} {\ov{\BH}} _{\btau i}={}_{j} {\ov{\BH}} _{\btau i}={}_{\btau j} {\ov{\BH}} _{i}$ for any $(i,j)\in \Omega$.

Hence ${}_j\oH_i $ is an $\oH _j\mbox{-}\oH _i$-bimodule, which is free as a left $\oH _j$-module (and respectively, right $\oH _i$-module), with a basis $_j\LL_i$ (and respectively, $_j\RR_i$) defined in the following; cf. \cite[(2.12), (2.13)]{LW19}.
\begin{eqnarray}
     \label{basis of Hij left}
_j\LL_i&=&\left\{ \begin{array}{cc}
\{\alpha \mid (\alpha:i\rightarrow j)\in Q_1\} & \text{ if }i=\btau i, \btau j=j,\\
\{\alpha+\btau\alpha \mid(\alpha:i\rightarrow j)\in Q_1\} & \text{ if }i=\btau i, \btau j\neq j,\\
\{\alpha,\btau \alpha \mid (\alpha:i\rightarrow j)\in Q_1\} & \text{ if }i\neq \btau i,\btau j=j,\\
\{\alpha+\btau \alpha \mid (\alpha:i\rightarrow j) \text{ or }(\alpha:i\rightarrow \btau j)\in Q_1\} & \text{ if }i\neq \btau i,\btau j\neq j;\label{eqn:basis of L}
\end{array}\right.\\
_j\RR_i&=& \left\{ \begin{array}{cc}\label{basis of Hij right}
\{\alpha \mid (\alpha:i\rightarrow j)\in Q_1\} & \text{ if }i=\btau i,  \btau j=j,\\
\{\alpha,\btau \alpha \mid (\alpha:i\rightarrow j)\in Q_1\} & \text{ if }i=\btau i, \btau j\neq j,\\
\{\alpha+\btau \alpha \mid (\alpha:i\rightarrow j)\in Q_1\} & \text{ if }i\neq \btau i,\btau j=j,\\
\{\alpha+\btau \alpha \mid (\alpha:i\rightarrow j) \text{ or }(\alpha:i\rightarrow \btau j)\in Q_1\} & \text{ if }i\neq \btau i,\btau j\neq j.\end{array}\right. \label{eqn:basis of R}
\end{eqnarray}
%Then $_jL_i$ (respectively, $_jR_i$) is a basis of ${}_j\BH_i $ as a left $\BH _j$-modules (respectively, as a right $\BH _i$-modules).

%
%
\subsubsection{}

We denote
\begin{equation}  \label{eq:ovOmega}
\overline{\Omega}:=\{(i,j)\in \ci \times \ci \mid (i,j)\in\Omega\text{ or }(i,\btau j)\in\Omega\}.
\end{equation}
Recall that $\ci $ is a (fixed) subset of $Q_0$ formed by the representatives of all $\btau$-orbits. The tuple $(\oH_i,\,_j\oH_i):=(\oH_i,\,_j\oH_i)_{i\in\ci ,(i,j)\in\ov{\Omega}}$  is called a \emph{modulation} of $(Q, \btau)$ and is denoted by $\cm(Q, \btau)$.

A representation $(N_i,N_{ji}):=(N_i,N_{ji})_{i\in\ci ,(i,j)\in\ov{\Omega}}$ of $\cm(Q, \btau)$ is defined by assigning to each $i\in \ci$ a finite-dimensional ${\ov{\BH}} _i$-module $N_i$ and to each $(i,j)\in \overline{\Omega}$ an ${\ov{\BH}} _j$-morphism
$N_{ji}:{}_j{\ov{\BH}}_i \otimes_{{\ov{\BH}} _i} N_i\rightarrow N_j$. A morphism $f:L\rightarrow N$ between representations $L=(L_i,L_{ji})$ and $N=(N_i,N_{ji})$ of $\cm(Q, \btau)$ is a tuple $f=(f_i)_{i\in \ci}$ of ${\ov{\BH}} _i$-morphisms $f_i:L_i\rightarrow N_i$ such that the following diagram is commutative for each $(i,j)\in\overline{\Omega}$:
\[\xymatrix{_j{\ov{\BH}}_i\otimes_{{\ov{\BH}} _i} L_i \ar[r]^{1\otimes f_i} \ar[d]^{L_{ij}}&  _j{\ov{\BH}}_i \otimes_{{\ov{\BH}} _i} N_i\ar[d]^{N_{ij}}\\
 L_j\ar[r]^{f_j} & N_j}\]

\begin{proposition}
   \label{prop:modulated representation}
The categories $\rep(\cm(Q, \btau))$ and $\rep(\ov{Q},\ov{I})$ are isomorphic.
\end{proposition}

\begin{proof}
The proof is the same as \cite[Proposition 2.16]{LW19}, hence omitted here.
\end{proof}

\subsubsection{ }

The materials in this subsection are inspired by \cite{GLS17,LW21}  and will be used in \S\ref{subsec:APR} to define reflection functors for $\imath$quivers.

Let $Q^*$ be the quiver constructed from $Q$ by reversing all the arrows $\alpha:i\rightarrow j$ such that $\btau i\neq j$. For any $i,j\in\I$ such that $\btau i\neq j$, we have $(i,j)\in\Omega$ if and only if $(j,i)\in \Omega^*:=\Omega(Q^*)$. For any $\alpha:i\rightarrow j$ in $Q$ such that $\btau i\neq j$, denote by $\widetilde{\alpha}:j\rightarrow i$ the corresponding arrow in $Q^*$.   Then $\btau$ induces an involution $\btau^*$ of $Q^*$. Clearly, $\btau^* i=\btau i$ for any vertex $i\in Q_0$. Then similarly we can define $\Lambda^*=\K Q^*\otimes_\K R_2$, and an involution $\btau^{*\sharp}$ for $\Lambda^*$, and its $\btau^{*\sharp}$-fixed point subalgebra $(\Lambda^*)^\imath$. Note that $\oH$ is also a subalgebra of $(\Lambda^*)^\imath$.

It is worth noting that $\ci$ is also a subset of representatives of $\btau^*$-orbits. In this way, one can define $\ov{\Omega}^*$ (cf. \eqref{eq:ovOmega} for $\ov{\Omega}$).

For any $(j,i)\in\Omega^*$, we can define
$_i\oH_j$ as follows:
\begin{align*}
_i\oH_j:= \oH_i\Span_k\{\widetilde{\alpha},\btau^* \widetilde{\alpha} \mid(\widetilde{\alpha}:j\rightarrow i)\in Q^*_1\text{ or }(\widetilde{\alpha}:j\rightarrow \btau^* i)\in Q^*_1\}\oH_j;
\end{align*}

Recall from \eqref{basis of Hij left}--\eqref{basis of Hij right} the basis $_i\LL_j$ (and respectively, $_i\RR_j$) for ${}_i\oH_j $ as a left $\oH _i$-module (and respectively, right $\oH _j$-module).
Let $_j\LL_i^*$ and $_j\RR_i^*$ be the dual bases of $\Hom_{\oH_j}(_j \oH_i, \oH_j)$ and $\Hom_{\oH_i}(_j\oH_i,\oH_i)$, respectively. Denote by $b^*$ the corresponding dual basis vector for any $b\in{_j\LL_i}$ or $b\in{_j\RR_i}$.

Since ${}_i\oH_j$ and $\Hom_{\oH_j}(_j\oH_i,\oH_j)$ are right free $\oH_j$-modules with bases given by ${}_i\RR_j$ and ${}_j\LL_i^*$ respectively, there is a right $\oH_j$-module isomorphism
\[
\rho:\, _i\oH_j \longrightarrow \Hom_{\oH_j}(_j\oH_i,\oH_j)
\]
such that $\rho(\widetilde{b})=b^*$ for any $b\in{}_j\LL_i$.
It is then routine to check that $\rho$ is actually an $\oH_i$-$\oH_j$-bimodule isomorphism.

Similarly, there is an $\oH_i$-$\oH_j$-bimodule isomorphism
$$\lambda: \, _i\oH_j\longrightarrow \Hom_{\oH_i}(_j\oH_i,\oH_i).$$
These two isomorphisms satisfy that $\rho(_i\RR_j)={}_j\LL_i^*$ and $\lambda(_i\LL_j)={}_j\RR_i^*$. We sometimes identify the spaces $\Hom_{\oH_j}(_j\oH_i,\oH_j)$, $_i\oH_j$ and $\Hom_{\oH_i}(_j\oH_i,\oH_i)$ via $\rho$ and $\lambda$.

If $N_j$ is an $\oH_j$-module, then we have a natural isomorphism of $\oH_i$-modules
$$\Hom_{\oH_j}(_j\oH_i,N_j)\longrightarrow\, _i\oH_j\otimes_{\oH_j}N_j$$
defined by
$$f\mapsto \sum_{b\in_j\LL_i} b^*\otimes f(b).$$
Furthermore, for any $\oH_i$-module $L_i$, there is a natural isomorphism of $\K$-vector spaces:
$$\Hom_{\oH_j}(_j\oH_i\otimes_{\oH_i} L_i, N_j)\longrightarrow \Hom_{\oH_i}(L_i,\Hom_{\oH_j}(_j\oH_i,N_j)).$$
Composing the two maps above, we obtain the following.

\begin{lemma}[cf. \text{\cite[Lemma 3.2]{LW21}}]
  \label{lem:ad}
There exists a canonical $k$-linear isomorphism
\begin{align*}
{\rm ad}_{ji}  ={\rm ad}_{ji}(L_i,N_j)  : & \Hom_{\oH_j}(_j\oH_i\otimes_{\oH_i} L_i,N_j)\longrightarrow \Hom_{\oH_i}(L_i,\,_i\oH_j\otimes_{\oH_j}N_j)
\\
{\rm ad}_{ji}: & f\mapsto \big(f^\vee:l\mapsto \sum_{b\in_j\LL_i} b^*\otimes f(b\otimes l) \big).
\end{align*}
The inverse ${\rm ad}_{ji}^{-1}$ is given by
${\rm ad}_{ji}^{-1} (g) = \big(g^\vee:h\otimes l\mapsto \sum_{b\in_j\LL_i} b^*(h)(g(l))_b \big),$ where the elements $(g(l))_b\in N_j$ are uniquely determined by
$g(l)=\sum_{b\in_j\LL_i} b^*\otimes(g(l))_b.$
\end{lemma}

%\begin{proof}
%The proof is the same as \cite[Lemma 3.2]{LW21}, and hence omitted here.
%\end{proof}

%
%
\subsection{Reflection functors}
  \label{subsec:APR}

In this subsection, we shall introduce the reflection functors in the setting of $\imath$quivers.

Let $(Q, \btau)$ be a virtually acyclic $\imath$quiver. Without loss of generality, we assume $Q$ to be connected and of rank $\ge 2$. Recall $\Omega=\Omega(Q)$ is the orientation of $Q$. For any sink $\ell \in Q_0$, define the quiver $s_\ell (Q)$ by reversing all the arrows ending to $\ell $. Note that in this case, we have $n_{\ell,\tau \ell}=0$.
By definition, $\ell $ is a sink of $Q$ if and only if $\btau \ell $ is a sink of $Q$.  Define the quiver
\begin{align}  \label{eq:QQ}
Q' = \bs_\ell Q =\left\{ \begin{array}{cc} s_\ell (Q) & \text{ if } \btau \ell =\ell ,\\ s_{\ell }s_{\btau\ell  }(Q) &\text{ if }\btau \ell \neq \ell .  \end{array}\right.
\end{align}
Note that $s_\ell s_{\btau \ell }(Q)=s_{\btau \ell }s_\ell (Q)$.
Then $\btau$ also induces an involution $\btau$ on the quiver $Q'$. In this way, we can define $\Lambda'=\K Q'\otimes_k R_2$ with an involution $\btau^\sharp$, and denote the $\btau^\sharp$-fixed point subalgebra by $\Lambda'^\imath =\bs_{\ell }\Lambda^{\imath}$. Note that
$\bs_{\ell }\Lambda^{\imath}=\bs_{\btau \ell }\Lambda^{\imath}$
for any sink $\ell \in Q_0$.
The quiver $\ov{Q'}$ of $\bs_\ell \Lambda^{\imath}$ can be constructed from $\ov{Q}$ by reversing all the arrows ending to $\ell $ and $\btau \ell $. Denote by $\Omega':=\Omega(Q')$ the orientation of $Q'$.

We shall define a {\rm reflection functor} associated to a sink $\ell \in Q_0$ (compare \cite{GLS17})
\begin{align}
F_\ell ^+:\mod(\Lambda^{\imath}):=\rep(\ov{Q},\ov{I})\longrightarrow\mod(\bs_\ell \iLa):=\rep(\ov{Q'},\ov{I'}),
\end{align}
in \eqref{eq:F+} below.
Using Proposition \ref{prop:modulated representation}, we identify the category $\rep(\ov{Q},\ov{I})$ with $\rep(\cm(Q, \btau))$, and respectively, $\rep(\ov{Q'},\ov{I'})$ with $\rep(\cm(Q',\btau))$.

Without loss of generality, we assume that the sink $\ell \in \ov{\I}_\btau$. Let $L=(L_i,L_{ji})\in \rep(\cm(Q, \btau))$. Then $L_{ji}:\,_j\oH_i\otimes_{\oH_i} L_i\rightarrow L_j$ is a $\oH_j$-morphism for any $(i,j)\in\ov{\Omega}$. Denote by
\[
L_{\ell ,{\rm in}}:= (L_{\ell i})_i:\bigoplus_{i\in\ov{\Omega}(-,\ell )}  \,_\ell \oH_i\otimes_{\oH_i} L_i\longrightarrow L_\ell.
\]

Let $N_\ell :=\ker(L_{\ell ,{\rm in}})$. Note that $\oH_\ell=\BH_\ell$ is finite-dimensional by our assumption. We have $\dim N_\ell<\infty$. By definition, there exists an exact sequence
\begin{align}
\label{def:reflection}
0\longrightarrow N_\ell \longrightarrow \bigoplus_{i\in\ov{\Omega}(-,\ell )}  \,_\ell \oH_i\otimes_{\oH_i} L_i\xrightarrow{L_{\ell ,{\rm in}}} L_\ell .
\end{align}
Denote by
$(N_{i\ell }^\vee)_i$ the inclusion map $N_\ell \rightarrow \bigoplus_{i\in\ov{\Omega}(-,\ell )}\,_\ell \oH_i\otimes_{\oH_i} L_i $.

For any $L\in \rep(\cm(Q, \btau))$, define
\begin{align}
  \label{eq:F+}
F_\ell ^+(L)=(N_r,N_{rs})\in \rep(\cm(Q',\btau)),
\end{align}
where
\[
N_r:=\left\{ \begin{array}{cc} L_r & \text{ if }r\neq \ell ,
\\
N_{\ell } &\text{ if }r=\ell,
\end{array} \right.
\qquad
 N_{rs}:=\left\{ \begin{array}{cc} L_{rs} &\text{ if } (s,r)\in \ov{\Omega} \text{ with }r\neq \ell ,
\\(N_{r\ell }^\vee)^\vee & \text{ if } (s,r)\in\ov{\Omega}^* \text{ and } s=\ell.
  \end{array}  \right.
  \]
Here $(N_{r\ell }^\vee)^\vee={\rm ad}_{r\ell }^{-1}(N_{r\ell }^\vee)$; see Lemma~\ref{lem:ad}.

Dually, associated to any source $\ell \in Q_0$, we have a reflection functor
\begin{align}
F_\ell ^-:\mod(\Lambda^{\imath}) \longrightarrow\mod(\bs_\ell \Lambda^{\imath}).
\end{align}

The above constructions are obviously functorial; cf. \cite[\S3.2]{LW21}. It is straightforward to show that $F_\ell^+$ is left exact, and $F_\ell^-$ is right exact. Both functors are covariant, $k$-linear and additive.

\subsection{Torsion pairs}

Let $\ell\in Q$ be a sink. For any $A\in\{\Lambda^\imath,\bs_\ell\Lambda^\imath\}$, and any $j\in Q_0$ let
\begin{align*}
\ct_j^A:=\{X\in\mod(A)\mid \Hom_A(X,S_j\oplus S_{\btau j})=0\},\\
\cs_j^A:=\{X\in\mod(A)\mid \Hom_A(S_j\oplus S_{\btau j},X)=0\}.
\end{align*}
For any $M\in\mod(A)$, we denote by $\sub_j(M)$ the largest submodule $U$ of $M$ supported at $j$ and $\btau j$; by $\fac_j(M)$ the largest quotient module $M/V$ of $M$ supported at $j$ and $\btau j$.

\begin{proposition}
\label{prop:adjoint}
For any sink vertex $\ell\in Q$ the following hold:
\begin{enumerate}
\item The pair $(F_\ell^-,F_\ell^+)$ is a pair of adjoint functors, i.e., there is a functorial isomorphism
\begin{align*}
\Hom_{\Lambda^\imath}(F_\ell^-(M),N)\cong \Hom_{\bs_\ell\Lambda^\imath}(M,F_\ell^+(N)).
\end{align*}
\item The adjunction morphisms $\id\rightarrow F_\ell^+F_\ell^-$ and $F_\ell^-F_\ell^+\rightarrow \id$ can be inserted in functorial short exact sequences
$$0\longrightarrow \sub_\ell\longrightarrow \id \longrightarrow F_\ell^+F_\ell^-\longrightarrow0,\qquad 0\longrightarrow F_\ell^-F_\ell^+\longrightarrow \id\longrightarrow \fac_\ell\longrightarrow0.$$
\end{enumerate}
\end{proposition}

\begin{proof}
The proof is the same as \cite[Proposition 9.1]{GLS17}. For the sake of completeness, we give the proof here.

For (1), it is enough to construct a pair of mutual inverses between $\Hom_{\Lambda^\imath}(F_\ell^-(M),N)$ and $\Hom_{\bs_\ell\Lambda^\imath}(M,F_\ell^+(N))$, which are functorial in $M$ and $N$. The construction is as follows.

Let $Q'=\bs_\ell (Q)$ and $M=(M_r,M_{rs})\in \rep(\cm(Q',\btau))$. Let $N=(N_r,N_{rs})\in \rep\cm(Q,\btau)$.
Denote by $(U_r,U_{rs}):=F_\ell^-(M)\in\rep(\cm(Q,\btau))$, and by $(V_r,V_{rs}):=F_\ell^+(N)\in\rep(\cm(Q',\btau))$. In particular, $U_r=M_r$ and $V_r=N_r$ if $r\neq\ell$.
For any $f=(f_r)_{r\in\I_\btau}:M\rightarrow F_\ell^+(N)$, where $f_r:M_r\rightarrow V_r$, by definition, we have the following commutative diagram with exact rows:
\[\xymatrix{& M_\ell \ar[rr]^{(M_{i\ell}^\vee)_i\qquad\qquad} \ar[d]^{f_\ell}  &&  \bigoplus_{i\in\ov{\Omega}(-,\ell)} {}_\ell\oH_i \otimes_{\oH_i} M_i \ar[rr]^{\qquad\qquad(U_{\ell i})_i} \ar[d]^{(\Id\otimes f_i)_i} & & U_\ell \ar@{.>}[d]^{g_\ell}\ar[r] &0
\\
0\ar[r] & V_\ell \ar[rr]^{(V_{i\ell}^\vee)_i\qquad\qquad}  && \bigoplus_{i\in\ov{\Omega}(-,\ell)} {}_\ell\oH_i \otimes_{\oH_i} N_i\ar[rr]^{\qquad\qquad(N_{\ell i})_i} && N_\ell &  }\]
Then there exists a unique morphism $g_\ell:U_\ell\rightarrow N_\ell$ such that the above diagram commutes. Define $g:=(g_r)_{r\in\I_\btau}$ such that
$$g_r=\left\{ \begin{array}{cc} f_r & \text{ if }r\neq \ell ,
\\
g_{\ell } &\text{ if }r=\ell.
\end{array}  \right. $$
Clearly, $g:F_\ell^-(M)=U\rightarrow N$ is a morphism of representations.

Conversely, for any $g:F_\ell^-(M)\rightarrow N$, one can construct a morphism $f:M\rightarrow F_\ell^+(N)$. The above constructions are naturally functorial in $M$ and $N$, and then (1) follows.

For (2),  we only prove the second one. For any $N=(N_r,N_{rs})\in \rep\cm(Q,\btau)$, by definition, it is obvious that the adjunction $\alpha_N:F_\ell^-F_\ell^+(N)\rightarrow N$ is injective, and $\coker(\alpha_N)=\coker(N_{\ell,{\rm in}})$ which is supported at $\ell$ and $\btau(\ell)$. Since $\ell$ is a sink, we have $\coker(N_{\ell,{\rm in}})\cong \fac_\ell(N)$, and then obtain the desired short exact sequence
\begin{align*}
0\longrightarrow F_\ell^-F_\ell^+(N)\stackrel{\alpha_N}{\longrightarrow} N\longrightarrow \fac_\ell(N)\longrightarrow0.
\end{align*}
The proposition is proved.
\end{proof}

We have the following corollaries of Proposition \ref{prop:adjoint}.

\begin{corollary}
We have
the following equivalence of subcategories:
\begin{align}
\label{eqn:YY}
F^+_\ell: \ct_\ell^{\Lambda^\imath}\stackrel{\simeq}{\longrightarrow} \cs_\ell^{\bs_\ell\Lambda^\imath},
\end{align}
with its inverse given by $F^-_\ell$.
\end{corollary}

\begin{corollary}
\begin{enumerate}
\item For $M,N\in\ct^{\Lambda^\imath}_\ell$, $F_\ell^+$ induces an isomorphism
\begin{align}
\Ext^1_{\Lambda^\imath}(M,N)\cong \Ext^1_{\bs_\ell \Lambda^\imath}(F_\ell^+(M),F_\ell^+(N)).
\end{align}
\item For $M,N\in\cs^{\bs_\ell\Lambda^\imath}_\ell$, $F_\ell^-$ induces an isomorphism
\begin{align}
\Ext^1_{\bs_\ell\Lambda^\imath}(M,N)\cong \Ext^1_{ \Lambda^\imath}(F_\ell^-(M),F_\ell^-(N)).
\end{align}
\end{enumerate}
\end{corollary}

Let $\ct:=\ct_\ell^{\Lambda^\imath}$, and let $\cf$ be the extension closed subcategory of $\mod(\Lambda^\imath)$ generated by $S_\ell$ and $S_{\btau \ell}$.

\begin{lemma}
\label{lem: torsion pair}
{\quad}
\begin{itemize}
\item[(a)] $(\ct,\cf)$ is a torsion pair in $\mod(\Lambda^\imath)$;
\item[(b)] For any $M\in \mod(\Lambda^\imath)$, there exists a short exact sequence
\[
0 \longrightarrow M \longrightarrow T_M^0 \longrightarrow T_M^1 \longrightarrow 0
\]
with $T_M^0,T_M^1\in \ct$ and $T_M^1\in\cp^{\leq1}(\Lambda^\imath)$.
\end{itemize}
\end{lemma}

\begin{proof}
(a) First, $\Hom_{\Lambda^\imath}(\ct,\cf)=0$.
For any $M\in\mod(\Lambda^\imath)$, let $tM$ be the maximal submodule of $M$ such that its top $\Top M$ satisfies $(e_{\ell}+e_{\btau\ell}) \Top M=0$. Since $\ell$ (also $\btau\ell$) is a sink, we have $M/tM$ concentrated at the full subquiver formed by $\ell$ and $\btau \ell$, and then it is a natural $\BH_\ell$-module. So $M/tM\in\cf$.
It follows an exact sequence $0\rightarrow tM\rightarrow M\rightarrow M/tM\rightarrow0$. So $(\ct,\cf)$ is a torsion pair in $\mod(\Lambda^\imath)$.

(b)
First, we prove it for any $M\in \cf$. Note that $\cf=\add\{S_\ell,S_{\btau \ell},\E_{\ell},\E_{\btau \ell}\}$. Without loss of  generality, we assume $M$ is indecomposable.

Case (1) $\underline{M=S_\ell, S_{\btau \ell}}$.
We only prove for $M=S_\ell$. As $\ell$ is a sink, there exist at least one arrow $\alpha:j\rightarrow \ell$ in $Q_0$. So there exists a string module $X$ with its string $\ell\xleftarrow{\alpha} j\xrightarrow{\varepsilon_j}\btau j$.
Then $X,\E_j\in\ct$, and there exists a short exact sequence
$0\rightarrow S_\ell\rightarrow X\rightarrow \E_j\rightarrow0$.

Case (2) $\underline{M=\E_\ell, \E_{\btau \ell}}$.
We only prove for $M=\E_\ell$. We have the following exact sequence $0\rightarrow S_{\btau \ell}\rightarrow \E_\ell\rightarrow S_{\ell}\rightarrow0$.
Since $\pd\E_j\leq1$, and by Case (1), we have the following commutative diagram
\[\xymatrix{S_{\btau\ell}\ar[r] \ar@{=}[d]& \E_\ell \ar[r] \ar[d] & S_\ell \ar[d] \\
S_{\btau \ell} \ar[r] & Y\ar[r] \ar[d] & X\ar[d] \\
 & \E_j \ar@{=}[r] & \E_j }\]
Consider the short exact sequence in the second row. Note that $\btau$ induces an equivalence $\btau$ of $\mod(\Lambda^\imath)$.
From Case (1), we obtain an exact sequence $0\rightarrow S_{\btau \ell}\rightarrow \btau(X)\rightarrow \E_{\btau j}\rightarrow0$ with $\btau(X),\E_{\btau j}\in\ct$.
Then we have the pushout diagram
\[\xymatrix{S_{\btau \ell} \ar[r] \ar[d] & Y\ar[r] \ar[d] & X\ar@{=}[d] \\
\btau(X) \ar[r] \ar[d] & W\ar[r] \ar[d] & X\\
\E_{\btau j} \ar@{=}[r] & \E_{\btau j}}\]
Since $\ct$ is closed under extensions, we have $W\in\ct$ by $X,\btau(X)\in\ct$.
Combining these two diagrams, we have the exact sequences
$0\rightarrow \E_\ell\rightarrow W\rightarrow U\rightarrow0$ and $0\rightarrow \E_{\btau j}\rightarrow U\rightarrow \E_{j}\rightarrow0$.
Hence, $U\in\ct\cap \cp^{\leq1}(\Lambda^\imath)$.

Next, for general $M$, we have the exact sequence
$0\rightarrow tM\rightarrow M\rightarrow M/tM\rightarrow0$ with $tM\in\ct$, $M/tM\in \cf$.
From the above, we obtain an exact sequence $0\rightarrow M/tM\rightarrow W\rightarrow T_M^1\rightarrow0$ with $W,T_M^1\in\ct$ and  $T_M^1\in\cp^{\leq1}(\Lambda^\imath)$.
By a similar argument as above, from these exact sequences one obtains the desired resolution $0\rightarrow M\rightarrow T_M^0\rightarrow T_M^1\rightarrow0$.
\end{proof}

Let $\cy=\cs^{\bs_\ell\Lambda^\imath}_\ell$, and $\cx$ be the extension closed subcategory of $\mod(\bs_\ell\Lambda^\imath)$ generated by $S_\ell,S_{\btau \ell}$. Note that $\bs_\ell\Lambda^\imath$ is a $1$-Gorenstein algebra.
Dually, we have the following lemma.
\begin{lemma}
\label{lem: torsion pair 2}
We have the following.
\begin{itemize}
\item[(a)] $(\cx,\cy)$ is a torsion pair in $\mod(\bs_\ell\Lambda^\imath)$;
\item[(b)] for any $M\in \mod(\bs_\ell\Lambda^\imath)$, there exists a short exact sequence
\[
0 \longrightarrow Y_M^1 \longrightarrow Y_M^0 \longrightarrow M \longrightarrow 0
\]
with $Y_M^0,Y_M^1\in \cy$ and $Y_M^1\in \cp^{\leq1}(\bs_\ell\Lambda^\imath)$.
\end{itemize}
\end{lemma}

We have a $\Z$-linear isomorphism $\dimv\colon K_0(\mod(\K Q)) \rightarrow \Z^\I$, which sends an isoclass to its dimension vector. By identifying $i \in \I$ with a simple root $\alpha_i$ and thus $\Z^\I$ with the root lattice (of a Kac-Moody algebra $\fg$), we have simple reflection $s_i$ acting on $\Z^\I$; see \S\ref{subsec:QG} for more detail. Then  we denote
\begin{align*}
\bs_i= \left\{
\begin{array}{ll}
s_{i}, & \text{ if } i=\btau i
\\
s_is_{\btau i}, & \text{ if } i\neq \btau i.
\end{array}
\right.
\end{align*}

\begin{lemma}
\label{lem:reflecting dimen}
Let $(Q, \btau)$ be an $\imath$quiver with a sink $\ell$. Let $L\in\mod(\K Q)\subseteq \mod(\iLa)$ be an indecomposable module. Then either $F_\ell^+(L)=0$ (equivalently $L\cong S_\ell$ or $S_{\btau \ell}$) or $F_\ell^+(L)$ is indecomposable with
$\dimv F_\ell^+(L)=\bs_\ell(\dimv L)$.
\end{lemma}

\begin{proof}
Consider an indecomposable module $L=(L_i,L_{ji})\in\cm(Q, \btau)$. If $L\cong S_\ell$ or $S_{\btau \ell}$, by definition of $F_\ell^+$, we have $F_\ell^+(L)=0$.

Otherwise, we have $L\in \ct$. As $F_\ell^+:\ct\rightarrow \cy$ is an equivalence (see \eqref{eqn:YY}), we have $\End_{\bs_\ell\Lambda^\imath}(F_\ell^+(L))\cong \End_{\iLa}(L)$ which is a local algebra since $L$ is indecomposable. So $F_\ell^+(L)$ is also indecomposable.

Let $F^+_{\ell}(L)=(N_i,N_{ji})\in\cm(Q, \btau)$.
Since $L$ is indecomposable and $\ell$ is a sink, it follows that the morphism $L_{\ell,{\rm in}}$ in \eqref{def:reflection} is surjective, and then \eqref{def:reflection} becomes a short exact sequence. We have
\begin{align*}
\dim_\K (e_{\ell}N_\ell) &=\sum_{(\alpha:i\rightarrow \ell)\in Q_1} \dim_\K (e_i L_i)-\dim_\K (e_{\ell} L_\ell),
  \\
\dim_\K (e_{\btau\ell}N_\ell)& =\sum_{(\alpha:i\rightarrow \btau\ell)\in Q_1} \dim_\K (e_{\btau i} L_i)-\dim_\K (e_{\btau\ell} L_\ell).
\end{align*}
Since $N_i=L_i$ for $i\neq \ell$, from the above, we conclude that $\dimv F^+_{\ell}(L)=\bs_\ell(\dimv L)$.
\end{proof}

Dual results also hold for $F_\ell^-$ and for any source $\ell$ of $Q$.

%%%%%%%%%%%%%%%%

%%
\subsection{$\imath$Hall algebras}

In this subsection we consider $k=\F_q$, and set
\[
\sqq=\sqrt{q}.
\]

Generalizing \cite{LP21}, the first author defined a (twisted) semi-derived Hall algebra of a 1-Gorenstein algebra \cite{Lu19}. The $\imath$Hall algebra $\tMH$ for $\imath$quiver $(Q,\btau)$ is by definition the twisted semi-derived Hall algebra for the module category of the $\imath$quiver algebra $\iLa$ \cite{LW19,LW20a}. We recall it here briefly.

Let $\ch(\Lambda^\imath)$ be the Ringel-Hall algebra of $\Lambda^\imath$, i.e.,
$$\ch(\Lambda^\imath)=\bigoplus_{[M]\in\Iso(\mod(\Lambda^\imath))}\Q(\sqq)[M],$$
with the multiplication defined by (see \cite{Br13})
\[
[M]\diamond [N]=\sum_{[M]\in\Iso(\mod(\Lambda^\imath))}\frac{|\Ext^1(M,N)_L|}{|\Hom(M,N)|}[L].
\]

For any three objects $X,Y,Z$, let
\begin{align}
F_{XY}^Z &= \big|\{L\subseteq Z, L \cong Y\text{ and }Z/L\cong X\} \big|
\notag \\
&= \frac{\big|\Ext^1(X,Y)_Z|}{|\Hom(X,Y)\big|} \cdot \frac{|\aut(Z)|}{|\aut(X)| |\aut(Y)|}
\qquad \text{(Riedtman-Peng formula)}.
\label{Ried-P}
\end{align}

Define $I$ to be the two-sided ideal of $\ch(\Lambda^\imath)$ generated by
\begin{align}
  \label{eq:ideal}
&\{[K]-[K'] \mid \res_\BH(K)\cong\res_\BH(K'),  K,K'\in\cp^{<\infty}(\Lambda^\imath)\} \bigcup
\\\notag
&\{[L]-[K\oplus M]\mid \exists \text{ exact sequence } 0 \longrightarrow K \longrightarrow L \longrightarrow M \longrightarrow 0, K\in\cp^{<\infty}(\Lambda^\imath)\}.
\end{align}
%Here $\cp^{<\infty}(\Lambda^\imath)$ is the subcategory of $\mod(\Lambda^\imath)$ formed by modules of finite projective dimensions.

Consider the following multiplicatively closed subset $\cs$ of $\ch(\Lambda^\imath)/I$:
\begin{equation}
  \label{eq:Sca}
\cs = \{ a[K] \in \ch(\Lambda^\imath)/I \mid a\in \Q(\sqq)^\times, K\in \cp^{<\infty}(\Lambda^\imath)\}.
\end{equation}

The semi-derived Hall algebra of $\Lambda^\imath$ \cite{Lu19} is defined to be the localization
$$\cs\cd\ch(\Lambda^\imath):= (\ch(\Lambda^\imath)/I)[\cs^{-1}].$$

Let $\langle\cdot,\cdot \rangle_Q$ be the Euler form of $Q$.
%Via the restriction functor $\res: \mod(\Lambda^{\imath})\rightarrow\mod (\K Q)$, w
We define the $\imath$Hall algebra (i.e., a twisted semi-derived Hall algebra) $\tMH$ \cite[\S4.4]{LW19} to be the $\Q(\sqq)$-algebra on the same vector space as $\utMH$ but with twisted multiplication given by
\begin{align}
   \label{eqn:twsited multiplication}
[M]* [N] =\sqq^{\langle \res(M),\res(N)\rangle_Q} [M]\diamond[N].
\end{align}

\subsection{Symmetries of $\imath$Hall algebras}
  \label{sec:symmetryHall}

We shall formulate the reflection functor associated to a sink $\ell \in Q_0$, which induces an $\imath$Hall algebra isomorphism $\Gamma_{\ell}: \tMH \stackrel{\cong}{\rightarrow}  \tMHl$.

Let $(Q, \btau)$ be an $\imath$quiver.  We assume $Q$ to be connected and of rank $\ge 2$. Let $\ell$ be a sink in $Q$.
Recall $Q'= \bs_\ell Q$ from \eqref{eq:QQ}. As in \S\ref{subsec:APR}, $\btau$ induces an involution of $Q'$ which is also denoted by $\btau$. Let $\bs_\ell\iLa$ denote the $\imath$quiver algebra associated to $(Q', \btau)$.
Recall Lemma \ref{lem: torsion pair}. Similar to the proof of \cite[Theorem A.22]{LW19}, we have the following.

\begin{lemma}
  \label{lem:Upsilon}
  Let $\ell$ be a sink of $Q$. Then we have an isomorphism of algebras:
\begin{align}
  \label{eqn:reflection functor 1}
\Gamma_{\ell}:\utMH&\stackrel{\cong}{\longrightarrow} \cs\cd\ch(\bs_\ell\iLa)\\\notag
[M]&\mapsto q^{-\langle T_M,M\rangle} [F_\ell^+(T_M)]^{-1}\diamond [F_\ell^+(X_M)],
\end{align}
where $M\in\mod(\Lambda^\imath)$ and $X_M,T_M\in \ct$ ($T_M\in\cp^{\leq1}(\Lambda^\imath)$) fit into a short exact sequence
$0\rightarrow M\rightarrow X_M\rightarrow T_M\rightarrow0$.
\end{lemma}

\begin{proof}
By a similar proof to \cite[Theorem A.22]{LW19}, $\Gamma_\ell$ is well defined and is an algebra morphism.

On the other hand, by Lemma \ref{lem: torsion pair 2}, we obtain a morphism
$\Gamma_\ell^-: \cs\cd\ch(\bs_\ell\iLa)\rightarrow \utMH$, which maps
$N\mapsto q^{-\langle U_N, N  \rangle}[F_\ell^-(U_N)]^{-1}\diamond [F_\ell^-(Y_N)]$,
where $N$, $Y_N,U_N\in\cy$ (and $U_N\in \cp^{\leq1}(\bs_\ell(\Lambda^\imath))$) fit into a short exact sequence
$0\rightarrow U_N\rightarrow Y_N\rightarrow N\rightarrow0$.

Since $F_\ell^+:\ct\longrightarrow \cy$ is an equivalence with $F_\ell^-$ as its inverse, we have
\begin{align*}
\Gamma_\ell^-\circ \Gamma_\ell([M]) &=  q^{-\langle T_M,M\rangle} [T_M]^{-1}\diamond [X_M]=[M],
\\
\Gamma_\ell\circ \Gamma_\ell^-([N]) &=  q^{-\langle U_N, N  \rangle}[U_N]^{-1}\diamond [Y_N]=[N]
\end{align*}
for any $M\in\mod(\Lambda^\imath)$, $N\in\mod(\bs_\ell\Lambda^\imath)$. It follows that $\Gamma_\ell$ and $\Gamma_\ell^-$ are  inverses to each other.
\end{proof}

It is well known that the Cartan matrix $C$ for $Q$ is the matrix of the symmetric bilinear form $(\cdot,\cdot)_Q$ defined by
\[
(x,y)_Q:=\langle x,y\rangle_Q+\langle y,x\rangle_Q
\]
for any $x,y\in K_0(\mod(\K Q))$. Here $\langle\cdot,\cdot\rangle_Q$ is the Euler form of $\K Q$. Recall $\sqq =\sqrt{q}$.

\begin{theorem}%[\text{\cite[Theorem 4.3]{LW21}}]
  \label{thm:Gamma}
The isomorphism $\Gamma_\ell$ in (\ref{eqn:reflection functor 1}) induces the following isomorphism of $\imath$Hall algebras: %, also denoted by $\Gamma_{\ell}$:
\begin{eqnarray}
\Gamma_{\ell}:\tMH & \stackrel{\cong}{\longrightarrow} & \tMHl,
  \label{eqn:reflection functor 2} \\
\, [M]&\mapsto& \sqq^{\langle \res (T_M),\res(M)\rangle_Q}q^{-\langle T_M,M\rangle} [F_{\ell}^+(T_M)]^{-1}* [F_{\ell}^+(X_M)].
 \notag
\end{eqnarray}
%where $0\rightarrow M\rightarrow X_M\rightarrow T_M\rightarrow0$ is a short exact sequence with $X_M\in \Fac T$ and $T_M\in \add T$.
\end{theorem}

\begin{proof}
Similar to \cite[Lemma 4.2]{LW21}, the Euler forms match, and then the result follows by Lemma~\ref{lem:Upsilon} (to which we refer for notations).
\end{proof}

For any $i\in \I$, denote by $S_i$ (respectively, $S_i'$) the simple $\iLa$-module (respectively, $\bs_\ell\iLa$-module), denote by $\E_i$ (respectively, $\E_i'$) the
generalized simple $\iLa$ (respectively, $\bs_\ell\iLa$-module). We similarly define $\E_\alpha,\E_\beta'$ for $\alpha\in K_0(\mod(\K Q))$, $\beta\in K_0(\mod(\K(\bs_\ell Q)))$
in the (twisted) semi-derived Hall algebras (where $\ell$ is a sink of $Q$).

Recall the root lattice $\Z^{\I}=\Z\alpha_1\oplus\cdots\oplus\Z\alpha_n$, and we have an isomorphism of abelian groups $\Z^\I\rightarrow K_0(\mod(\K Q))$, $\alpha_i\mapsto \widehat{S_i}$. This isomorphism induces the action of the reflection $s_i$ on $K_0(\mod(\K Q))$.
Thus for $\alpha\in K_0(\mod (\K Q))$ and $i\in\I$, we can make sense $[\E_{\bs_i\alpha}] \in \tMH$. Similarly,  we have $[\E'_{\bs_i\alpha}]\in \tMHl$.

\begin{proposition} [\text{\cite[Proposition 4.4]{LW21}}]
   \label{prop:reflection}
Let $(Q, \btau)$ be an $\imath$quiver, and $\ell\in Q_0$ be a sink.
Then the isomorphism $\Gamma_{\ell}:\tMH\xrightarrow{\cong}  \tMHl$ sends
\begin{align}
\Gamma_{\ell}([M])&= [F_{\ell}^+(M)], \quad \forall M\in\ct,
 \label{eqn:reflection 1}
  \\
\Gamma_{\ell}([S_\ell])&= \left\{
\begin{array}{cc}\sqq[\E'_\ell]^{-1}* [S_{\btau \ell}'], &\text{ if }\btau \ell\neq \ell,
\\ \,[\E'_\ell]^{-1}* [S_{\btau \ell}'], &\text{ if }\btau \ell= \ell, \end{array}\right.
\label{eqn:reflection 2}
 \\
\Gamma_{\ell}([S_{\btau \ell}])&=
\sqq [\E'_{\btau \ell}]^{-1} * [S_{\ell}'], \quad\quad \text{ if }\btau \ell\neq \ell,
\label{eqn:reflection 3}
 \\
\Gamma_{\ell}([\E_\alpha])&= [\E'_{\bs_{\ell}\alpha}],
\quad \forall\alpha\in K_0(\mod(\K Q)).\label{eqn:reflection 4}
\end{align}
\end{proposition}

Similarly, one can give the formulas of $\Gamma^-_{\ell}:\tMH\xrightarrow{\cong}  \tMHl$ for any source $\ell$ of an $\imath$quiver $(Q,\btau)$.

\begin{remark}
 \label{rem:FT}
Similar to \cite{SV99, XY01}, for a sink $\ell \in \I$, there exists a Fourier transform (which is an algebra isomorphism)
$%\begin{align*}
\texttt{FT}_\ell : \tMHl \rightarrow \tMH,
$
%\end{align*}
which maps $[S'_j]\mapsto [S_j]$, $[\E_j']\mapsto [\E_j]$ for each $j\in \I$; compare Theorem \ref{thm:Ui=iHall}.
The composition of $\Gamma_\ell$ with $\texttt{FT}_\ell$ gives us an automorphism
$\texttt{FT}_\ell\circ\Gamma_\ell: \tMH \rightarrow \tMH.$
\end{remark}

\section{Formula for a reflection functor $\Gamma_i$ ($i =\btau i$)}
  \label{sec:formulaRF}

In this section, we establish a closed formula for the action of the reflection functor $\Gamma_i$ with $i =\btau i$.

\subsection{$\imath$Divided powers}

Let $v$ be an indeterminate.
Define the quantum integers, quantum (double) factorials, and quantum binomial coefficients, for $r \in \N$ and $m \in \Z$,
\begin{align}  \label{eq:binom}
\begin{split}
 [m]= [m]_v =\frac{\bv^m-\bv^{-m}}{\bv-\bv^{-1}},
 &\qquad
 [r]^{!}= [r]_v^! =\prod_{i=1}^r [i]_v,
  \\
[2r]^{!!} = [2r]_v^{!!} =\prod_{i=1}^r [2i]_v,
& \qquad
 \qbinom{m}{r} =\frac{[m][m-1]\ldots [m-r+1]}{[r]^!}.
 \end{split}
\end{align}
We often need to specialize them by substituting $v$ with $\sqq =\sqrt{q}$ below.

For a $\iLa$-module $M$, we shall write
\[
[l M] =[\underbrace{M\oplus \cdots \oplus M}_{l}],
\qquad
[M]^{l} = \underbrace{[M]* \cdots * [M]}_{l}.
\]

Let $\Z_2=\{\ov{0},\ov{1}\}$.
Following \cite{LW20a}, we define the $\imath$divided power of $[S_i]$ in $\tMH$ as follows:
\begin{align}
\label{eq:idividedHallodd}
&[S_i]_{\odd}^{(m)}:=\frac{1}{[m]_{\sqq}^!}\left\{ \begin{array}{ll} [S_i]\prod_{j=1}^k ([S_i]^2+\sqq^{-1}(\sqq^2-1)^2[2j-1]_{\sqq}^2 [\E_i] ) & \text{if }m=2k+1,
\\
\prod_{j=1}^k ([S_i]^2+\sqq^{-1}(\sqq^2-1)^2[2j-1]_{\sqq}^2[\E_i]) &\text{if }m=2k; \end{array}\right.
 \\
 \label{eq:idividedHallev}
&[S_i]_{\ev}^{(m)}:= \frac{1}{[m]_{\sqq}^!}\left\{ \begin{array}{ll} [S_i]\prod_{j=1}^k ([S_i]^2+\sqq^{-1}(\sqq^2-1)^2[2j]_{\sqq}^2[\E_i] ) &\text{if }m=2k+1,\\
\prod_{j=1}^{k} ([S_i]^2+\sqq^{-1}(\sqq^2-1)^2[2j-2]_{\sqq}^2[\E_i]) &\text{if }m=2k. \end{array}\right.
\end{align}

We recall the expansion formula of the $\imath$divided powers in terms of an $\imath$Hall basis. See \eqref{eq:binom} for notation $[2k]_\sqq^{!!}$.

\begin{lemma}
[\text{\cite[Propositions 6.4--6.5]{LW20a}}]
 \label{lem:iDPev}
For any $m \in \N$, $\ov{p}\in\Z_2$, we have
\begin{align}
&[S_i]^{(m)}_{\ov{p}}=\begin{cases}
\sum_{k=0}^{\lfloor\frac{m}{2}\rfloor} \frac{\sqq^{k(k-1)-\binom{m-2k}{2}}}{[m-2k]_{\sqq}^{!}[2k]_\sqq^{!!}}  (\sqq-\sqq^{-1})^k
[(m-2k)S_i]*[\E_i]^k, & \text{ if }\ov{m}=\ov{p};
 \label{eq:2mev}
 \\
 \\
 \sum_{k=0}^{\lfloor\frac{m}{2}\rfloor}  \frac{\sqq^{k(k+1) -\binom{m-2k}{2}}}{[m-2k]_{\sqq}^{!}[2k]_\sqq^{!!}} (\sqq-\sqq^{-1})^k
 [(m-2k)S_i]*[\E_i]^k, & \text{ if }\ov{m}\neq\ov{p}.
\end{cases}
\end{align}
\end{lemma}

\subsection{Formulas of $\Gamma_i$ for $i=\btau i$}

Below is the first main result of this paper.

\begin{theorem}
\label{thm:Braidsplit}
Let $(Q,\btau)$ be an $\imath$quiver.
For any sink $i\in Q_0$ such that $i=\btau i \neq j$, we have
\begin{align}
\label{eqn:Braidsplit}
\Gamma_i([S_j])=&\sum_{r+s=-c_{ij}}(-1)^r\sqq^{r}(1-\sqq^2)^{c_{ij}} [S'_i]_{\ov{p}}^{(r)}*[S'_j]*[S'_i]_{\ov{c_{ij}}+\ov{p}}^{(s)}\\
\notag&+ (-1)^{p}  \sum_{t\geq1}\sum_{\stackrel{r+s+2t=-c_{ij} }{\ov{r}=\ov{p} }}   \sqq^{r}(1-\sqq^2)^{c_{ij}+2t}[S'_i]_{\ov{p}}^{(r)}*[S'_j]*[S'_i]_{\ov{c_{ij}}+\ov{p}}^{(s)}*[\E'_i]^t.
\end{align}
\end{theorem}

The proof of Theorem \ref{thm:Braidsplit} will occupy the remainder of this section and Appendix~\ref{app:prop}. Let $a=-c_{ij}$. Then
\eqref{eqn:Braidsplit} is equivalent to the following formulas \eqref{eqn:braidsplitev}--\eqref{eqn:braidsplitodd} (where the 2 cases for $\ov{p}=\ov 0, \ov 1$ are separated):
\begin{align}
\label{eqn:braidsplitev}
\Gamma_i([S_j])=&\sum_{r+s=a}(-1)^r\sqq^{r}(1-\sqq^2)^{-a} [S'_i]_{\ev}^{(r)}*[S'_j]*[S'_i]_{\ov{a}}^{(s)}\\
\notag&+\sum_{t\geq1}\sum_{\stackrel{r+s+2t=a}{2\mid r}} \sqq^{r}(1-\sqq^2)^{-a+2t}[S'_i]_{\ev}^{(r)}*[S'_j]*[S'_i]_{\ov{a}}^{(s)}*[\E'_i]^t,
\\
\label{eqn:braidsplitodd}
\Gamma_i([S_j])=&\sum_{r+s=a}(-1)^r\sqq^{r}(1-\sqq^2)^{-a} [S'_i]_{\odd}^{(r)}*[S'_j]*[S'_i]_{\ov{1}+\ov{a}}^{(s)}\\
\notag& - \sum_{t\geq1}\sum_{\stackrel{r+s+2t=a}{2\nmid r}} \sqq^{r}(1-\sqq^2)^{-a+2t}[S'_i]_{\odd}^{(r)}*[S'_j]*[S'_i]_{\ov{1}+\ov{a}}^{(s)}*[\E'_i]^t.
\end{align}

The proof of \eqref{eqn:braidsplitev} for $j =\btau j$ will be given in \S\ref{subsec:SSS}--\S\ref{subsec:contribution d=0} and Section~\ref{sec:identities}, while a similar proof of \eqref{eqn:braidsplitodd} for $j =\btau j$ can be found in Appendix~\ref{app:braidsplitodd}.
This proves Theorem \ref{thm:Braidsplit} for $j =\btau j$.

We then explain how Theorem \ref{thm:Braidsplit} for $j \neq \btau j$ is reduced to the case for $j =\btau j$.

\subsection{Reduction for the formula \eqref{eqn:braidsplitev} with $\btau j=j$}
\label{subsec:braidsplitev}

For the case $\btau j=j$, it is enough to consider the rank 2 $\imath$quiver $Q$ with trivial involution $\btau=\Id$, as shown in the left figure of \eqref{diag:split KM}.  Here $a=-c_{ij}$. Then the quiver $\ov{Q'}$ of $\bs_i\Lambda^\imath$ is shown in the right figure of \eqref{diag:split KM}.
\begin{center}\setlength{\unitlength}{0.7mm}
\begin{equation}
\label{diag:split KM}
\begin{picture}(50,13)(0,0)
\put(-40,0){\begin{picture}(50,13)
\put(-23,0){$\ov{Q}=$}
\put(0,-3){$i$}
%\put(4,1){\vector(1,0){20}}
\put(12,-1){\vector(-1,0){9}}
\put(12.5,-2.2){$a$}
\put(16,-1){\line(1,0){8}}

\put(25,-3){$j$}
\color{purple}
\qbezier(-1,1)(-3,3)(-2,5.5)
\qbezier(-2,5.5)(1,9)(4,5.5)
\qbezier(4,5.5)(5,3)(3,1)
\put(3.1,1.4){\vector(-1,-1){0.3}}
\qbezier(24,1)(22,3)(23,5.5)
\qbezier(23,5.5)(26,9)(29,5.5)
\qbezier(29,5.5)(30,3)(28,1)
\put(28.1,1.4){\vector(-1,-1){0.3}}
\put(-1,10){$\varepsilon_i$}
\put(24,10){$\varepsilon_j$}
\end{picture}}
\put(60,0){\begin{picture}(50,13)
\put(-23,0){$\ov{Q'}=$}
\put(0,-3){$i$}
%\put(4,1){\vector(1,0){20}}
\put(4,-1){\line(1,0){8}}
\put(12.5,-2.2){$a$}
\put(16,-1){\vector(1,0){8}}

\put(25,-3){$j$}
\color{purple}
\qbezier(-1,1)(-3,3)(-2,5.5)
\qbezier(-2,5.5)(1,9)(4,5.5)
\qbezier(4,5.5)(5,3)(3,1)
\put(3.1,1.4){\vector(-1,-1){0.3}}
\qbezier(24,1)(22,3)(23,5.5)
\qbezier(23,5.5)(26,9)(29,5.5)
\qbezier(29,5.5)(30,3)(28,1)
\put(28.1,1.4){\vector(-1,-1){0.3}}
\put(-1,10){$\varepsilon_i$}
\put(24,10){$\varepsilon_j$}
\end{picture}}
\end{picture}
\vspace{0.2cm}
\end{equation}
\end{center}

For any $kQ'$-module $M$ with dimension vector $n\widehat{S'_i}+\widehat{S'_j}$, it can be decomposed to be $(S'_i)^{\oplus u_M}\oplus N$ with $N$ indecomposable. We denote
\begin{align}
\mathcal{I}_{k} &=\{[M]\in\Iso(\mod(kQ'))\mid \exists N\subseteq M \text{ s.t. }N\cong S'_j, M/N\cong (S'_i)^{\oplus k }\},
\\
p(a,r,s,t) &=-s(a+t)+2ra+(u_M-t+2s-r)(t-r) +(s-r)^2
 \label{eq:p}
\\
&\qquad\qquad +{s-r \choose 2} + (t-r)^2 +{t-r \choose 2} +r(s+t) -{r+1 \choose 2}+1.
\notag
\end{align}

\begin{lemma}[\text{\cite[Proposition 7.3]{LW20a}}]
  \label{lem:SSSM}
For any $s,t\in\N$, we have
\begin{align}
&[s S'_i]*[S'_j]*[t S'_i]
  \label{eq:SSS}\\
&=\sum_{r=0}^{\min\{ s,t\}}\sum_{[M]\in \mathcal{I}_{s+t-2r}}
\sqq^{p(a,r,s,t)} (\sqq -\sqq^{-1})^{s+t-r+1} \frac{[s]_\sqq^{!} [t]_\sqq^{!}}{[r]_\sqq^{!}}
\qbinom{u_M}{t-r}_\sqq \frac{[M]}{|\aut(M)|}*[\E'_i]^r.
\notag
\end{align}
\end{lemma}

To prove \eqref{eqn:braidsplitev}, we shall compute the RHS of  \eqref{eqn:braidsplitev} in the rank 2 $\imath$quiver algebra associated to \eqref{diag:split KM} above, with the help of  Lemma~\ref{lem:iDPev} and Lemma~\ref{lem:SSSM}.

\subsection{Computation of $[S'_i]_{\ev}^{(r)}*[S'_j] *[S'_i]_{\ov{a}}^{(s)}$}
\label{subsec:SSS}

Let us first compute $[S'_i]_{\ev}^{(r)}*[S'_j] *[S'_i]_{\ov{a}}^{(s)}$, depending on the parity of $r$.
\subsubsection{$r$ is even}

%{\bf Case (I)}. Assume $r$ is even.
For any $s\geq0$ such that $r+s+2t=a$ with $t\geq0$, we have by Lemma~\ref{lem:iDPev} (noting $\ov{s} =\ov{a}$) and Lemma~\ref{lem:SSSM} that
\begin{align*}
[S'_i]_{\ev}^{(r)}*[S'_j] *[S'_i]_{\ov{a}}^{(s)}
&=\sum_{k=0}^{\frac{r}{2}}\frac{\sqq^{k(k-1)-\binom{r-2k}{2}} \cdot (\sqq-\sqq^{-1})^k }{[r-2k]_\sqq^{!}[2k]_\sqq^{!!}} [(r-2k)S'_i]*[\E'_i]^k*[S'_j]*
\\
&\times \sum_{m=0}^{\lfloor\frac{s}{2}\rfloor} \frac{\sqq^{m(m-1)-\binom{s-2m}{2}}\cdot (\sqq-\sqq^{-1})^{m}}{[s-2m]_\sqq^{!} [2m]_\sqq^{!!}} [(s-2m)S'_i]*[\E'_i]^m
\\
&= \sum_{k=0}^{\frac{r}{2}}\sum_{m=0}^{\lfloor\frac{s}{2}\rfloor} \frac{\sqq^{k(k-1)+m(m-1)-\binom{r-2k}{2} -\binom{s-2m}{2}}\cdot (\sqq-\sqq^{-1})^{k+m} }{[r-2k]_\sqq![s-2m]_\sqq![2k]_\sqq^{!!}[2m]_\sqq^{!!}}
 \\
&\qquad\qquad \times [(r-2k)S'_i]*[S'_j]*[(s-2m)S'_i]*[\E'_i]^{k+m}
\\
&= \sum_{k=0}^{\frac{r}{2}}\sum_{m=0}^{\lfloor\frac{s}{2}\rfloor} \sum_{n=0}^{\min\{r-2k,s-2m\}} \sum_{[M]\in\mathcal{I}_{r+s-2k-2m-2n}}
\frac{\sqq^{k(k-1)+m(m-1)-\binom{r-2k}{2} -\binom{s-2m}{2}} }{[r-2k]_\sqq![s-2m]_\sqq![2k]_\sqq^{!!}[2m]_\sqq^{!!}}
 \\
&\qquad\qquad \times \sqq^{p(a,n,r-2k,s -2m)}
(\sqq -\sqq^{-1})^{r+s -k -m -n+1 } \frac{[r-2k]_\sqq^{!} [s -2m]_\sqq^{!}}{[n]_\sqq^{!}}
\\
&\qquad\qquad \times \qbinom{u_M}{s -2m-n}_\sqq \frac{[M]}{|\aut(M)|}*[\E'_i]^{n+k+m}.
\end{align*}
This can be simplified to be
\begin{align*}
[S'_i]_{\ev}^{(r)}*[S'_j] *[S'_i]_{\ov{a}}^{(s)}
&= \sum_{k=0}^{\frac{r}{2}}\sum_{m=0}^{\lfloor\frac{s}{2}\rfloor} \sum_{n=0}^{r-2k} \sum_{[M]\in\mathcal{I}_{r+s-2k-2m-2n}}
\\
&\qquad  \frac{\sqq^{\cz(a,r,s,k,m,n)}
(\sqq-\sqq^{-1})^{r+s-k-m-n+1} }{[n]_\sqq^{!} [2k]_\sqq^{!!}[2m]_\sqq^{!!}}
\qbinom{u_M}{s -2m-n}_\sqq \frac{[M] *[\E'_i]^{n+k+m}}{|\aut(M)|}
\end{align*}
where
\begin{equation}
\label{eq:z}
\cz(a,r,s,k,m,n) ={k(k-1)+m(m-1)-\binom{r-2k}{2} -\binom{s-2m}{2} +p(a,n,r-2k, s -2m)},
\end{equation}
and
\begin{align*}
p& (a,n,r-2k, s -2m)
\\
&= -(r-2k)(a+s-2n)+2a n+(u_M-s+2m+2r-4k-n)(s-2m-n)\\
&\quad +(r-2k-n)^2+{ r-2k-n \choose 2} + (s-2m-n)^2+ { s-2m-n \choose 2} \\
&\quad +n(r+s-2m-2k) -{ n+1 \choose 2}+1.
\end{align*}

\subsubsection{$r$ is odd}

%{\bf Case (II)}. Assume $r$ is odd.
For any $s\geq0$ such that $r+s+2t=a$ with $t\geq0$, we have by Lemma~\ref{lem:iDPev} (noting $\ov{s} \not=\ov{a}$) and Lemma~\ref{lem:SSSM}
\begin{align*}
[S'_i]_{\ev}^{(r)}*[S'_j] *[S'_i]_{\ov{a}}^{(s)}
=&\sum_{k=0}^{\lfloor\frac{r}{2}\rfloor}\frac{\sqq^{k(k+1)-\binom{r-2k}{2}} \cdot (\sqq-\sqq^{-1})^k }{[r-2k]_\sqq![2k]_\sqq^{!!}} [(r-2k)S'_i]*[\E'_i]^k*[S'_j]*\\
&\times \sum_{m=0}^{\lfloor\frac{s}{2}\rfloor} \frac{\sqq^{m(m+1)-\binom{s-2m}{2}}\cdot (\sqq-\sqq^{-1})^{m}}{[s-2m]_\sqq^{!}[2m]_\sqq^{!!}} [(s-2m)S'_i]*[\E'_i]^m\\
=&\sum_{k=0}^{\lfloor\frac{r}{2}\rfloor}\sum_{m=0}^{\lfloor\frac{s}{2}\rfloor} \sum_{n=0}^{\min\{r-2k,s-2m\}} \sum_{[M]\in\mathcal{I}_{r+s-2k-2m-2n}}
 \frac{\sqq^{k(k+1)+m(m+1)-\binom{r-2k}{2} -\binom{s-2m}{2}}  }{[r-2k]_\sqq![s-2m]_\sqq![2k]_\sqq^{!!}[2m]_\sqq^{!!}}\\
&\quad  \times \sqq^{p(a,n,r-2k, s -2m)}
(\sqq -\sqq^{-1})^{r+s -k -m -n+1 } \frac{[r-2k]_\sqq^{!} [s -2m]_\sqq^{!}}{[n]_\sqq^{!}}
\\ &
\quad  \times \qbinom{u_M}{s -2m-n}_\sqq \frac{[M]}{|\aut(M)|}*[\E'_i]^{n+k+m}.
\end{align*}
This can be simplified to be
\begin{align*}
&[S'_i]_{\ev}^{(r)}*[S'_j] *[S'_i]_{\ov{a}}^{(s)}
= \sum_{k=0}^{\lfloor\frac{r}{2}\rfloor}\sum_{m=0}^{\lfloor\frac{s}{2}\rfloor} \sum_{n=0}^{r-2k} \sum_{[M]\in\mathcal{I}_{r+s-2k-2m-2n}}
  \\
& \qquad\qquad \frac{\sqq^{\cz(a,r,s,k,m,n)+2k+2m}
 (\sqq-\sqq^{-1})^{r+s -k -m -n+1} }{[n]_\sqq^{!} [2k]_\sqq^{!!}[2m]_\sqq^{!!}}
\qbinom{u_M}{s -2m-n}_\sqq \frac{[M] *[\E'_i]^{n+k+m}}{|\aut(M)|}.
\end{align*}

\subsection{Computation of RHS\eqref{eqn:braidsplitev} }

Summing up the computations of $[S'_i]_{\ev}^{(r)}*[S'_j] *[S'_i]_{\ov{a}}^{(s)}$ in \S\ref{subsec:SSS}, we obtain
\begin{align}
\text{RHS} ~\eqref{eqn:braidsplitev}
   \label{eq:long}
 &= \sum_{r=0,2\mid r}^{a}\sum_{k=0}^{\frac{r}{2}} \sum_{m=0}^{ \lfloor\frac{a-r}{2}\rfloor } \sum_{n=0}^{r-2k} \sum_{[M]\in\mathcal{I}_{a-2k-2m-2n}}
 (-1)^a (\sqq-\sqq^{-1})^{-k-m-n+1}
 \\
& \qquad \times \frac{\sqq^{r+\cz(a,r,a-r,k,m,n)-a} }{[n]_\sqq^{!} [2k]_\sqq^{!!}[2m]_\sqq^{!!}}
\qbinom{u_M}{a-r -2m-n}_\sqq \frac{[M] *[\E'_i]^{n+k+m}}{|\aut(M)|}
 \notag \\
&-\sum_{r=0,2\nmid r}^{a} \sum_{k=0}^{\lfloor\frac{r}{2}\rfloor}\sum_{m=0}^{\lfloor\frac{a-r}{2}\rfloor} \sum_{n=0}^{r-2k} \sum_{[M]\in\mathcal{I}_{a-2k-2m-2n}}
 (-1)^a  (\sqq-\sqq^{-1})^{-k-m-n+1}
\notag \\
& \qquad \times \frac{\sqq^{r+\cz(a,r,a-r,k,m,n)+2k+2m-a} }{[n]_\sqq^{!} [2k]_\sqq^{!!}[2m]_\sqq^{!!}}
\qbinom{u_M}{a-r -2m-n}_\sqq \frac{[M] *[\E'_i]^{n+k+m}}{|\aut(M)|}
 \notag \\
&+\sum_{t\geq1}\sum_{r=0,2\mid r}^{a-2t}\sum_{k=0}^{\frac{r}{2}}\sum_{m=0}^{\lfloor \frac{a-r}{2}\rfloor-t} \sum_{n=0}^{r-2k} \sum_{[M]\in\mathcal{I}_{a-2t-2k-2m-2n}}
 (-1)^a (\sqq-\sqq^{-1})^{-k-m-n+1}
\notag \\
& \qquad \times \frac{\sqq^{r+ \cz(a,r,a-2t-r,k,m,n)  -a+2t} }{[n]_\sqq^{!} [2k]_\sqq^{!!}[2m]_\sqq^{!!}}
\qbinom{u_M}{a-2t-r -2m-n}_\sqq \frac{[M] *[\E'_i]^{n+k+m+t}}{|\aut(M)|}.
\notag
\end{align}

Fix
\[
d =t+k+m+n
\]
 and fix $[M]   \in\mathcal{I}_{a-2d}$. For $u_M=0=d$, $M$ is indecomposable which is isomorphic to $F_i^+(S_j)$ by Lemma~ \ref{lem:reflecting dimen}. In this case, the coefficient of $[M]$  of the RHS of \eqref{eq:long} is $1$ by noting that $|\Aut(M)|=q-1$ in this case.  By \eqref{eqn:reflection 1},  we have reduced the proof of  \eqref{eqn:braidsplitev} to proving the coefficient of $\frac{[M] *[\E'_i]^{d}}{|\aut(M)|}$ of the RHS of \eqref{eq:long} is 0, for any given $[M] \in \mathcal I_{a-2d}$ such that not both $d$ and $u_M$ are 0.

We have $0\le n=d-k-m-t \leq r-2k$. Set $u=u_M$, and
\begin{align}
\label{eq:defA}
A(a,d,u)
&:= \sum_{t\geq0}\sum_{r=0,2\mid r}^{a-2t}\sum_{k=0}^{\frac{r}{2}}\sum_{m=0}^{\lfloor\frac{a-r}{2}\rfloor-t}
\delta\{0\le n \le r-2k\}
\\
& \quad \times \frac{\sqq^{r+\cz(a,r,a-2t-r,k,m,n)-a+2t}
(\sqq-\sqq^{-1})^{-k-m-n+1} }{[n]_\sqq^{!} [2k]_\sqq^{!!}[2m]_\sqq^{!!}}
 \qbinom{u}{a-2t-r -2m-n}_\sqq \notag
  \\
&\quad -\sum_{r=0,2\nmid r}^{a} \sum_{k=0}^{\lfloor\frac{r}{2}\rfloor}\sum_{m=0}^{\lfloor\frac{a-r}{2}\rfloor}
  \delta\{0\le n \le r-2k\}
\notag \\
& \quad \times \frac{\sqq^{r+\cz(a,r,a-r,k,m,n)+2k+2m-a}
 (\sqq-\sqq^{-1})^{ -k -m -n+1} }{[n]_\sqq^{!} [2k]_\sqq^{!!}[2m]_\sqq^{!!}}
\qbinom{u}{a -r-2m-n}_\sqq,
\notag
\end{align}
 %any $a\geq0$, $0\leq d\leq \frac{a}{2}$, $0\leq u \leq a-2d$, %$d$ and $u_M$ not both zero,
where $\cz(\cdot, \cdot, \cdot, \cdot, \cdot, \cdot)$ is given in \eqref{eq:z}; also see \eqref{eq:p} for $p(\cdot, \cdot, \cdot, \cdot)$. Here $\delta\{X\}=1$ if the statement $X$ holds and $=0$ if $X$ is false.

Then the coefficient of $\frac{[M] *[\E'_i]^{d}}{|\aut(M)|}$ of the RHS of \eqref{eq:long} is $(-1)^a A(a,d,u_M)$. Summarizing, we have reached the following reduction toward the proof of \eqref{eqn:braidsplitev}.

\begin{proposition}
 \label{prop:GammaA}
The identity \eqref{eqn:braidsplitev} is equivalent to the following identity
\begin{equation}
\label{eq:A=0}
A(a,d,u) =0,
\end{equation}
for non-negative integers $a,d,u$ subject to the constraints:
\begin{equation}  \label{constraints}
0\leq d\leq \frac{a}{2}, \quad 0\leq u \leq a-2d, \quad  d \text{ and $u$ not both zero}.
\end{equation}
\end{proposition}
\subsection{Reduction for the identity \eqref{eq:A=0}  }
\label{subsec:A=0}

We shall denote the 2 summands in $A =A(a,d,u)$ in \eqref{eq:defA} as $A_0, A_1$, and thus
\[
A =A_0 -A_1.
\]
We shall denote
\begin{align}
w &=r-2k-n. %,  \qquad a=a.
\end{align}
Set
\begin{align}
d= k+m+n +t
\end{align}
in the $A_0$ side, and $d =k+m+n$ in the $A_1$ side.  Then
\[
a -2t -r -2m -n=a -2d -w
\]
in the $A_0$ side, and
\[
a  -r -2m -n=a -2d -w
\]
in the $A_1$ side.
(Below we focus on the $A_0$ side, and the corresponding formulas for $A_1$ side can be obtained by setting $t=0$.)
Note that
\begin{align*}
&p(a,n,r-2k, a-2t-r -2m)
\\
=&{\small -(r-2k)(2a-2t-r-2n)+2a n+(u-a+2t+3r+2m-4k-n)(a-2t-r-2m-n)  } \\
& +(r-2k-n)^2+{ r-2k-n \choose 2} + (a-2t-r-2m-n)^2+ { a-2t-r-2m-n \choose 2}  \\
&+n(a-2t-2m-2k) -{ n+1 \choose 2}+1.
\end{align*}
One can write
$p(a, n, r-2k, a-2t-r -2m) =p(a, n, n+w, a-2d-w +n) =p' +p''$, where $p'$ depends on $n$, while $p''$ does not. Then one can show that
\[
p' =an + \binom{n}{2},
\]
which is independent of $t$.
Note
\begin{align}
\label{eq:zz}
\cz: =\cz(a,r,a-2t-r,k,m,n)
&= k(k-1)+m(m-1)-\binom{r-2k}{2} -\binom{a-2t-r-2m}{2} \\
&\quad +p(a,n,r-2k, a-2t-r -2m).
\notag
\end{align}
We rewrite $\cz$ in \eqref{eq:zz} as $\cz =\cz_1 +\cz_2$, where $\cz_2$ depends only on $a, d, w$ but do not depend on $k, m, n, t$. One can show that on the $A_0$ side
\begin{align}
\label{eq:zz2}
\cz_1 = 2 \binom{t+n+1}{2} -\binom{n}{2} -2dt -2km,
\end{align}
and thus on the $A_1$ side, $\cz_1 =2 \binom{n+1}{2} -\binom{n}{2}  -2km.$

Will the above preparations, we return to \eqref{eq:defA}. Noting $r=w +n +2k$, we rewrite
\[
r +\cz +2t = x_1 +x_2,
\]
where $x_2$ depends only on $a, d, w$, and
\begin{align}
 \label{eq:x1}
x_1 = t^2 -2dt +t +2nt + \binom{n+1}{2}  -2km -2m.
\end{align}
Setting $t =0$ in \eqref{eq:x1}, we obtain on the $A_1$ side that
\[
r+\cz +2k +2m = x_1' +x_2,
\]
where
\begin{align}
x_1' =\binom{n+1}{2} -2km +2k.
\end{align}
Note $r= 2k +n +w \equiv n+w \pmod 2$.
Therefore, proving the identity \eqref{eq:A=0} for $d>0$ amounts to showing that the coefficients of $\qbinom{u}{a-2d-w}_\sqq$ from \eqref{eq:defA} are 0, for $d>0$ (and fixed $u, a, w$); that is,  for $d>0$,
\begin{align}
 \label{eq:coeff}
\sum_{\stackrel{w+n \text{ even}}{t+k+m+n =d}}
\frac{\sqq^{t^2 -2dt +t  +2nt +\binom{n+1}{2} -2km -2m}}{[n]^!_\sqq [2k]^{!!}_\sqq [2m]^{!!}_\sqq} (\sqq -\sqq^{-1})^{t}
- \sum_{\stackrel{w+n \text{ odd}}{k+m+n =d}}
\frac{\sqq^{\binom{n+1}{2} -2km +2k}}{[n]^!_\sqq [2k]^{!!}_\sqq [2m]^{!!}_\sqq} =0.
\end{align}
It is understood here and below that all $t, k, m, n$ in the summations above are in $\N$.
Summarizing, we have obtained the following reduction toward the proof of \eqref{eq:A=0}.

\begin{proposition}
 \label{prop:ACD}
The identity \eqref{eq:A=0} that $A(a,d,u) =0$, for $d>0$, is equivalent to the $v$-binomial identity \eqref{eq:coeff}.
\end{proposition}
We shall prove the identity \eqref{eq:coeff} in Section~\ref{sec:identities}.

\subsection{The identity \eqref{eq:A=0} for $d=0$}
\label{subsec:contribution d=0}
%%\vspace{2mm}

By \eqref{constraints}, we must have $u>0$ when $d=0$.

In this case, we have $k=m=n=t=0$, and a direct computation shows that the power $z$ can be simplified to be $z = au -uw +1 $, and thus  $r+z -a+2t = r+z -a +2k +2m = au +1 -a + (1-u)w$. Recall $A(a, 0, u)$ from \eqref{eq:defA}. Then, for $0<u \leq a$, we have
\begin{align}
   \label{eq:d=0}
A(a, 0, u)
&= \sqq^{au +1-a} (\sqq-\sqq^{-1}) \sum_{w \ge 0}  (-1)^w \sqq^{(1-u)w} \qbinom{u}{a -w}_\sqq
\\
&\stackrel{(i)}{=} (-1)^{a} \sqq (\sqq-\sqq^{-1}) \sum_{x= 0}^u (-1)^x \sqq^{(u-1)x}  \qbinom{u}{x}_\sqq
\stackrel{(ii)}{=}  0,
\notag
\end{align}
 where we have changed variables $x=a-w$ in the identity (i) and replaced the upper bound of the summation for $x$ from $a$ to $u$ (thanks to $u \le a$ and $\qbinom{u}{x}_\sqq=0$ for $x>u$); the $v$-binomial identity \eqref{1.3.1} (with $z=-1$) was used in (ii) above.

\subsection{Proof of Theorem \ref{thm:Braidsplit} }

The identity \eqref{eq:A=0} follows now by combining the identity \eqref{eq:coeff} for $d>0$ (which is proved in Section~\ref{sec:identities}) and the identity \eqref{eq:d=0}, thanks to  Proposition~\ref{prop:ACD}. Then by Proposition~\ref{prop:GammaA}, the formula \eqref{eqn:braidsplitev} follows.

The formula \eqref{eqn:braidsplitodd} is proved similarly by a reduction to an analogous identity \eqref{eq:B=0}; see Appendix~\ref{app:braidsplitodd}.
%The proof of the formula \eqref{eqn:braidsplitodd} is very similar and hence will be skipped (the detail can be found in Appendix~{\bf B} in arXiv v1).

For the case $\btau j=j$, the formula \eqref{eqn:Braidsplit} follows from \eqref{eqn:braidsplitev}--\eqref{eqn:braidsplitodd}.

For the case $\btau j\neq j$, it is enough to consider a general rank 2 $\imath$quiver $(Q,\btau)$ such that $i=\btau i$ and $\btau j \not=j$, as shown in the left figure of \eqref{diag:qsplit KM} below. Then the $\imath$quiver $\ov{Q'}$ of $\bs_i\Lambda^\imath$ is shown in the right figure of \eqref{diag:qsplit KM}.
\begin{center}\setlength{\unitlength}{1mm}
\vspace{.1cm}
 \begin{equation}
 \label{diag:qsplit KM}
  \begin{picture}(50,10)(0,-10)
  \put(-20,0){
 \begin{picture}(50,10)
\put(0,-2){$j$}
\put(20,-2){$\btau j$}

  \put(-15,-10){$\ov{Q}=$}

\put(7,-12){\vector(1,-2){3.5}}
\put(5.5,-12){$a$}
\put(6,-10){\line(-1,2){3.5}}

\put(3,1){\line(1,0){7.5}}
\put(10.5,0){$r$}
\put(12,1){\vector(1,0){7.5}}

\put(10.5,-1.5){\vector(-1,0){7.5}}
\put(10.5,-2.5){$r$}
\put(19.5,-1.5){\line(-1,0){7.5}}

\put(15.5,-12){\vector(-1,-2){3.5}}
\put(15,-12){$a$}
\put(16.5,-10){\line(1,2){3.5}}

\put(10,-22){$i$}
\color{purple}
\qbezier(3,2)(11,5)(19,2)
\put(19,2){\vector(2,-1){0.2}}

\qbezier(3,-2)(11,-5)(19,-2)

\put(3,-2){\vector(-2,1){0.2}}
\put(10,4){$^{\varepsilon_j}$}
\put(10,-5){$_{\varepsilon_{\btau j}}$}
\put(10,-28){$_{\varepsilon_i}$}
\begin{picture}(50,23)(-10,19)
\color{purple}
\qbezier(-1,-1)(-3,-3)(-2,-5.5)
\qbezier(-2,-5.5)(1,-9)(4,-5.5)
\qbezier(4,-5.5)(5,-3)(3,-1)
\put(3.1,-1.4){\vector(-1,1){0.3}}
\end{picture}
\end{picture}}
 \put(50,0){
 \begin{picture}(50,10)
\put(0,-2){$j$}
\put(20,-2){$\btau j$}

  \put(-15,-10){$\ov{Q'}=$}

\put(10,-18){\line(-1,2){3}}
\put(5.5,-12){$a$}
\put(6,-10){\vector(-1,2){3.5}}

\put(3,1){\line(1,0){7.5}}
\put(10.5,0){$r$}
\put(12,1){\vector(1,0){7.5}}

\put(10.5,-1.5){\vector(-1,0){7.5}}
\put(10.5,-2.5){$r$}
\put(19.5,-1.5){\line(-1,0){7.5}}

\put(12.5,-18){\line(1,2){3}}
\put(15,-12){$a$}
\put(16.5,-10){\vector(1,2){3.5}}

\put(10,-22){$i$}
\color{purple}
\qbezier(3,2)(11,5)(19,2)
\put(19,2){\vector(2,-1){0.2}}

\qbezier(3,-2)(11,-5)(19,-2)

\put(3,-2){\vector(-2,1){0.2}}
\put(10,4){$^{\varepsilon_j}$}
\put(10,-5){$_{\varepsilon_{\btau j}}$}
\put(10,-28){$_{\varepsilon_i}$}
\begin{picture}(50,23)(-10,19)
\color{purple}
\qbezier(-1,-1)(-3,-3)(-2,-5.5)
\qbezier(-2,-5.5)(1,-9)(4,-5.5)
\qbezier(4,-5.5)(5,-3)(3,-1)
\put(3.1,-1.4){\vector(-1,1){0.3}}
\end{picture}
\end{picture}}
\end{picture}
\vspace{1.8cm}
\end{equation}
\end{center}
By the same argument as in \cite[Proposition 9.4]{LW20a}, the computations involved in proving the formula \eqref{eqn:Braidsplit} with $j \not=\btau j$ are the same as for the $\imath$quiver in \eqref{diag:split KM} with $j =\btau j$.

The proof for \eqref{eqn:Braidsplit} and thus Theorem \ref{thm:Braidsplit} is completed (modulo the proof of the identity \eqref{eq:coeff}, which will be given in Section~\ref{sec:identities}).

%%%%%%%
%%%%%%%
\section{Several quantum binomial identities}
 \label{sec:identities}

 In this section, we shall prove the identity \eqref{eq:coeff}, by establishing several additional $v$-binomial identities (which seem to be new). We shall switch notations from $\sqq$ to a general parameter $v$ in this section.

\subsection{Reformulating \eqref{eq:coeff}}

With $\sqq$ replaced by $v$, we shall denote the first (and respectively, second) summand in \eqref{eq:coeff} as $D_{\bar w}$ (and respectively, $C_{\bar w}$), for $\bar w\in \{0,1\}$ with $w\equiv \bar w \pmod 2$. That is,
\begin{align}
 \label{eq:coeff0}
D_0 &=
\sum_{\stackrel{n \text{ even}}{t+k+m+n =d}}
\frac{v^{t^2 -2dt +t +2nt +\binom{n+1}{2} -2km -2m}}{[n]^{!}_v [2k]^{!!}_v [2m]^{!!}_v} (v -v^{-1})^{t},
\\
C_0 &= \sum_{\stackrel{n \text{ odd}}{k+m+n =d}}
\frac{v^{\binom{n+1}{2} -2km  +2k}}{[n]^{!}_v [2k]^{!!}_v [2m]^{!!}_v},
\\
 \label{eq:coeff1}
D_1 &=
\sum_{\stackrel{n \text{ odd}}{t+k+m+n =d}}
\frac{v^{t^2 -2dt +t +2nt +\binom{n+1}{2} -2km -2m}}{[n]^{!}_v [2k]^{!!}_v [2m]^{!!}_v} (v -v^{-1})^{t},
\\
 \label{eq:coeff2}
C_1 &= \sum_{\stackrel{n \text{ even}}{k+m+n =d}}
\frac{v^{\binom{n+1}{2} -2km  +2k}}{[n]^{!}_v [2k]^{!!}_v [2m]^{!!}_v}.
\end{align}
Proving the identity \eqref{eq:coeff} amounts to showing that
\begin{align}  \label{eq:DC0}
D_0 -C_0 =0,
\qquad
D_1 -C_1 =0.
\end{align}

\subsection{Identities for $D_0 - D_1$ and $C_0 - C_1$}

To establish \eqref{eq:DC0}, we shall prove new $v$-binomial identities which lead to the stronger statement that
$D_0=D_1 =C_0 =C_1$ in Theorem~\ref{thm:DDCC}.

Since an identity involving summations over $n$ with a fixed parity (such as $D_i, C_i$) seems hard to prove directly, our strategy is to proceed by proving identities regarding the combinations $D_0 \pm D_1$ or $C_0 \pm C_1$.

By definition we have
\begin{align}
 \label{eq:B01}
D_0 \pm D_1 &=
\sum_{t+k+m+n =d} (\pm 1)^n
\frac{v^{t^2 -2dt +t +2nt +\binom{n+1}{2} -2km -2m}}{[n]^{!}_v [2k]^{!!}_v [2m]^{!!}_v} (v -v^{-1})^{t},
\\
 \label{eq:C01}
C_1 \pm C_0 &=
\sum_{k+m+n =d}  (\pm 1)^n
\frac{v^{\binom{n+1}{2}} v^{-2km  +2k}}{[n]^{!}_v [2k]^{!!}_v [2m]^{!!}_v}.
\end{align}

\begin{lemma} [\text{\cite[Lemmas 8.2, 8.3]{LW20a}}]
 \label{lem:binomial}
For $p \ge 0$ and $d \ge 1$,  the following identities hold:
\begin{align}
[p]^{!}_v\sum_{\stackrel{k,m \in \N}{k+m =p}} \frac{v^{-2km +2m} }{[2k]^{!!}_v [2m]^{!!}_v} &= v^{\frac{p(3-p)}2},
  \label{eq:km1}
\\
\label{eq:kmrd}
\sum_{\stackrel{k,m,r \in \N}{k+m+r =d}} (-1)^r \frac{v^{{r+1 \choose 2} -2(k-1)m}}{[r]^{!}_v [2k]^{!!}_v [2m]^{!!}_v} & =0.
\end{align}
\end{lemma}

In the notation of $C_0, C_1$, the second formula in Lemma~\ref{lem:binomial} can read as
\begin{align}
 \label{eq:CC}
  C_1 -C_0=0.
\end{align}

%The identity $D_0-D_1=0$ is reformulated as follows.

\begin{proposition}
  \label{prop:DD}
For $d\ge 1$, the identity $D_0-D_1=0$ holds; that is,
\begin{align}
 \label{eq:B-=0}
\sum_{\stackrel{t,k,m,n \in \N}{t+k+m+n =d}} (-1)^n
\frac{v^{t^2 -2dt +t +2nt +\binom{n+1}{2} -2km -2m}}{[n]^{!} [2k]^{!!} [2m]^{!!}} (v -v^{-1})^{t}=0.
\end{align}
\end{proposition}

\begin{proof}
Recall the following standard $v$-binomial identity \cite[1.3.1(c)]{Lus93}:
\begin{align} \label{1.3.1}
\sum_{n=0}^{d} v^{n(d-1)} \qbinom{d}{n} z^{n} &=  \prod_{j=0}^{d-1} (1 +v^{2j} z).
\end{align}

By \eqref{eq:km1} (with $k,m$ switched), we have
\begin{align}
  \label{eq:km6}
\sum_{\stackrel{k,m \in \N}{k+m =p}} \frac{v^{-2km +2k} }{[2k]^{!!} [2m]^{!!}} &=  \frac{v^{\frac{p(3-p)}2}}{[p]^{!}}.
\end{align}

For any fixed $t$ with $0\le t \le d$, using the identity \eqref{eq:km6}, we compute
\begin{align}
 \label{eq:tkm}
\sum_{\stackrel{k,m,n \in \N}{k+m+n =d-t} }
& (\pm 1)^n
 \frac{v^{t^2 -2dt +t  +2nt +\binom{n+1}{2} -2km -2m}}{[n]^{!} [2k]^{!!} [2m]^{!!}}
 \\
&=  \sum_{\stackrel{p,n \in \N}{p+n =d-t} }
 (\pm 1)^n  \frac{v^{t^2 -2dt +t +2nt +\binom{n+1}{2} - \frac{p(p+1)}2}}{[n]^{!} [p]^{!}}
  \notag \\
&=\frac{v^{(1-d)t  -\frac{(d-t)(d+t+1)}2} }{[d-t]^{!}}
\sum_{n=0}^{d-t} v^{n(d-t-1)} \qbinom{d-t}{n} (\pm v^{2t+2})^n
  \notag \\
&\stackrel{(*)}{=} \frac{ v^{(1-d)t  -\frac{(d-t)(d+t+1)}2} }{[d-t]^{!}} \prod_{j=t+1}^{d} (1 \pm v^{2j})
  \notag \\
&= \frac{ v^{(1-d)t} }{[d-t]^{!}} \prod_{j=t+1}^{d} (v^{-j} \pm v^{j}),
  \notag
\end{align}
where the second last equality $(*)$ follows by \eqref{1.3.1} with $z=\pm v^{2t+2}$.

Using now the minus sign version of \eqref{eq:tkm}, we compute
\begin{align*}
\sum_{\stackrel{t,k,m,n \in \N}{t+k+m+n =d}}
& (-1)^n
\frac{v^{t^2 -2dt +t +2nt +\binom{n+1}{2} -2km -2m}}{[n]^{!} [2k]^{!!} [2m]^{!!}} (v -v^{-1})^{t}
\\
&= \sum_{t=0}^{d} (-1)^{d-t}  \frac{ v^{(1-d)t} }{[d-t]^{!}} \prod_{j=t+1}^{d} (v^j -v^{-j})   (v -v^{-1})^{t}
\\
&= \sum_{t=0}^{d} (-1)^{t}   v^{(1-d)t}  \qbinom{d}{t} \cdot (-1)^d  (v -v^{-1})^{d}
=0,
\end{align*}
where for the last equality we have used \eqref{1.3.1} with $z=-1$ (and a switch $v \leftrightarrow v^{-1}$).
\end{proof}

\subsection{Identities for $D_1+ D_0$ and $C_1 + C_0$}

We shall also prove the following $v$-binomial identities for $C_1 +C_0$ and $D_0 +D_1$; cf. \eqref{eq:B01}--\eqref{eq:C01}. %, from which the identity $D_0 +D_1 =C_0 +C_1$ follows.

\begin{proposition}
  \label{prop:DD=CC}
For $d\ge 1$, the following identities hold:
\begin{align}
 \label{eq:C=0}
\sum_{\stackrel{k,m,n \in \N}{k+m+n =d} }
& \frac{v^{\binom{n+1}{2} -2km  +2k}}{[n]^! [2k]^{!!} [2m]^{!!}}
= \frac{2v^d (v +v^{-1}) (v^2 +v^{-2}) \cdots (v^{d-1} +v^{1-d})}{[d]!},
\\
 \label{eq:B+=0}
\sum_{\stackrel{t,k,m,n \in \N}{t+k+m+n =d} }
& \frac{v^{t^2 -2dt +t +2nt +\binom{n+1}{2} -2km -2m}}{[n]^! [2k]^{!!} [2m]^{!!}} (v -v^{-1})^{t}
\\
&= \frac{2v^d (v +v^{-1}) (v^2 +v^{-2}) \cdots (v^{d-1} +v^{1-d})}{[d]!}.
\notag
\end{align}
In particular, the identity $D_0 +D_1 =C_0 +C_1$ holds.
\end{proposition}

\begin{proof}
Using the identity \eqref{eq:km6}, we obtain
\begin{align*}
\sum_{\stackrel{k,m,n \in \N}{k+m+n =d} }
 \frac{v^{\binom{n+1}{2} -2km  +2k}}{[n]^! [2k]^{!!} [2m]^{!!}}
&= \sum_{{n+p =d} } v^{\binom{n+1}{2}} v^{\frac{p(3-p)}2} \frac{1}{[n]^! [p]^!}
\\
&= v^{d-\binom{d}{2}} \frac{1}{[d]^!}  \sum_{n=0}^{d} v^{n(d-1)} \qbinom{d}{n}
= v^{d-\binom{d}{2}} \prod_{j=0}^{d-1} (1 +v^{2j}),
\end{align*}
where the last equality follows by \eqref{1.3.1} with $z=1$.
This proves the identity \eqref{eq:C=0}.

Summing the plus sign version of the identity \eqref{eq:tkm} over $0\le t \le d$, we have reduced the proof of \eqref{eq:B+=0} to the following identity:
\begin{align}
 \label{eq:B3}
\sum_{t=0}^d
& \frac{v^{(1-d)t} (v^{t+1} +v^{-t-1}) (v^{t+2} +v^{-t-2}) \cdots (v^{d} +v^{-d})}{[d-t]!} (v -v^{-1})^{t}
\\
&= \frac{2v^d (v +v^{-1}) (v^2 +v^{-2}) \cdots (v^{d-1} +v^{1-d})}{[d]!}.
\notag
\end{align}
The identity \eqref{eq:B3} does not seem to be easy to prove directly; it follows as a special case when $k=d$ of a more general identity in Proposition~\ref{prop:B4} below (by a change of variable $t=d-s$).
\end{proof}

\begin{proposition}
  \label{prop:B4}
Let $d\ge 1$. Then the following identity holds, for $0\le k \le d$:
\begin{align}
 \label{eq:B4}
\sum_{s=0}^k
& \frac{v^{(1-d)(d-s)} (v^{d-s+1} +v^{s-d-1}) (v^{d-s+2} +v^{s-d-2}) \cdots (v^{d} +v^{-d})}{[s]!} (v -v^{-1})^{d-s}
\\
&= v^{d(k+1-d)} (v -v^{-1})^{d-k}  \frac{(v^{d-k} +v^{k-d}) (v^{d-k+1} +v^{k-d-1}) \cdots (v^{d-1} +v^{1-d}) }{[k]!}.
\notag
\end{align}
\end{proposition}

\begin{proof}
We prove by induction on $k$. The base case for $k=0$ is clear.
By the inductive assumption, we have
\begin{align*}
\text{LHS}\eqref{eq:B4}
&= v^{d(k-d)} (v -v^{-1})^{d-k+1}  \frac{(v^{d-k+1} +v^{k-1-d}) (v^{d-k+2} +v^{k-d-2}) \cdots (v^{d-1} +v^{1-d}) }{[k-1]!}
\\
&\quad+\frac{v^{(1-d)(d-k)} (v^{d-k+1} +v^{k-d-1}) (v^{d-k+2} +v^{k-d-2}) \cdots (v^{d} +v^{-d})}{[k]!} (v -v^{-1})^{d-k}
\\
&= v^{d(k-d)} (v-v^{-1})^{d-k}  \frac{(v^{d-k+1} +v^{k-1-d}) (v^{d-k+2} +v^{k-d-2}) \cdots (v^{d-1} +v^{1-d})}{[k]!}
\\
&\quad \times \Big(v^k-v^{-k} +v^{d-k}(v^d+v^{-d}) \Big)
\\
&= \text{RHS}\eqref{eq:B4}.
\end{align*}
The proposition is proved.
\end{proof}

Recall $D_0, C_0, D_1, C_1$ from \eqref{eq:coeff0}--\eqref{eq:coeff2}.
\begin{theorem}
 \label{thm:DDCC}
For $d\ge 1$, we have
\[
C_0=C_1=D_0=D_1
= \frac{v^d (v +v^{-1}) (v^2 +v^{-2}) \cdots (v^{d-1} +v^{1-d})}{[d]!}.
\]
\end{theorem}

\begin{proof}
%The identity \eqref{DDCC}
Follows by Propositions~\ref{prop:DD}--\ref{prop:DD=CC} and the identity \eqref{eq:CC}.
\end{proof}
The identities \eqref{eq:DC0} and then \eqref{eq:coeff} follow from Theorem~\ref{thm:DDCC}.

%%%%%%%
%%%%%%%%%%
\section{Formula for a reflection functor $\Gamma_i$ ($i \not=\btau i$)}
 \label{sec:formulaRF2}

In this section, we establish a closed formula for the action of the reflection functor $\Gamma_i$ with $i \not =\btau i$.

\subsection{Formulas of $\Gamma_i$ for $i\neq \btau i$}

For any $m\in\N$ and $i \in \I$ such that $i \neq \tau i$, we define
\begin{align}
 \label{eq:SDP}
[S_i]^{(m)}:= \frac{[S_i]^m}{[m]_\sqq^!}.
\end{align}
We also set $[S_i]^{(m)}=0$, for $m<0$. The main result of this section is the following.

\begin{theorem}
\label{thm:braidqsplit}
Let $(Q,\btau)$ be an $\imath$quiver.
For any sink $i\in Q_0$ such that $c_{i,\btau i}=0$ and $j \neq i,\btau i$, we have
\begin{align}
\label{eq:braidqsplit}
\Gamma_i([S_j])
= &\sum^{-\max(c_{ij},c_{\tau i,j})}_{u=0} \sum^{-c_{ij}-u}_{r=0} \sum_{s=0}^{-c_{\tau i,j}-u}
(-1)^{r+s} \sqq^{-r-s+(r-s)u} (\sqq-\sqq^{-1})^{c_{ij}+c_{\tau i,j}+2u}
\\
 &\qquad \times [S'_{i}]^{(-c_{ij}-r-u)}*
  [S'_{\tau i}]^{(-c_{\tau i,j}-s-u)} *[S'_j]* [S'_{\tau i}]^{(s)} *[S'_{ i}]^{(r)} *[\E'_{\tau i}]^u.
  \notag
\end{align}
\end{theorem}
It follows by the assumption $c_{i,\btau i}=0$ that $i\neq \btau i$.

\begin{remark}
In case $c_{i,\btau i} \in 2\Z_{\le -1}$, the rank one $\imath$-subquiver corresponding to $\{i, \tau i\}$ is not of finite type, and the reflection functor $\Gamma_i$ is not defined. This is consistent with the structure of the relative Weyl group $W^{\btau}$ in \eqref{eq:Wtau}--\eqref{def:si} and Lemma~\ref{lem:iWeyl} below.
\end{remark}

The proof of Theorem \ref{thm:braidqsplit} will occupy the remainder of this section. The proof of \eqref{eq:braidqsplit} for $j \neq\btau j$ will be given in \S\ref{subsec:buildblock}--\S\ref{subsec:Tab=0}. We then explain in \S\ref{subsec:braidqsplitjfixed} how Theorem \ref{thm:braidqsplit} for $j = \btau j$ reduces to the case for $j \neq\btau j$.

\subsection{Summands of RHS\eqref{eq:braidqsplit}}
%{A rank 2 $\imath$quiver algebra}
\label{subsec:buildblock}

Assume $\btau j\neq j$. In this case, the proof of Theorem~\ref{thm:braidqsplit} reduces to the consideration of the rank 2 $\imath$quiver $(Q,\btau)$ as shown in the left figure of \eqref{eq:figKM2}. The quiver $\ov{Q'}$ of $\bs_i\Lambda^\imath$ (where $Q'={\bf s}_i Q$) is shown in the right figure of \eqref{eq:figKM2}.
Here $a=-c_{ij}=-c_{\btau i,\btau j}$, $b=-c_{j,\btau i}=-c_{i,\btau j}$, and $-c_{\btau j,\btau j}=2r$. %Then the quiver $\ov{Q}$ of $\Lambda^\imath$ is

\begin{equation}
\label{eq:figKM2}
\xymatrix{&i   \ar@[purple][rr]<0.5ex>^{\textcolor{purple}{\varepsilon_i}}  && {\btau i} \ar@[purple][ll]<0.5ex>^{\textcolor{purple}{\varepsilon_{\btau i}}}   \\
\ov{Q}=&\\
&j\ar[uu]|-{a} \ar[uurr]|-{\textcolor{white}{b}}\ar@/^0.5pc/@[purple][rr]<0.5ex>^{\textcolor{purple}{\varepsilon_j}}  \ar@<0.5ex>[rr]|-{r}& &{\btau j} \ar[uull]|-{b}  \ar[uu]|-{a} \ar@<0.5ex>[ll]|-{r} \ar@/^0.5pc/@[purple][ll]<0.5ex>^{\textcolor{purple}{\varepsilon_{\btau j}}}}  \qquad\qquad\qquad \xymatrix{&i\ar[dd]|-{a} \ar@[purple][rr]<0.5ex>^{\textcolor{purple}{\varepsilon_i}} \ar[ddrr]|-{\textcolor{white}{b}}&& {\btau i} \ar@[purple][ll]<0.5ex>^{\textcolor{purple}{\varepsilon_{\btau i}}}  \ar[dd]|-{a} \ar[ddll]|-{b} \\
\ov{Q'}=&\\
&j \ar@/^0.5pc/@[purple][rr]<0.5ex>^{\textcolor{purple}{\varepsilon_j}}  \ar@<0.5ex>[rr]|-{r}& &{\btau j}\ar@<0.5ex>[ll]|-{r} \ar@/^0.5pc/@[purple][ll]<0.5ex>^{\textcolor{purple}{\varepsilon_{\btau j}}}}
\end{equation}

The computations in this subsection will be performed in $\tMHi$.
We denote
\begin{align}
\mathcal{I}'_{k,l} =\{[M]&\in\Iso(\mod(kQ'))\mid \exists N\subseteq M \text{ s.t. }N\cong S_j, M/N\cong S_i^{\oplus k} \oplus S_{\btau i}^{ \oplus l}\},
\label{eq:Ik}
\\
p'(t,d,m,n) &=(m-d)^2+{ m-d \choose 2} +(n-d)^2 + {n-d \choose 2}
 \label{eq:p'}
\\
&\qquad\qquad+d(m+n-d)+ {d \choose 2}+ (t-(n-d))(n-d), \notag
\\
\label{eq:w}
w(m_1,m_2,n_1,&n_2,d,e,t_1,t_3)=c_{ij}(m_1-2d)+c_{\btau i,j}(n_1-2e) -m_1m_2-n_1n_2
\\
&\qquad\qquad+ 2(n_1-e)(n_2-d) +(e-d)(m_1+m_2-n_1-n_2)
\notag
\\
\notag
&\qquad\qquad+2(m_1-d)(m_2-e)+p'(t_3,d,m_1,n_2)+p'(t_1,e,n_1,m_2)+1.
\end{align}

For any $[M]\in\Iso(\mod(kQ'))$ with dimension vector $m\widehat{S'_i}+n\widehat{S'_{\btau i}}+\widehat{S'_j}$, there exists a unique indecomposable $kQ'$-module $N$ such that $M\cong N\oplus (S'_i)^{\oplus t^M_1}\oplus (S'_{\btau i})^{\oplus t^M_3}$ for some unique $t^M_1,t^M_3 \in \N$.

\begin{proposition}
\label{prop:build-block}
For any $m_1,m_2,n_1,n_2\in\N$, we have
\begin{align}
\label{eq:built1}
&[(S'_i)^{\oplus m_1}\oplus (S'_{\btau i})^{\oplus n_1}]*[S'_j]* [(S'_i)^{\oplus m_2} \oplus (S'_{\btau i})^{\oplus n_2}]\\\notag
=&\sum_{e=0}^{\min( n_1,m_2)}\sum_{d=0}^{\min(n_2,m_1 )}\sum_{[M]\in\mathcal{I}'_{m_1+m_2-d-e,n_1+n_2-d-e}} \sqq^{ w(m_1,m_2,n_1,n_2,d,e,t_1^M,t_3^M)}  \\\notag
& \qquad \times (\sqq-\sqq^{-1})^{m_1+m_2+n_1+n_2-d-e+1}
 \frac{[m_1]_\sqq^![n_2]_\sqq^!}{[d]_\sqq^!} \frac{[m_2]_\sqq^![n_1]_\sqq^!}{[e]_\sqq^!}
 \\
  &\qquad \times \qbinom{t^M_3}{n_2-d}_\sqq \qbinom{t^M_1}{m_2-e}_\sqq\frac{[M]}{|\aut(M)|} * [\E_i]^{ d}* [\E_{\btau i}]^{e}.\notag
\end{align}
\end{proposition}

The proof of Proposition~\ref{prop:build-block} is long and can be found in Appendix~\ref{app:prop}.

%
%\subsection{Summands of RHS \eqref{eq:braidqsplit}}

\begin{lemma}
We have
\begin{align}
\label{eq:RHSsummbraidsplit}
&[(S'_i)^{\oplus (-c_{ij}-u-r)}]*[(S'_{\btau i})^{\oplus (-c_{\tau i,j}-u-s)}]*[S'_j]* [(S'_{\btau i})^{\oplus s}]*[(S'_i)^{\oplus r}] \\\notag
&=\sum_{y=0}^{\min(-c_{ij}-u-r,-c_{\btau i,j}-u-s)}\sum_{x=0}^{\min(r,s)} \sum_{e=0}^{\min( -c_{\btau i,j}-u-s-y,r-x)}
\\\notag
&\sum_{d=0}^{\min(s-x,-c_{ij}-u-r-y )}\sum_{[M]\in\mathcal{I}'_{-c_{ij}-u-y-x-d-e,-c_{\btau i,j}-u-y-x-d-e}}\\\notag
&(\sqq-\sqq^{-1})^{-c_{ij}-c_{\btau i,j}-2u-x-y-d-e+1}\sqq^{(r-s)(x-y)+x(r+s-x)+ y(-c_{ij}-c_{\btau i,j}-2u-r-s-y) +\binom{x}{2}+\binom{y}{2}}
\\\notag
&\times\sqq^{w(-c_{ij}-u-r-y,r-x, -c_{\btau i,j}-u-s-y,s-x,d,e,t_1^M,t_3^M)}   \qbinom{t_3^M}{s-x-d}_\sqq \qbinom{t_1^M}{r-x-e}_\sqq
\\\notag
& \times \frac{[r]_\sqq^![s]_\sqq^![-c_{\tau i,j}-u-s]_\sqq^![-c_{ij}-u-r]_\sqq^!}{[d]_\sqq^![e]_\sqq^![x]_\sqq^![y]_\sqq^!}\frac{[M]}{|\aut(M)|} *[\E_i]^{ d+y}* [\E_{\btau i}]^{e+x}.
\end{align}
\end{lemma}

\begin{proof}
First, the following formula holds by a direct computation (cf. \eqref{eq:rank d}):
\begin{align*}
[(S'_{\btau i})^{\oplus s}]*[(S'_i)^{\oplus r}]%=
%&\sum_{x=0}^{\min(s,r)} \prod_{i=0}^{x-1} \frac{(q^s-q^i)(q^r-q^i)}{q^x-q^i} [(S'_{\btau i})^{\oplus (s-x)} \oplus (S'_i)^{\oplus(r-x)}\oplus \E_{\btau i}^{\oplus x}]\\
%=&\sum_{x=0}^{\min(s,r)} \sqq^{rx-sx} \prod_{i=0}^{x-1} \frac{(q^s-q^i)(q^r-q^i)}{q^x-q^i}
%[(S'_i)^{\oplus (r-x)} \oplus (S'_{\btau i})^{\oplus(s-x)}]* [\E_{\btau i}^{\oplus x}]
%\\
&=\sum_{x=0}^{\min(s,r)} \sqq^{x(r-s)}\sqq^{x(r+s-x)+\binom{x}{2}}(\sqq-\sqq^{-1})^{x}\qbinom{r}{x}_\sqq \qbinom{s}{x}_\sqq[x]_\sqq^! \\
&\qquad\qquad\qquad\times[(S'_i)^{\oplus (r-x)} \oplus (S'_{\btau i})^{\oplus(s-x)}]* [\E_{\btau i}^{\oplus x}].
\end{align*}
Hence, we also have
\begin{align*}
 &[(S'_i)^{\oplus (-c_{ij}-u-r)}]*[(S'_{\btau i})^{\oplus (-c_{\btau i,j}-u-s)}]
 %=&\sum_{y=0}^{\min(-c_{ij}-u-r,-c_{\btau i,j}-u-s)} \prod_{i=0}^{y-1} \frac{(q^{-c_{ij}-u-r}-q^i)(q^{-c_{\btau i,j}-u-s}-q^i)}{q^y-q^i} [(S'_i)^{\oplus (-c_{ij}-u-r-y)} \oplus (S'_{\btau i})^{\oplus(-c_{\btau i,j}-u-s-y)}\oplus \E_i^{\oplus y}] \notag
 %\\
% = \sum_{y=0}^{\min(-c_{ij}-u-r,-c_{\btau i,j}-u-s)} \sqq^{(c_{ij}-c_{\btau i,j}+ r-s)y}\times
  %\\
%  &\prod_{i=0}^{y-1} \frac{(q^{-c_{ij}-u-r}-q^i)(q^{-c_{\btau i,j}-u-s}-q^i)}{q^y-q^i} [(S'_i)^{\oplus (-c_{ij}-u-r-y)} \oplus (S'_{\btau i})^{\oplus(-c_{\btau i,j}-u-s-y)}]*[ \E_i^{\oplus y}].
\\
=&\sum_{y=0}^{\min(-c_{ij}-u-r,-c_{\btau i,j}-u-s)}\sqq^{y(c_{ij}-c_{\btau i,j} +r-s)}\sqq^{y(-c_{ij}-c_{\btau i,j}-2u-r-s-y)+\binom{y}{2}} (\sqq-\sqq^{-1})^{y}\qbinom{-c_{ij}-u-r}{y}_\sqq
 \\
&\qquad\qquad\qquad \times\qbinom{-c_{\btau i,j}-u-s}{y}_\sqq[y]_\sqq^! [(S'_i)^{\oplus (-c_{ij}-u-r-y)} \oplus (S'_{\btau i})^{\oplus(-c_{\btau i,j}-u-s-y)}]*[ \E_i^{\oplus y}].
\end{align*}
Therefore we obtain that
%\begin{align*}
%&[(S'_i)^{\oplus (-c_{ij}-u-r-y)} \oplus (S'_{\btau i})^{\oplus(-c_{\btau i,j}-u-s-y)}]*[S'_j]* [(S'_i)^{\oplus (r-x)} \oplus (S'_{\btau i})^{\oplus(s-x)}]*[ \E_i^{\oplus y}]* [\E_{\btau i}^{\oplus x}]
 %\\
%=& \sum_{e=0}^{\min( -c_{\btau i,j}-u-s-y,r-x)}\sum_{d=0}^{\min(s-x,-c_{ij}-u-r-y )}\sum_{[M]\in\mathcal{I}'_{-c_{ij}-u-y-x-d-e,-c_{\btau i,j}-u-y-x-d-e}}
% \\
% &\sqq^{ c_{ij}(-c_{ij}-u-r-y-2d)+c_{\btau i,j}(-c_{\btau i,j}-u-s-y- 2e) +(c_{ij}+u+r+y)(r-x)+(c_{\btau i,j}+u+s+y)(s-x) +(e-d)(c_{\btau i,j}-c_{ij}) } \\\notag
%&\sqq^{ 2(-c_{ij}-u-r-y-d)(r-x-e)+ 2(-c_{\btau i,j}-u-s-y-e)(s-x-d)} (\sqq-\sqq^{-1})^{-c_{ij}-c_{\btau i,j}-2u-2x-2y-d-e}
%\\
%&\sqq^{p(t_3^M,d,-c_{ij}-u-r-y,s-x)} \frac{[-c_{ij}-u-r-y]_\sqq^![s-x]_\sqq^!}{[d]_\sqq^!}\qbinom{t_3^M}{s-x-d}_\sqq \\
%&\sqq^{p(t_1^M,e,-c_{\btau i,j}-u-s-y,r-x)} \frac{[r-x]_\sqq^![-c_{\btau i,j}-u-s-y]_\sqq^!}{[e]_\sqq^!} \qbinom{t_1^M}{r-x-e}_\sqq\frac{(q-1)}{|\aut(M)|} [M]*[\E_i]^{ d}* %[\E_{\btau i}]^{e}.
%\end{align*}
\begin{align*}
&[(S'_i)^{\oplus (-c_{ij}-u-r)}]*[(S'_{\btau i})^{\oplus (-c_{\btau i,j}-u-s)}]*[S'_j]* [(S'_{\btau i})^{\oplus s}]*[(S'_i)^{\oplus r}] \\
=&  \sum_{y=0}^{\min(-c_{ij}-u-r,-c_{\btau i,j}-u-s)} \sum_{x=0}^{\min(r,s)} \sqq^{x(r-s)+y(c_{ij}-c_{\btau i,j}+r-s)}\sqq^{x(r+s-x)+ y(-c_{ij}-c_{\btau i,j}-2u-r-s-y) +\binom{x}{2}+\binom{y}{2}}
\\
& \times(\sqq-\sqq^{-1})^{x+y}\qbinom{r}{x}_\sqq \qbinom{s}{x}_\sqq[x]_\sqq^! \qbinom{-c_{ij}-u-r}{y}_\sqq\qbinom{-c_{\btau i,j}-u-s}{y}_\sqq[y]_\sqq^!
\\
 &[(S'_i)^{\oplus (-c_{ij}-u-r-y)} \oplus (S'_{\btau i})^{\oplus(-c_{\btau i,j}-u-s-y)}]*[ \E_i^{\oplus y}]*[S'_j]* [(S'_i)^{\oplus (r-x)} \oplus (S'_{\btau i})^{\oplus(s-x)}]* [\E_{\btau i}^{\oplus x}]
\\
=&  \sum_{y=0}^{\min(-c_{ij}-u-r,-c_{\btau i,j}-u-s)} \sum_{x=0}^{\min(r,s)} \sqq^{(r-s)(x-y)}\sqq^{x(r+s-x)+ y(-c_{ij}-c_{\btau i,j}-2u-r-s-y) +\binom{x}{2}+\binom{y}{2}}
\\
& \times(\sqq-\sqq^{-1})^{x+y}\qbinom{r}{x}_\sqq \qbinom{s}{x}_\sqq[x]_\sqq^! \qbinom{-c_{ij}-u-r}{y}_\sqq\qbinom{-c_{\btau i,j}-u-s}{y}_\sqq[y]_\sqq^!
\\
& [(S'_i)^{\oplus (-c_{ij}-u-r-y)} \oplus (S'_{\btau i})^{\oplus(-c_{\btau i,j}-u-s-y)}]*[S'_j]* [(S'_i)^{\oplus (r-x)} \oplus (S'_{\btau i})^{\oplus(s-x)}]*[ \E_i^{\oplus y}]* [\E_{\btau i}^{\oplus x}].
%\\
%=&\sum_{y=0}^{\min(-c_{ij}-u-r,-c_{\btau i,j}-u-s)}\sum_{x=0}^{\min(r,s)} \sum_{e=0}^{\min( -c_{\btau i,j}-u-s-y,r-x)}\sum_{d=0}^{\min(s-x,-c_{ij}-u-r-y )}\sum_{[M]\in\mathcal{I}'_{-c_{ij}-u-y-x-d-e,-c_{\btau i,j}-u-y-x-d-e}}\\
%&\sqq^{(r-s)(x-y)+x(r+s-x)+ y(-c_{ij}-c_{\btau i,j}-2u-r-s-y) +\binom{x}{2}+\binom{y}{2}} (\sqq-\sqq^{-1})^{x+y}
%\\
%& \qbinom{r}{x}_\sqq \qbinom{s}{x}_\sqq[x]_\sqq^! \qbinom{-c_{ij}-u-r}{y}_\sqq\qbinom{-c_{\btau i,j}-u-s}{y}_\sqq[y]_\sqq^!\qbinom{t_3^M}{s-x-d}_\sqq \qbinom{t_1^M}{r-x-e}_\sqq
%\\
%&\sqq^{w(-c_{ij}-u-r-y,r-x, -c_{\btau i,j}-u-s-y,s-x,d,e,t_1^M,t_3^M)} (\sqq-\sqq^{-1})^{-c_{ij}-c_{\btau i,j}-2u-2x-2y-d-e+1}
%\\
%& \frac{[r-x]_\sqq^![s-x]_\sqq^![-c_{\btau i,j}-u-s-y]_\sqq^![-c_{ij}-u-r-y]_\sqq^!}{[d]_\sqq^![e]_\sqq^!} \frac{[M]}{|\aut(M)|} *[\E_i]^{ d+y}* [\E_{\btau i}]^{e+x}
%\\
%=&\sum_{y=0}^{\min(-c_{ij}-u-r,-c_{\btau i,j}-u-s)}\sum_{x=0}^{\min(r,s)} \sum_{e=0}^{\min( -c_{\btau i,j}-u-s-y,r-x)}\sum_{d=0}^{\min(s-x,-c_{ij}-u-r-y )}\sum_{[M]\in\mathcal{I}'_{-c_{ij}-u-y-x-d-e,-c_{\btau i,j}-u-y-x-d-e}}\\
%&\sqq^{(r-s)(x-y)+x(r+s-x)+ y(-c_{ij}-c_{\btau i,j}-2u-r-s-y) +\binom{x}{2}+\binom{y}{2}}
%\\
%&\sqq^{w(-c_{ij}-u-r-y,r-x, -c_{\btau i,j}-u-s-y,s-x,d,e,t_1^M,t_3^M)} (\sqq-\sqq^{-1})^{-c_{ij}-c_{\btau i,j}-2u-x-y-d-e+1}
%\\
%&  \qbinom{t_3^M}{s-x-d}_\sqq \qbinom{t_1^M}{r-x-e}_\sqq  \frac{[r]_\sqq^![s]_\sqq^![-c_{\tau i,j}-u-s]_\sqq^![-c_{ij}-u-r]_\sqq^!}{[d]_\sqq^![e]_\sqq^![x]_\sqq^![y]_\sqq^!}\frac{[M]}{|\aut(M)|} *[\E_i]^{ d+y}* [\E_{\btau i}]^{e+x}.
\end{align*}
The lemma now follows from the above computation and applying \eqref{eq:built1}.
\end{proof}

\subsection{Reduction of the formula \eqref{eq:braidqsplit}}

Recall $\mathcal{I}'_{k,l}$ from \eqref{eq:Ik}, $p'(\cdot, \cdot, \cdot, \cdot)$ from \eqref{eq:p'} and the function $w$ from \eqref{eq:w}.
Then
\begin{align*}
[S'_i]^{(l)}:=\frac{[S'_i]^{*l}}{[l]^!_\sqq}=\sqq^{-\frac{l(l-1)}{2}}\frac{[(S'_i)^{\oplus l}]}{[l]^!_\sqq}.
\end{align*}
Using \eqref{eq:RHSsummbraidsplit}, we compute
\begin{align}
%&\sum^{-\max(c_{ij},c_{\tau i,j})}_{u=0} \sum^{-c_{ij}-u}_{r=0} \sum_{s=0}^{-c_{\tau i,j}-u}
%(-1)^{r+s} \sqq^{-r-s+(r-s)u} (\sqq-\sqq^{-1})^{c_{ij}+c_{\tau i,j}+2u}
 \text{RHS} \eqref{eq:braidqsplit} \label{eq:S5}
 %&[S'_{i}]^{(-c_{ij}-r-u)} [S'_{\tau i}]^{(-c_{\tau i,j}-s-u)} *[S'_j]* [S'_{\tau i}]^{(s)} *[S'_{ i}]^{(r)} *[\E'_{\tau i}]^u. \notag
 %\notag \\
 = &\sum^{-\max(c_{ij},c_{\tau i,j})}_{u=0} \sum^{-c_{ij}-u}_{r=0} \sum_{s=0}^{-c_{\tau i,j}-u}\sum_{y=0}^{\min(-c_{ij}-u-r,-c_{\tau i,j}-u-s)}\sum_{x=0}^{\min(r,s)}
  \\
  &\sum_{e=0}^{\min( -c_{\tau i,j}-u-s-y,r-x)}\sum_{d=0}^{\min(s-x,-c_{ij}-u-r-y )}\sum_{[M]\in\mathcal{I}'_{-c_{ij}-u-y-x-d-e,-c_{\tau i,j}-u-y-x-d-e}}
  \notag \\\notag
& (-1)^{r+s} \sqq^{(r-s)u-r-s} \sqq^{- {r \choose 2} -{s \choose 2}- { -c_{ij}-u-r \choose }-{-c_{\btau i,j}-u-s \choose 2}}
 \\
 &\notag\times\sqq^{(r-s)(x-y)+x(r+s-x)+ y(-c_{ij}-c_{\btau i,j}-2u-r-s-y) +\binom{x}{2}+\binom{y}{2}}
 (\sqq-\sqq^{-1})^{-x-y-d-e+1}
\\\notag
&\times\sqq^{w(-c_{ij}-u-r-y,r-x, -c_{\btau i,j}-u-s-y,s-x,d,e,t_1^M,t_3^M)}   \qbinom{t_3^M}{s-x-d}_\sqq \qbinom{t_1^M}{r-x-e}_\sqq
\\\notag
& \times \frac{1}{[d]_\sqq^![e]_\sqq^![x]_\sqq^![y]_\sqq^!}
\frac{[M]}{|\aut(M)|}*[\E_1]^{ d+y}* [\E_3]^{e+x+u}.
\notag
 \end{align}

Introduce
\begin{align*}
C& (f,g,t_1,t_3)
\\\notag
& := \sum^{\min(a,b)}_{u=0} \sum^{a-u}_{r=0} \sum_{s=0}^{b-u}\sum_{y=0}^{\min(a-u-r,b-u-s)}\sum_{x=0}^{\min(r,s)}
\\
&\delta(0\leq g-y\leq \min(a-u-r-y,s-x))\delta(0\leq f-x-u\leq \min(b-u-s-y,r-x))\notag
\\\notag
& (-1)^{r+s}  \sqq^{- {r \choose 2} -{s \choose 2}- { a-u-r \choose }-{b-u-s \choose 2}}
 \sqq^{(r-s)(x-y+u)-r-s+x(r+s-x)+ y(a+b-2u-r-s-y) +\binom{x}{2}+\binom{y}{2}}
 \\
 \notag
&\times\sqq^{w(a-u-r-y,r-x, b-u-s-y,s-x,g-y,f-x-u,t_1,t_3)} (\sqq-\sqq^{-1})^{-f-g+u+1}
%& (-1)^{r+s} \sqq^{-r-s+(r-s)u} \sqq^{-{r \choose 2}-{s \choose 2}-{a-u-r \choose 2}-{b-u-s \choose 2}}
%\sqq^{(r-s)(x-y)} \sqq^{(r+s-x)x +{x \choose 2} +(a+b-2u-r-s-y)y+{y \choose 2}}  \\\notag
%&\sqq^{- a(a-u-r+y-2g)-b(b+u-s-y- 2f+2x) +(-a+u+r+y)(r-x)+(-b+u+s+y)(s-x) +(f-g-x+y-u)(a-b) } \\\notag
%&\sqq^{ 2(a-u-r-g)(r-f+u)+ 2(b-s-y-f+x)(s-x-g+y)} \sqq^{p'(t_1,f-x-u,b-u-s-y,r-x)+p'(t_3,g-y,a-u-r-y,s-x)}
\\\notag
&\times \qbinom{t_3}{s-x-g+y}_\sqq \qbinom{t_1}{r-f+u}_\sqq \cdot
 \frac{1}{[g]_\sqq^! [f-u]_\sqq^!} \qbinom{g}{y}_\sqq \qbinom{f-u}{x}_\sqq,
%\frac{1}{[g-y]_\sqq^![f-x-u]_\sqq^![x]_\sqq^![y]_\sqq^!}.
\notag
\end{align*}
where $f,g,t_1,t_3$ are subject to the constraints
\begin{align}
\label{cond}
0\leq f,g\leq \min(a,b), \; 0\leq t_1\leq a-f-g, \; 0\leq t_3\leq b-f-g.
\end{align}
Recall $a=-c_{ij}$, $b=-c_{\tau i,j}$. For fixed $M$ (i.e., fixed $t_1^M, t_3^M$) and fixed $f:=e+x+u$, $g:=d+y$, the coefficient of $\frac{[M]}{|\Aut(M)|}*[\E_1]^{g}* [\E_3]^{f}$ in the sum \eqref{eq:S5} is $C(f,g,t_1^M,t_3^M)$. Note $e=f-x-u$, and $d=g-y$.

\vspace{2mm}
{\bf Claim.} The $\delta$ functions in $C(f,g,t_1,t_3)$ are all equal to 1.

Let us prove the Claim. Since the product of $v$-binomials above is $0$ whenever $s-x-g+y<0$ or $r-u-f<0$ or $g-y<0$ or $f-x-u<0$, $\delta(g-y\leq \min(a-u-r-y,s-x))$ can be replaced by $\delta(g-y\leq  a-u-r-y)$ while $\delta(f-x-u\leq \min(b-u-s-y,r-x))$ can be replaced by $\delta(f-x-u\leq b-u-s-y)$.

But if $g-y >  a-u-r-y$, or equivalently, $a-f-g< r-f+u$, then it follows by \eqref{cond} that $t_1 \leq a-f-g <r-f+u$ and the $v$-binomial product in $C(f,g,t_1,t_3)$ is $0$. So $\delta(g-y\leq \min(a-u-r-y,s-x))$ is removable from the above summations unconditionally.

Similarly, if $f-x-u > b-u-s-y$, or equivalently, $b-f-g <s-x-g+y$, then it follows by \eqref{cond} that $t_3\leq b-f-g < s-x-g+y$ and the $v$-binomial product in $C(f,g,t_1,t_3)$ is $0$. So  $\delta(f-x-u\leq b-u-s-y)$ is removable as well.
The Claim is proved.
\vspace{2mm}

Denote
\begin{align}
T (f,g,t_1,t_3) \label{eq:Tab}
&= \sum^{\min(a,b)}_{u=0} \sum^{a-u}_{r=0} \sum_{s=0}^{b-u}\sum_{y=0}^{\min(a-u-r,b-u-s)}\sum_{x=0}^{\min(r,s)} (-1)^{r+s} (\sqq-\sqq^{-1})^{u}
 \\
&\quad \times
  \sqq^{2(fr+gs)-\frac{1}{2}u(u+1)+u(f-x+t_1)-x(2g +f+t_3)+y(g+t_3)+t_3s-s+t_1r-r}
\notag \\
&\quad \times \qbinom{t_3}{s-x-g+y}_\sqq \qbinom{t_1}{r-f+u}_\sqq
  \frac{1}{[g]_\sqq^! [f-u]_\sqq^!} \qbinom{g}{y}_\sqq \qbinom{f-u}{x}_\sqq.
\notag
\end{align}
Recall the function $w$ defined in \eqref{eq:w}.
By a direct computation of the $\sqq$-powers in $C(f,g,t_1,t_3)$ above, we can rewrite
\begin{align}
\label{eq:C=T}
&C(f,g,t_1,t_3)
= \sqq^{-(a+b)(f+g)+4fg+\frac{1}{2}f(3f+1)+\frac{1}{2}g(3g+1)-t_1f-t_3g +1} (\sqq-\sqq^{-1})^{-f-g+1}   T(f,g,t_1,t_3).
\end{align}

Fix $[M]   \in\mathcal{I}'_{a-f-g,b-f-g}$. For $f=g=t_1=t_3=0$, $M$ is indecomposable which is isomorphic to $F_i^+(S'_j)$ by Lemma~ \ref{lem:reflecting dimen}. In this case, the coefficient of $[M]$  of RHS\eqref{eq:braidqsplit} is equal to $(q-1)C(0,0,0,0)=T(0,0,0,0)$ by \eqref{eq:C=T} and noting that $|\Aut(M)|=q-1$. Summarizing, we have reached the following reduction toward the proof of \eqref{eq:braidqsplit}.

\begin{proposition}
  \label{prop:GammaT}
The formula  \eqref{eq:braidqsplit} is equivalent to the identities
\begin{align*}
  %\label{eq:Tab=0}
T(0,0,0,0) =1, \quad \text { and } \; T(f,g,t_1,t_3)=0,
\end{align*}
for non-negative integers $f,g,t_1,t_3$ subject to the following constraints:
\begin{align}
 \label{eq:fgtt}
f,g\leq \min(a,b), \;\;
t_1\leq a-f-g, \;\;
t_3\leq b-f-g, \;\;
\text{ and not all $f,g,t_1,t_3$ are zero.}
\end{align}
%$0\leq f,g\leq \min(a,b)$, $0\leq t_1\leq a-f-g, 0\leq t_3\leq b-f-g$ and not all $f,g,t_1,t_3$ are zero.
\end{proposition}

\subsection{Proof of the identities \eqref{eq:Tab=0}}
\label{subsec:Tab=0}

This subsection is devoted to the proof of the following.
\begin{proposition}
  \label{prop:T}
The following identity holds:
\begin{align}
 T(f,g,t_1,t_3)& =0,
\label{eq:Tab=0}
\end{align}
for non-negative integers $f,g,t_1,t_3$ satisfying the conditions \eqref{eq:fgtt}. Moreover, we have $T(0,0,0,0) =1$.
\end{proposition}

Below we shall replace $\sqq$ by the free variable $v$. We shall change variables
\[
s_3 =s -x +y -g,
\qquad
r_1 =r +u -f.
\]
In other words, we have $s =s_3 +x -y +g$, and $r= r_1 -u +f$.
We can rewrite the sign in \eqref{eq:Tab} as
\begin{align}  \label{eq:sign}
(-1)^{r+s} =(-1)^{g+f} (-1)^{u +x}  (-1)^{r_1}  (-1)^{s_3} (-1)^{y}.
\end{align}
We rewrite the $v$-power in \eqref{eq:Tab} as
\begin{align*}
& {2(fr+gs)-\frac{1}{2}u(u+1)+u(f-x+t_1)-x(2g +f+t_3)+y(g+t_3)+t_3s-s+t_1r-r}
\\
&=
(2f^2+2g^2 +gt_3-g +ft_1 -f)  +  P,
\end{align*}
where
\begin{align}   \label{eq:power}
P =u(1-f) -\frac12 u(u+1)  -(f+u+1) x  +(2f +t_1 -1) r_1 + (2g +t_3 -1) s_3 +(1-g)y.
\end{align}

Looking closely, we see $y$ runs from 0 to $g$ freely in \eqref{eq:Tab}.
By pulling out the parts relevant to $y$ in \eqref{eq:Tab}, we obtain the following factor of \eqref{eq:Tab}:
\begin{align*}
\sum_{y=0}^g (-1)^y v^{(1-g)y} \qbinom{g}{y} =
\begin{cases}
1, & \text{ if } g=0,
\\
0, &\text{ if } g>0.
\end{cases}
\end{align*}
%which is equal to $1$ if $g=0$, and equal to $0$ if $g>0$;
Hence the identity \eqref{eq:Tab=0} holds when $g>0$.

From now on, we shall assume $g=0$. We observe the summation over $s_3$ is taken freely from $0$ to $t_3$ in \eqref{eq:Tab}. By pulling out the parts relevant to $s_3$ in \eqref{eq:Tab} and using \eqref{eq:sign}--\eqref{eq:power}, we obtain the following factor of \eqref{eq:Tab}:
\begin{align*}
\sum_{s_3=0}^{t_3} (-1)^{s_3} %v^{2g s_3}
v^{(t_3-1)s_3} \qbinom{t_3}{s_3}
=\prod_{a=0}^{t_3-1} (1-v^{2a %+2g
})
 =
\begin{cases}
1, & \text{ if } t_3=0,
\\
0, &\text{ if } t_3>0.
\end{cases}
\end{align*}
%which is equal to $1$ if $t_3=0$, and equal to $0$ if $t_3>0$;
Hence the identity \eqref{eq:Tab=0} holds when $t_3>0$.

Now we shall assume $t_3 =g=0$.
We observe the summation over $r_1$ is taken freely from $0$ to $t_1$ in \eqref{eq:Tab}. By pulling out the parts relevant to $r_1$ in \eqref{eq:Tab} and using \eqref{eq:sign}--\eqref{eq:power}, we obtain the following factor of \eqref{eq:Tab}
\begin{align}
 \label{eq:tg0}
\sum_{r_1=0}^{t_1} (-1)^{r_1}  v^{2f r_1} v^{(t_1-1)r_1} \qbinom{t_1}{r_1}
=\prod_{a=0}^{t_1 -1} (1-v^{2a +2f}),
\end{align}
which follows by \eqref{1.3.1} with $z= - v^{2f}$; the RHS of \eqref{eq:tg0} is $0$ if $f=0$ (and in this case we must have $t_1>0$ since not all $f,g,t_1,t_3$ are zero by  \eqref{eq:fgtt}). Hence the identity \eqref{eq:Tab=0} holds when $f=0$.

Now we shall assume $t_3 =g=0$ and $f>0$.
By pulling out the remaining double summations over $u, x$ (where we first sum over $x$ freely from $0$ to $f-u$) in \eqref{eq:Tab} and using \eqref{eq:sign}--\eqref{eq:power}, we obtain the following factor of \eqref{eq:Tab}:
\begin{align*}
&\sum_{u=0}^{f} \sum_{x=0}^{f-u} (-1)^{u+x}  v^{u(1-f) -\frac12 u(u+1)  -(f+u+1) x} \frac1{[x]![f-u-x]!} (v-v^{-1})^u
\\
&= \sum_{u=0}^{f} (-1)^u v^{u(1-f) -\frac12 u(u+1)}  (v-v^{-1})^u\sum_{x=0}^{f-u} (-1)^{x}  v^{-(f+u+1) x} \frac1{[x]![f-u-x]!}
\\
&= \sum_{u=0}^{f} (-1)^u v^{u(1-f) -\frac12 u(u+1)} \frac{(v-v^{-1})^u}{[f-u]!} \sum_{x=0}^{f-u} (-1)^{x}  v^{-2fx}v^{(f-u-1) x}  \qbinom{f-u}{x}
\\
&\stackrel{(*)}{=} \sum_{u=0}^{f} (-1)^u v^{u(1-f) -\frac12 u(u+1)} \frac{(v-v^{-1})^u}{[f-u]!} \prod_{a=0}^{f-u-1} (1 -v^{2a -2f})
\\
&\stackrel{(**)}{=}  \sum_{u=0}^{f} (-1)^u v^{u(1-f)} \qbinom{f}{u} \cdot v^{ -\frac12 f(f+1)}  (v-v^{-1})^f =0.
\end{align*}
The identity ($*$) above follows by \eqref{1.3.1} with $z=-v^{2f}$, and the equality ($**$)   follows from %$\qbinom{f}{u} = \frac{[f][f-1]\cdots [u+1]}{[f-u]!}$ and
\[
\prod_{a=0}^{f-u-1} (1 -v^{2a -2f})  =v^{-\frac12 f(f+1) +\frac12 u(u+1)} (v-v^{-1})^{f-u} \cdot [f][f-1]\cdots [u+1].
\]

This completes the proof of the identity \eqref{eq:Tab=0}.
In addition, we read off from the above proof that $T(0,0,0,0)=1$. Proposition~\ref{prop:T} is proved.

%%
%\subsection{Formulas for braid group action}

%\begin{lemma}
%We have the following automorphisms.
% \begin{align}
%\TT'_{i,e}(B_j)
%= &\sum^{-\max(c_{ij},c_{\tau i,j})}_{u=0} \; \sum^{-c_{ i,j}-u}_{r=0} \; \sum_{s=0}^{-c_{\tau i,j}-u}
%\\
%& (-1)^{r+s} v^{e\big(r-s+(-c_{ij}-r-s-u)u\big)} B_i^{(r)} B_{\tau i}^{(-c_{\tau i,j}-u-s)} B_j B_{\tau i}^{(s)} B_i^{(-c_{ij}-r-u)}\tk_{\tau i}^u,
%\notag
%\end{align}
% \begin{align}
%\TT''_{i,-e}(B_j)
%= &\sum^{-\max(c_{ij},c_{\tau i,j})}_{u=0} \; \sum^{-c_{ i,j}-u}_{r=0} \; \sum_{s=0}^{-c_{\tau i,j}-u}
%\\
%&(-1)^{r+s} v^{e\big(r-s+(-c_{ij}-r-s-u)u\big)} \tk_{ i}^u B_i^{(-c_{ij}-r-u)}  B_{\tau i}^{(s)} B_j  B_{\tau i}^{(-c_{\tau i,j}-u-s)} B_i^{(r)}.
%\notag
%\end{align}
%\end{lemma}

%
%
\subsection{Proof of Theorem~ \ref{thm:braidqsplit}}
\label{subsec:braidqsplitjfixed}

For the case $j \neq \tau j$, the formula \eqref{eq:braidqsplit} (or Theorem~ \ref{thm:braidqsplit}) follows by Proposition~\ref{prop:GammaT} and Proposition~\ref{prop:T}.

For the remaining case $j = \btau j$, the proof of Theorem~\ref{thm:braidqsplit} is reduced to the consideration of the rank 2 $\imath$quiver $(Q,\btau)$ as shown in the left figure of \eqref{diag:qsplit KM2}. The quiver $\ov{Q'}$ of $\bs_i\Lambda^\imath$ (associated to $Q' =\bs_i Q$)  is shown in the right figure of \eqref{diag:qsplit KM2}.
Here $a=-c_{ij}=-c_{\btau i,j}$. %Then the quiver $\ov{Q}$ of $\Lambda^\imath$ is
\begin{center}\setlength{\unitlength}{1mm}
\vspace{.1cm}
\begin{equation}
\label{diag:qsplit KM2}
\begin{picture}(50,10)(0,-10)

 \put(-20,0){
 \begin{picture}(50,10)
\put(0,-2){$i$}
\put(20,-2){$\btau i$}

  \put(-15,-10){$\ov{Q}=$}

\put(10,-18){\line(-1,2){3}}
\put(5.5,-12){$a$}
\put(6,-10){\vector(-1,2){3.5}}

\put(12.5,-18){\line(1,2){3}}
\put(15,-12){$a$}
\put(16.5,-10){\vector(1,2){3.5}}

\put(10,-22){$j$}
\color{purple}

\put(2.5,1){\vector(1,0){17}}
\put(19.5,-1){\vector(-1,0){17}}
\put(10,1){$^{\varepsilon_i}$}
\put(10,-3){$_{\varepsilon_{\btau i}}$}
\put(10,-28){$_{\varepsilon_j}$}
\begin{picture}(50,23)(-10,19)
\color{purple}
\qbezier(-1,-1)(-3,-3)(-2,-5.5)
\qbezier(-2,-5.5)(1,-9)(4,-5.5)
\qbezier(4,-5.5)(5,-3)(3,-1)
\put(3.1,-1.4){\vector(-1,1){0.3}}
\end{picture}
\end{picture}}

\put(50,0){
\begin{picture}(50,10)
\put(0,-2){$i$}
\put(20,-2){$\btau i$}

\put(-15,-10){$\ov{Q'}=$}
\put(7,-12){\vector(1,-2){3.5}}
\put(5.5,-12){$a$}
\put(6,-10){\line(-1,2){3.5}}

\put(15.5,-12){\vector(-1,-2){3.5}}
\put(15,-12){$a$}
\put(16.5,-10){\line(1,2){3.5}}
\put(10,-22){$j$}

\color{purple}
\put(2.5,1){\vector(1,0){17}}
\put(19.5,-1){\vector(-1,0){17}}
\put(10,1){$^{\varepsilon_i}$}
\put(10,-3){$_{\varepsilon_{\btau i}}$}
\put(10,-28){$_{\varepsilon_j}$}
\begin{picture}(50,23)(-10,19)
\color{purple}
\qbezier(-1,-1)(-3,-3)(-2,-5.5)
\qbezier(-2,-5.5)(1,-9)(4,-5.5)
\qbezier(4,-5.5)(5,-3)(3,-1)
\put(3.1,-1.4){\vector(-1,1){0.3}}
\end{picture}
\end{picture}}
\end{picture}
\vspace{1.8cm}
\end{equation}
\end{center}
By the same argument as in \cite[Proposition 9.4]{LW20a}, the computations involved in proving \eqref{eq:braidqsplit} with $j =\btau j$ are the same as for the $\imath$quiver in the diagram \eqref{eq:figKM2} for $j \neq \btau j$ and $b=a$.

%%%%%%%%%%%
%%%%%%%%%%%

\section{Symmetries of $\imath$quantum groups}
  \label{sec:braid}

\subsection{Quantum groups}
  \label{subsec:QG}

Let $Q$ be a quiver (without loops) with vertex set $Q_0= \I$.
Recall that $n_{ij}$ is the number of edges connecting vertices $i$ and $j$. Let $C=(c_{ij})_{i,j \in \I}$ be the symmetric generalized Cartan matrix of the underlying graph of $Q$, defined by $c_{ij}=2\delta_{ij}-n_{ij}.$ Let $\fg$ be the corresponding Kac-Moody Lie algebra, with the Chevalley involution denoted by $\omega$. Let $\alpha_i$ ($i\in\I $) be the simple roots of $\fg$, and denote the root lattice by $\Z^{\I}:=\Z\alpha_1\oplus\cdots\oplus\Z\alpha_n$. The {\em simple reflection} $s_i:\Z^{\I}\rightarrow\Z^{\I}$ is defined to be $s_i(\alpha_j)=\alpha_j-c_{ij}\alpha_i$, for $i,j\in \I$.
Denote the Weyl group by $W =\langle s_i\mid i\in \I\rangle$.

Let $\btau$ be an involution of $Q$, which induces an involution on $\fg$ again denoted by $\btau$. We shall define the {\em restricted Weyl group} associated to the quasi-split symmetric pair $(\fg, \fg^{\omega\tau})$ to be the following subgroup $W^{\btau}$ of $W$:
\begin{align}
  \label{eq:Wtau}
W^{\btau} =\{w\in W\mid \btau w =w \btau\}
\end{align}
where $\btau$ is regarded as an automorphism of $\Aut(C)$. In finite type, it is well known that the restricted Weyl group defined this way coincides with the one arising from real groups (cf., e.g., \cite{KP11}).

Recall the subset $\I_\tau$ of $\I$ from \eqref{eq:ci}, and define
\begin{align}
 \label{eq:Itau2}
\ov{\I}_\btau:=\{i\in\I_\btau\mid c_{i,\btau i}=0 \text{ or }2 \}.
\end{align}
In our setting, $\ov{\I}_\btau$ consists of exactly those $i \in \I_\btau$ such that the $\btau$-orbit of $i$ is of finite type. Note that $\ov{\I}_\btau=\I_\btau$ if $(Q,\btau)$ is acyclic. We denote by $\bs_{i}$, for $i\in\ov{\I}_\btau$, the following element of order 2 in the Weyl group $W$
\begin{align}
\label{def:si}
\bs_i= \left\{
\begin{array}{ll}
s_{i}, & \text{ if } i=\btau i
\\
s_is_{\btau i}, & \text{ if } i\neq \btau i.
\end{array}
\right.
\end{align}

\begin{lemma}
[\text{\cite[Appendix]{Lus03}}]
  \label{lem:iWeyl}
The restricted Weyl group $W^{\btau}$ can be identified with a Coxeter group with $\bs_i$ ($i\in \ov{\I}_\btau$) as its generators.
\end{lemma}

Let $\bv$ be an indeterminant. Write $[A, B]=AB-BA$. Then $\tU = \tU_\bv(\fg)$ is defined to be the $\Q(\bv)$-algebra generated by $E_i,F_i, \tK_i,\tK_i'$, $i\in \I$, where $\tK_i, \tK_i'$ are invertible, subject to the following relations:
\begin{align}
[E_i,F_j]= \delta_{ij} \frac{\tK_i-\tK_i'}{\bv-\bv^{-1}},  &\qquad [\tK_i,\tK_j]=[\tK_i,\tK_j']  =[\tK_i',\tK_j']=0,
\label{eq:KK}
\\
\tK_i E_j=\bv^{c_{ij}} E_j \tK_i, & \qquad \tK_i F_j=\bv^{-c_{ij}} F_j \tK_i,
\label{eq:EK}
\\
\tK_i' E_j=\bv^{-c_{ij}} E_j \tK_i', & \qquad \tK_i' F_j=\bv^{c_{ij}} F_j \tK_i',
 \label{eq:K2}
\end{align}
 and the quantum Serre relations, for $i\neq j \in \I$,
\begin{align}
& \sum_{r=0}^{1-c_{ij}} (-1)^r \left[ \begin{array}{c} 1-c_{ij} \\r \end{array} \right]  E_i^r E_j  E_i^{1-c_{ij}-r}=0,
  \label{eq:serre1} \\
& \sum_{r=0}^{1-c_{ij}} (-1)^r \left[ \begin{array}{c} 1-c_{ij} \\r \end{array} \right]  F_i^r F_j  F_i^{1-c_{ij}-r}=0.
  \label{eq:serre2}
\end{align}
Note that $\tK_i \tK_i'$ are central in $\tU$ for all $i$.
The comultiplication $\Delta: \widetilde{\U} \rightarrow \widetilde{\U} \otimes \widetilde{\U}$ is given by
\begin{align}  \label{eq:Delta}
\begin{split}
\Delta(E_i)  = E_i \otimes 1 + \tK_i \otimes E_i, & \quad \Delta(F_i) = 1 \otimes F_i + F_i \otimes \tK_{i}', \\
 \Delta(\tK_{i}) = \tK_{i} \otimes \tK_{i}, & \quad \Delta(\tK_{i}') = \tK_{i}' \otimes \tK_{i}'.
 \end{split}
\end{align}
%The Chevalley involution $\omega$ on $\tU$ is given by
%\begin{align}  \label{eq:omega}
%\omega(E_i)  = F_i,\quad  \omega(F_i) = E_i,\quad \omega(\tK_{i}) = \tK_{i}' , \quad \omega(\tK_{i}') =\tK_{i}, \quad \forall i\in \I.
%\end{align}

Analogously as for $\tU$, the quantum group $\bU$ is defined to be the $\Q(v)$-algebra generated by $E_i,F_i, K_i, K_i^{-1}$, $i\in \I$, subject to the  relations modified from \eqref{eq:KK}--\eqref{eq:serre2} with $\tK_i$ and $\tK_i'$ replaced by $K_i$ and $K_i^{-1}$, respectively. The comultiplication $\Delta$ is obtained by modifying \eqref{eq:Delta} with $\tK_i$ and $\tK_i'$ replaced by $K_i$ and $K_i^{-1}$, respectively (cf. \cite{Lus93}; beware that our $K_i$ has a different meaning from $K_i \in \U$ therein.)

\subsection{$\imath$Quantum groups}
  \label{subsec:iQG}

For a  (generalized) Cartan matrix $C=(c_{ij})$, let $\Aut(C)$ be the group of all permutations $\btau$ of the set $\I$ such that $c_{ij}=c_{\btau i,\btau j}$. An element $\btau\in\Aut(C)$ is called an \emph{involution} if $\btau^2=\Id$.

Let $\btau$ be an involution in $\Aut(C)$. We define $\widetilde{\bU}^\imath=\widetilde{\bU}'_\bv(\fk)$ to be the $\Q(v)$-subalgebra of $\tU$ generated by
\begin{equation}
  \label{eq:Bi}
B_i= F_i +  E_{\btau i} \tK_i',
\qquad \tk_i = \tK_i \tK_{\btau i}', \quad \forall i \in \I.
\end{equation}
Let $\tU^{\imath 0}$ be the $\Q(v)$-subalgebra of $\tUi$ generated by $\tk_i$, for $i\in \I$.
By \cite[Lemma 6.1]{LW19}, the elements $\tk_i$ (for $i= \btau i$) and $\tk_i \tk_{\btau i}$  (for $i\neq \btau i$) are central in $\tUi$.

Let $\bvs=(\vs_i)\in  (\Q(\bv)^\times)^{\I}$ be such that $\vs_i=\vs_{\btau i}$ for all $i$. % for each $i\in \I$ which satisfies $c_{i, \btau i}=0$.
Let $\Ui:=\Ui_{\bvs}$ be the $\Q(v)$-subalgebra of $\bU$ generated by
\[
B_i= F_i+\vs_i E_{\btau i}K_i^{-1},
\quad
k_j= K_jK_{\btau j}^{-1},
\qquad  \forall i \in \I, j \in \I\backslash\ci.
\]
It is known \cite{Let99, Ko14} that $\bU^\imath$ is a right coideal subalgebra of $\bU$ in the sense that $\Delta: \Ui \rightarrow \Ui\otimes \U$; and $(\bU,\Ui)$ is called a \emph{quantum symmetric pair} ({\em QSP} for short), as they specialize at $v=1$ to $(U(\fg), U(\fg^{\omega\tau}))$.
%where $\theta=\omega \circ \btau$, $\omega$ is the  Chevalley involution, and $\btau$ is understood here as an automorphism of $\fg$.

The algebras $\Ui_{\bvs}$, for $\bvs \in  (\Q(\bv)^\times)^{\I}$, are obtained from $\tUi$ by central reductions.

\begin{proposition} [\text{\cite[Proposition 6.2]{LW19}}]
(1) The algebra $\Ui$ is isomorphic to the quotient of $\tUi$ by the ideal generated by
\begin{align}   \label{eq:parameters}
\tk_i - \vs_i \; (\text{for } i =\btau i),
\qquad  \tk_i \tk_{\btau i} - \vs_i \vs_{\btau i}  \;(\text{for } i \neq \btau i).
\end{align}
The isomorphism is given by sending $B_i \mapsto B_i, k_j \mapsto \vs_{\btau j}^{-1} \tk_j, k_j^{-1} \mapsto \vs_{j}^{-1} \tk_{\btau j}, \forall i\in \I, j\in \I\backslash\ci$.

(2) The algebra $\widetilde{\bU}^\imath$ is a right coideal subalgebra of $\widetilde{\bU}$; that is, $(\widetilde{\bU}, \widetilde{\bU}^\imath)$ forms a QSP.
\end{proposition}

We shall refer to $\tUi$ and $\Ui$ as {\em (quasi-split) $\imath${}quantum groups}; they are called {\em split} if $\btau =\Id$.

For $i\in \I$ with $\btau i= i$, generalizing the constructions in \cite{BW18, BeW18}, we define the {\em $\imath${}divided powers} of $B_i$ to be (see also \cite{CLW21})
\begin{eqnarray}
&&\ff_{i,\odd}^{(m)}=\frac{1}{[m]!}\left\{ \begin{array}{ccccc} B_i\prod_{s=1}^k (B_i^2-v\tk_i[2s-1]^2 ) & \text{if }m=2k+1,\\
\prod_{s=1}^k (B_i^2-v\tk_i[2s-1]^2) &\text{if }m=2k; \end{array}\right.
  \label{eq:iDPodd} \\
&&\ff_{i,\ev}^{(m)}= \frac{1}{[m]!}\left\{ \begin{array}{ccccc} B_i\prod_{s=1}^k (B_i^2-v\tk_i[2s]^2 ) & \text{if }m=2k+1,\\
\prod_{s=1}^{k} (B_i^2-v\tk_i[2s-2]^2) &\text{if }m=2k. \end{array}\right.
 \label{eq:iDPev}
\end{eqnarray}
On the other hand, for $i\in \I$ with $i \neq \tau i$, we define the divided powers as in the quantum group setting: for $m\in \N$,
\begin{align}
   \label{eq:DP}
\ff_{i}^{(m)}= \frac{\ff_i^m}{[m]^!}.
\end{align}

We have the following {\em Serre presentation} of $\tUi$, with $\ov{p}_i\in \Z_2$ fixed for each $i\in \I$.

\begin{proposition} [\text{\cite[Theorem 4.2]{LW20a}; also cf. \cite{Let02}}]
 \label{prop:Serre}
The $\Q(v)$-algebra $\tUi$ has a presentation with generators $B_i$, $\tk_i$ $(i\in \I)$ and the relations \eqref{relation1}--\eqref{relation6} below: for $\ell \in \I$, and $i\neq j \in \I$,
\begin{align}
\tk_i \tk_\ell =\tk_\ell \tk_i,
\quad
\tk_i B_\ell & = v^{c_{\btau i,\ell} -c_{i \ell}} B_\ell \tk_i,
   \label{relation1}
\\
B_iB_{j}-B_jB_i &=0, \quad \text{ if }c_{ij} =0 \text{ and }\btau i\neq j,\label{relation2}
\\
\sum_{n=0}^{1-c_{ij}} (-1)^nB_i^{(n)}B_jB_i^{(1-c_{ij}-n)} &=0, \quad \text{ if } j \neq i\neq \btau i, \label{relation3}
\\
\sum_{n=0}^{1-c_{i,\btau i}} (-1)^{n+c_{i,\btau i}}B_i^{(n)}B_{\btau i}B_i^{(1-c_{i,\btau i}-n)}& =
\label{relation5}     \\
  \frac{1}{v-v^{-1}}  \Big(v^{c_{i,\btau i}} (v^{-2};v^{-2})_{-c_{i,\btau i}}    &
   B_i^{(-c_{i,\btau i})} \tk_i
   -(v^{2};v^{2})_{-c_{i,\btau i}}B_i^{(-c_{i,\tau i})} \tk_{\btau i}  \Big),
\text{ if } \btau i \neq i,
 \notag \\
\sum_{n=0}^{1-c_{ij}} (-1)^n  B_{i, \overline{p_i}}^{(n)}B_j B_{i,\overline{c_{ij}}+\overline{p}_i}^{(1-c_{ij}-n)} &=0,\quad   \text{ if }i=\btau i.
\label{relation6}
\end{align}
\end{proposition}

\subsection{$\imath$Quantum groups and $\imath$Hall algebras}
  \label{subsec:iQGiH}

Recall $\I_\btau$ from \eqref{eq:ci}.
Let $\bvs=(\vs_i)\in   (\Q(\sqq)^\times)^{\I}$ be such that $\vs_i=\vs_{\btau i}$ for each $i\in \I$  which satisfies $c_{i, \btau i}=0$. The \emph{reduced Hall algebra associated to $(Q,\btau)$} (or {\em reduced $\imath$Hall algebra}), denoted by $\rMH$, is defined to be the quotient $\Q(\sqq)$-algebra of $\tMH$ by the ideal generated by the central elements
\begin{align}
\label{eqn: reduce1}
[\E_i] +q \vs_i \; (\forall i\in \I \text{ with } i=\btau i), \text{ and }\; [\E_i]*[\E_{\btau i}] -\sqq^{c_{i,\btau i}}\vs_i\vs_{\btau i}\; (\forall i\in \I \text{ with }i\neq \btau i).
\end{align}

\begin{theorem}
[\text{\cite[Proposition 7.5, Theorem 7.7]{LW20a}}]
 \label{thm:Ui=iHall}
Let $(Q, \btau)$ be a virtually acyclic $\imath$quiver. Then there exists a $\Q({\sqq})$-algebra monomorphism
\begin{align}
   \label{eqn:psi morphism}
\widetilde{\Psi}: \tUi_{|v={\sqq}} &\longrightarrow \tMH,
\end{align}
which sends
\begin{align}
B_i \mapsto \frac{-1}{q-1}[S_{i}],\text{ if } i\in\ci,
&\qquad\qquad
\tk_j \mapsto - q^{-1}[\E_j], \text{ if }\btau j=j;
  \label{eq:split}
\\
B_{i} \mapsto \frac{{\sqq}}{q-1}[S_{i}],\text{ if }i\notin \ci,
&\qquad\qquad
\tk_j \mapsto \sqq^{\frac{-c_{j,\btau j}}{2}}[\E_j],\quad \text{ if }\btau j\neq j.
  \label{eq:extra}
\end{align}
Moreover, it induces an embedding $\Psi: \Ui_{|v={\sqq}}\stackrel{\simeq}{\longrightarrow} \rMH$, which sends $B_i$ as in \eqref{eq:split}--\eqref{eq:extra} and  $k_j \mapsto  \vs_{\btau j}^{-1}\sqq^{\frac{-c_{j,\btau j}}{2}}[\E_j], \text{ for } j \in \I\backslash\ci$.
\end{theorem}

Recall $[S_i]_{\ov{p}}^{(m)}$, for $i =\tau i$, defined in \eqref{eq:idividedHallodd}--\eqref{eq:idividedHallev} and $[S_i]^{(m)}$, for $i \neq \tau i$, defined in \eqref{eq:SDP}. For any $m\in\N$, the map $\widetilde{\Psi}$ in Theorem~\ref{thm:Ui=iHall} sends the $\imath$divided powers $B_{i,\ov p}^{(m)}$ in \eqref{eq:iDPodd}--\eqref{eq:iDPev} for $i =\tau i$ (cf. \cite[Lemma 6.3]{LW20a}) and the divided powers $B_{i}^{(m)}$ in \eqref{eq:DP} for $i \neq \tau i$ to
\begin{align}
%\widetilde{\Psi}(B_{i,\odd}^{(m)}) =\frac{[S_i]_{\odd}^{(m)}}{(1-\sqq^2)^m},\quad
 \widetilde{\Psi}(B_{i,\ov p}^{(m)})
 &= \frac{[S_i]_{\ov p}^{(m)}}{(1-\sqq^2)^m} \;\; \text{ for } \ov{p} \in \Z_2, \quad
 \quad \text{ if } i =\tau i,
   \label{eq:iDPpsi} \\
\widetilde{\Psi}(B_{i}^{(m)})
& =
\begin{cases}
\frac{[S_i]^{(m)}}{(1-\sqq^2)^m}     & \text{ for } i \in \I_\btau
\\
\frac{ \sqq^m [S_i]^{(m)}}{(\sqq^2-1)^m} & \text{ for } i \not \in \I_\btau,
  \end{cases}
 \qquad \text{ if } i \neq \tau i.
   \label{eq:DPpsi}
\end{align}

Let $\tCMH$ be the $\Q(\sqq)$-subalgebra (called the {\em composition algebra}) of $\tMHk$ generated by $[S_i]$ and $[\E_i]^{\pm 1}$, for $i\in\I$.

\begin{corollary} [\text{\cite[Corollary 9.9]{LW20a}}]
\label{cor: composition}
Let $(Q,\btau)$ be a virtually acyclic $\imath$quiver. Then there exists an algebra isomorphism:
$%\begin{align*}
\widetilde{\Psi}: \tUi_{|v={\sqq}} \stackrel{\cong}{\longrightarrow} \tCMH
$ %\end{align*}
given by \eqref{eq:split}--\eqref{eq:extra}.
\end{corollary}

\begin{corollary}
\label{cor:braidcomposition}
For any sink $i\in Q_0$, the isomorphism $\Gamma_i: \tMH\stackrel{\cong}{\rightarrow} \tMHl$ in \eqref{thm:Gamma} restricts to an isomorphism
\begin{align*}
\Gamma_i: \tCMH\stackrel{\cong}{\longrightarrow} \cc\widetilde{\ch}(k\bs_iQ,\btau).
\end{align*}
\end{corollary}

Following Ringel, we define a {\em generic composition subalgebra} $\tCMHg$ below. Let $\bfK$  be an infinite set of (nonisomorphic) finite fields, and let us choose for each $\K\in\bfK$ an element $\sqq_\K\in\C$ such that $\sqq_\K^2=|\K|$. Consider the direct product $\tCMHg:=\prod_{\K\in\bfK} \tCMH.$ We view $\tCMHg$ as a $\Q(v)$-module by mapping $v$ to $(\sqq_\K)_\K$.
As in \cite{Rin90,Gr95}, we have the following consequence of Corollary \ref{cor: composition}.

\begin{corollary}  [\text{\cite[Corollary 9.10]{LW20a}}]
\label{cor:isomgeneric}
Let $(Q,\btau)$ be a virtually acyclic $\imath$quiver. Then we have the following algebra isomorphism
$\widetilde{\Psi}: \tUi\longrightarrow \tCMHg$ defined by
\begin{align*}
B_j \mapsto \Big(\frac{-1}{|\K|-1}[S_{j}] \Big)_\K,\text{ if } j\in\ci,
&\qquad\qquad
\tk_i \mapsto \Big(- |\K|^{-1}[\E_i]\Big)_\K, \text{ if }i=\btau i;
  %\label{eq:splitg}
\\
B_{j} \mapsto \Big(\frac{{\sqq_\K}}{|\K|-1}[S_{j}]\Big)_\K,\text{ if }j\notin \ci,
&\qquad\qquad
\tk_i \mapsto \Big(\sqq_\K^{\frac{-c_{i,\btau i}}{2}}[\E_i]\Big)_\K,\quad \text{ if }i\neq \btau i.
 % \label{eq:extrag}
\end{align*}
\end{corollary}

As in Corollary \ref{cor:braidcomposition}, for any sink $i\in Q_0$, the isomorphism $\Gamma_i:\tCMH\stackrel{\cong}{\rightarrow} \cc\widetilde{\ch}(k\bs_iQ,\btau)$ induces an isomorphism
\begin{align}
\label{eq:Gammaigeneric}
\Gamma_i: \tCMHg\stackrel{\cong}{\longrightarrow} \cc\widetilde{\ch}(\bs_iQ,\btau).
\end{align}

\subsection{Symmetries of $\imath$quantum groups}
\label{subsec:BGiQG1}

Let $i$ be a sink of an $\imath$quiver $(Q,\btau)$. Similar to the isomorphism $\widetilde{\Psi}_{Q}$ in \eqref{eqn:psi morphism}, there exists an algebra isomorphism $\widetilde{\Psi}_{Q'}: \tUi_{|v={\sqq}}\rightarrow \tMHi$.

Recall $\ov{\I}_\btau$ from \eqref{eq:Itau2}. we define algebra automorphisms
\[
\TT''_{i,1} \in \Aut (\tUi), \qquad
\text{ for } i \in \ov{\I}_\btau,
\]
by Corollary~\ref{cor:isomgeneric} such that the following diagram commutes:
\begin{align}
  \label{eq:defT}
% \tTT_i = (\widetilde{\Psi}_{Q'})^{-1} \circ \Gamma_i \circ \widetilde{\Psi}_{Q}
\xymatrix{
\tUi  \ar[r]^{\TT''_{i,1}}  \ar[d]^{\widetilde{\Psi}_Q} & \tUi \ar[d]^{\widetilde{\Psi}_{Q'}} \\
\tCMHg \ar[r]^{\Gamma_i} &  \cc\widetilde{\ch}(\bs_iQ,\btau)
}
\end{align}
In other words, specializing at $v=\sqq$ we have the following commutative diagram:
\begin{equation}
  \label{eq:Up=T}
\xymatrix{
\tUi _{|v={\sqq}}\ar[r]^{\TT''_{i,1}}  \ar[d]^{\widetilde{\Psi}_Q} & \tUi _{|v={\sqq}}\ar[d]^{\widetilde{\Psi}_{Q'}}
 \\
\tCMH \ar[r]^{\Gamma_i} &  \cc\widetilde{\ch}(k\bs_iQ,\btau)
}
\end{equation}

%For any $i\in\I$, define
%\begin{align}
%\label{def:b}
%b_i:=\left\{ \begin{array}{ll} -v^2, & \text{ if }i=\btau i,
%\\
%1,&\text{ if }i\neq \btau i . \end{array}  \right.
%\end{align}
%For $\alpha=\sum_{i=1}^n d_i \alpha_i\in \Z\I$, define
%\begin{equation}
%  \label{def:b2}
%b_{\alpha}:= \prod_{i=1}^n b_i^{d_i}.
%\end{equation}

%Then by the isomorphism $\widetilde \psi$ \eqref{eq:extra} and the diagram \eqref{eq:Up=T}, the formula \eqref{eqn:reflection 4} is converted to be
%\begin{align}
%\tTT_i(\tk_\alpha) =\frac{b_{\bs_i\alpha}}{b_\alpha}\tk_{\bs_i \alpha},\quad \forall \alpha\in \Z^\I.
%\end{align}

\begin{theorem}
 \label{thm:BG}
We have a $\Q(v)$-algebra automorphism $\TT''_{i,1}$ of $\tUi$, for $i\in\bar{\I}_\btau$, such that
\begin{enumerate}
\item
$\underline{(i=\btau i)}:$ \;
$\TT''_{i,1}(\tk_j)= (-v^2 \tk_i)^{-c_{ij}}\tk_j$, and
\begin{align*}
\TT''_{i,1}(B_i) &=(-v^{2} \tk_{i} )^{-1}B_i,
\\
\TT''_{i,1}(B_j) &= \sum_{r+s=-c_{ij}}(-1)^{r}v^{r}B_{i,\ov{p}}^{(r)}B_jB_{i,\ov{c_{ij}}+\ov{p}}^{(s)}\\
& +\sum_{u\geq1}\sum_{\stackrel{r+s+2u=-c_{ij}}{
\ov{r}=\ov{p}}}(-1)^{r}v^{r}B_{i,\ov{p}}^{(r)}B_jB_{i,\ov{c_{ij}}+\ov{p}}^{(s)}(-v^2 \tk_i )^u, \qquad \text{ for }j\neq i;
\end{align*}
\item
$\underline{(i\neq \btau i)}:$ \;
$\TT''_{i,1}(\tk_j)= \tk_i^{-c_{ij}} \tk_{\btau i}^{-c_{\tau i,j}} \tk_j$,
\begin{align*}
\TT''_{i,1}(B_j) &=
\begin{cases}  -\tk_{i}^{-1}B_{\btau i},  & \text{ if }j=i \\
-B_i\tk_{\btau i}^{-1}  ,  &\text{ if }j=\btau i,\end{cases}
\end{align*}
and for $j\neq i,\btau i$,
\begin{align*}
\TT''_{i,1}(B_j)
&= \sum^{-\max(c_{ij},c_{\tau i,j})}_{u=0} \; \sum^{-c_{ i,j}-u}_{r=0} \; \sum_{s=0}^{-c_{\tau i,j}-u} (-1)^{r+s} v^{ r-s+(-c_{ij}-r-s-u)u } \\
&\qquad\qquad\qquad\qquad\qquad
 \times B_i^{(r)} B_{\tau i}^{(-c_{\tau i,j}-u-s)} B_j B_{\tau i}^{(s)} B_i^{(-c_{ij}-r-u)}\tk_{\tau i}^u.
\end{align*}
\end{enumerate}
\end{theorem}

\begin{proof}
By Theorem~\ref{thm:Ui=iHall}, Corollary~\ref{cor:isomgeneric}, and the commutative diagram \eqref{eq:defT}, we can transfer the automorphism $\Gamma_i$ and its properties (see Proposition~\ref{prop:reflection}) on $\imath$Hall algebra to an automorphism $\TT''_{i,1}$ on $\tUi$. The formulas for $\TT''_{i,1}$ then follow by \eqref{eq:iDPpsi}--\eqref{eq:DPpsi} and the formulas for $\Gamma_i$ in Theorem~ \ref{thm:Braidsplit} (for $i=\tau i$) and Theorem~ \ref{thm:braidqsplit} (for $i\not =\tau i$), respectively.
\end{proof}

\subsection{(Anti-)involutions on $\tUi$}

For quantum groups $\U$, different variants of braid group symmetries are related to each other via twisting by involutions and anti-involutions on $\U$ \cite[37.2.4]{Lus93}. We now formulate analogous (anti-)involutions in the $\imath$quantum group setting. The following lemma follows by inspection of the defining relations for $\tUi$ in Proposition~\ref{prop:Serre};  cf. \cite{CLW21c}. 

\begin{lemma}
\label{lem:bar}
(a) 
There exists a $\Q$-algebra involution $ \psi_\imath: \tUi\rightarrow \tUi$ (called a bar involution) such that
\[
\psi_\imath(v)=v^{-1}, \quad
\psi_\imath(\tk_i)=v^{c_{i,\btau i}}\tk_{\btau i}, \quad
\psi_\imath (B_i)=B_i,  \quad
\forall i\in \I.
\]
 
(b)
There exists a $\Q(v)$-algebra anti-involution $\sigma_\imath: \tUi\rightarrow \tUi$ such that
\begin{align}
  \label{eq:sigma}
\sigma_\imath(B_i)=B_{i}, \quad \sigma_\imath(\tk_i)= \tk_{\btau i},
%= \left\{\begin{array}{cc}k_{\btau i}, &\text{ if }  2\mid c_{i,\varrho i},\\
%-k_{\btau i}, & \text{ if } 2\nmid c_{i,\varrho i}.\end{array}\right.,
\quad \forall i\in \I.
\end{align}
\end{lemma}

%\begin{proof}
%If suffices to show that $\psi_\imath$ preserves all the defining relations for $\tUi$ in Proposition~\ref{prop:Serre}. Note that if $i=\btau i$, then $c_{i,\tau i}=2$. So $\psi_\imath (\tk_i)=v^2\tk_i$ in this case. So $\psi_\imath$ fixes $B_{i,\ov{p}}^{(n)}$ in \eqref{eq:iDPodd}--\eqref{eq:iDPev}, for any $i\in I$ so that $i=\btau i$, $\ov{p}\in\Z_2$ and $n\in\N$. Hence, one checks readily that $\psi_\imath$ preserves the relations \eqref{relation1}--\eqref{relation6} except \eqref{relation5}. For \eqref{relation5}, as its LHS is clearly fixed by $\psi_\imath$, it suffices to check its RHS is fixed by $\psi_\imath$ as follows:
%\begin{align*}
%&\psi_\imath \left(\frac{1}{v-v^{-1}}\left(v^{c_{i,\btau i}} (v^{-2};v^{-2})_{-c_{i,\btau i}} B_i^{(-c_{i,\btau i})} \tk_i \right. \left. -(v^{2};v^{2})_{-c_{i,\btau i}}B_i^{(-c_{i,\tau i})} \tk_{\btau i}  \right)\right)
%\\
%&= -\frac{1}{v-v^{-1}}\Big(v^{-c_{i,\btau i}} (v^{2};v^{2})_{-c_{i,\btau i}} v^{c_{i,\btau i}} B_i^{(-c_{i,\btau i})}\tk_{\btau i} - (v^{-2};v^{-2})_{-c_{i,\btau i}}v^{c_{i,\btau i}}B_i^{(-c_{i,\tau i})}\tk_i\Big)
%\\
%&= \frac{1}{v-v^{-1}}\Big( v^{c_{i,\btau i}}(v^{-2};v^{-2})_{-c_{i,\btau i}}B_i^{(-c_{i,\tau i})}\tk_i-(v^{2};v^{2})_{-c_{i,\btau i}} B_i^{(-c_{i,\btau i})}\tk_{\btau i} \Big).
%\end{align*}
%The lemma is proved. 
%\end{proof}

It follows by definition that
\begin{equation}   \label{eq:spsp}
\sigma_\imath \psi_\imath =\psi_\imath \sigma_\imath.
\end{equation}

\subsection{Automorphisms $\TT'_{i,e},\TT''_{i,e}$ on $\tUi$}

We shall present 3 more variants of the automorphism $\TT''_{i,1}$ in Theorem~\ref{thm:BG}.

\begin{theorem}
\label{thm:BG1}
For $i\in \bar{\I}_\btau$ and $e \in \{\pm 1\}$,
there are automorphisms $\TT''_{i,e}$ on $\tUi$  such that
\begin{align}
  \label{eq:psiT}
\psi_\imath \TT''_{i,e} \psi_\imath &=\TT''_{i,-e}.
\end{align}
Moreover,
\begin{enumerate}
\item
$\underline{(\btau i= i)}:$ \;
$\TT''_{i,e}(\tk_j)= (-v^{1+e} \tk_i)^{-c_{ij}}\tk_j$, and
\begin{align*}
\TT''_{i,e}(B_i) &= (-v^{1+e}\tk_{i})^{-1}B_i,
\\
\TT''_{i,e}(B_j) &= \sum_{r+s=-c_{ij}}(-1)^{r}v^{er}B_{i,\ov{p}}^{(r)}B_jB_{i,\ov{c_{ij}}+\ov{p}}^{(s)}\\
&+\sum_{u\geq1}\sum_{\stackrel{r+s+2u=-c_{ij}}{
\ov{r}=\ov{p}}}(-1)^{r+u}v^{er+eu}B_{i,\ov{p}}^{(r)}B_jB_{i,\ov{c_{ij}}+\ov{p}}^{(s)}(v\tk_i)^u, \quad \text{ for }j\neq i;
\end{align*}
\item
$\underline{(i\neq \btau i)}:$ \;
$\TT''_{i,e}(\tk_j)= \tk_i^{-c_{ij}} \tk_{\btau i}^{-c_{\tau i,j}} \tk_j$,
\begin{align*}
\TT''_{i,1}(B_j) &=
\begin{cases}
 -\tk_{i}^{-1}B_{\btau i},  & \text{ if }j=i\\
-B_i\tk_{\btau i}^{-1}  ,  &\text{ if }j=\btau i,
\end{cases}
\qquad\qquad
\TT''_{i,-1}(B_j)=
\begin{cases}
 -\tk_{\btau i}^{-1}B_{\btau i},  & \text{ if }j=i\\
-B_i\tk_{i}^{-1},  &\text{ if }j=\btau i,
\end{cases}
\end{align*}
and for $j\neq i,\btau i$,
\begin{align*}
\TT''_{i,1}(B_j)
&= \sum^{-\max(c_{ij},c_{\tau i,j})}_{u=0} \; \sum^{-c_{ i,j}-u}_{r=0} \; \sum_{s=0}^{-c_{\tau i,j}-u} (-1)^{r+s} v^{ r-s+(-c_{ij}-r-s-u)u } \\
&\qquad\qquad\qquad\qquad\qquad\qquad
 \times B_i^{(r)} B_{\tau i}^{(-c_{\tau i,j}-u-s)} B_j B_{\tau i}^{(s)} B_i^{(-c_{ij}-r-u)}\tk_{\tau i}^u,
\\
\TT''_{i,-1}(B_j)
= &\sum^{-\max(c_{ij},c_{\tau i,j})}_{u=0} \; \sum^{-c_{ i,j}-u}_{r=0} \; \sum_{s=0}^{-c_{\tau i,j}-u} (-1)^{r+s} v^{- (r-s+(-c_{ij}-r-s-u)u )} \\
&\qquad\qquad\qquad\qquad\qquad\qquad
 \times B_i^{(r)} B_{\tau i}^{(-c_{\tau i,j}-u-s)} B_j B_{\tau i}^{(s)} B_i^{(-c_{ij}-r-u)}\tk_{i}^u.
\end{align*}
\end{enumerate}
\end{theorem}

\begin{proof}
The formulas for $\TT''_{i,1}$ follows by Theorem~\ref{thm:BG}.
The map $\TT''_{i,-1}$ defined in \eqref{eq:psiT} as the $\psi_\imath$-conjugate of $\TT''_{i,1}$ is clearly an automorphism of $\tUi$. Then the formulas for $\TT''_{i,-1}$ can be verified readily by definition using formulas in Theorem~\ref{thm:BG} and Lemma~\ref{lem:bar}.
\end{proof}

\begin{theorem}
\label{thm:BG2}
For $i\in \bar{\I}_\btau$ and $e \in \{\pm 1\}$, there are automorphisms $\TT'_{i,e}$ on $\tUi$  such that
\begin{align}
 \label{eq:TTsig}
\TT_{i,e}'  &= \sigma_\imath \TT_{i,-e}'' \sigma_\imath,
\qquad
\psi_\imath \TT_{i,e}' \psi_\imath =\TT_{i,-e}'.
\end{align}
Moreover,
\begin{enumerate}
\item
$\underline{(\btau i = i)}:$ \;
$\TT'_{i,e}(\tk_j)=(-v^{1-e} \tk_i)^{-c_{ij}}\tk_j$, and
\begin{align*}
\TT'_{i,e}(B_i) &= (-v^{1-e}\tk_{i})^{-1}B_i,
\\
\TT'_{i,e}(B_j) &= \sum_{r+s=-c_{ij}}(-1)^{r}v^{-er}B_{i,\ov{c_{ij}}+\ov{p}}^{(s)}B_j B_{i,\ov{p}}^{(r)} \\
&+\sum_{u\geq1}\sum_{\stackrel{r+s+2u=-c_{ij}}{
\ov{r}=\ov{p}}}(-1)^{r+u}v^{-er-eu} B_{i,\ov{c_{ij}}+\ov{p}}^{(s)}B_j B_{i,\ov{p}}^{(r)} (v\tk_i)^u, \quad \text{ for }j\neq i.
\end{align*}
\item
$\underline{(i\neq \btau i)}:$ \;
$\TT'_{i,e}(\tk_j)= \tk_i^{-c_{ij}} \tk_{\btau i}^{-c_{\tau i,j}} \tk_j$, and
\begin{align*}
\TT'_{i,-1}(B_j) &=
\begin{cases}
 -B_{\btau i}\tk_{\btau i}^{-1},  & \text{ if }j=i \\
-\tk_{i}^{-1}B_{i},  &\text{ if }j=\btau i,
\end{cases}\qquad\qquad
\TT'_{i,1}(B_j)=
\begin{cases}  -B_{\btau i}\tk_{i}^{-1},  & \text{ if }j=i \\
-\tk_{\btau i}^{-1}B_{i},  &\text{ if }j=\btau i,
\end{cases}
\end{align*}
and for $j\neq i,\btau i$,
\begin{align*}
\TT'_{i,-1}(B_j)
&= \sum^{-\max(c_{ij},c_{\tau i,j})}_{u=0} \; \sum^{-c_{ i,j}-u}_{r=0} \; \sum_{s=0}^{-c_{\tau i,j}-u} (-1)^{r+s} v^{r-s+(-c_{ij}-r-s-u)u } \\
&\qquad\qquad\qquad\qquad\qquad\qquad
 \times  \tk_{i}^u  B_{i}^{(-c_{ij}-r-u)} B_{ \btau i}^{(s)}  B_j B_{\btau i}^{(-c_{\tau i,j}-u-s)} B_{ i}^{(r)},
\\
\TT'_{i,1}(B_j)
= &\sum^{-\max(c_{ij},c_{\tau i,j})}_{u=0} \; \sum^{-c_{ i,j}-u}_{r=0} \; \sum_{s=0}^{-c_{\tau i,j}-u} (-1)^{r+s} v^{- (r-s+(-c_{ij}-r-s-u)u )} \\
&\qquad\qquad\qquad\qquad\qquad\qquad
 \times  \tk_{\btau i}^u B_{i}^{(-c_{ij}-r-u)} B_{\btau i}^{(s)} B_j B_{\btau i}^{(-c_{\tau i,j}-u-s)}  B_{i}^{(r)}.
\end{align*}
\end{enumerate}
\end{theorem}

\begin{proof}
The map $\TT'_{i,e}$ defined in \eqref{eq:TTsig} as the $\sigma_\imath$-conjugate of $\TT''_{i,-e}$ is clearly an automorphism of $\tUi$. The formulas for $\TT'_{i,e}$ acting on the generators of $\tUi$ can be verified readily by definition using formulas in Lemma~\ref{lem:bar} and Theorem~\ref{thm:BG1}. The second formula in \eqref{eq:TTsig} follows from the first formula, \eqref{eq:spsp} and \eqref{eq:psiT}.
\end{proof}

\begin{theorem}
\label{thm:inverse}
We have $\TT_{i,e}' =(\TT_{i,-e}'')^{-1}$, for any $i\in\bar{\I}_\btau$ and $e=\{\pm1\}$.
\end{theorem}

\begin{proof}
It suffices to prove for $e=-1$ by \eqref{eq:TTsig}.
Given an $\imath$quiver $(Q,\btau)$ such that  $i$ is a sink, we have two isomorphisms, $\Gamma_i:\tMH  \stackrel{\cong}{\rightarrow}  \tMHi$ and $\Gamma_i^-:\tMHi\stackrel{\cong}{\rightarrow} \tMH$, which are inverses to each other; see Lemma \ref{lem:Upsilon} and its proof.

By the same arguments as in \S\ref{sec:formulaRF}--\S\ref{sec:formulaRF2}, one can obtain the explicit actions of $\Gamma_i^-$. In particular, similar to \S\ref{subsec:BGiQG1}, we have the following commutative diagram:
\begin{align}
  \label{eq:defT2}
% \tTT_i = (\widetilde{\Psi}_{Q'})^{-1} \circ \Gamma_i \circ \widetilde{\Psi}_{Q}
\xymatrix{
\tUi  \ar[r]^{\TT'_{i,-1}}  \ar[d]^{\widetilde{\Psi}_{Q'}} & \tUi \ar[d]^{\widetilde{\Psi}_{Q}} \\
 \ar[r]^{\Gamma_i^-}  \cc\widetilde{\ch}(\bs_iQ,\btau) & \tCMHg
}
\end{align}
Combining with \eqref{eq:defT}, we conclude that $\TT_{i,-1}'$ and $\TT_{i,1}''$ are inverses to each other.
\end{proof}

\begin{remark}
The results in Theorems~\ref{thm:BG1}--\ref{thm:BG2} verify substantially \cite[Conjecture~ 6.5]{CLW21} in case $i=\tau i$ and  \cite[Conjecture~3.7]{CLW21c} in case $i\not =\tau i$, for quasi-split $\imath$quantum groups $\tUi$ associated to symmetric generalized Cartan matrices with all $c_{i,\tau i}$ even. Actually, the notations $\TT'_{i,e}$ and $\TT''_{i,e}$ used in Theorem~\ref{thm:BG1}--\ref{thm:BG2} are swapped from the corresponding notations $\TT''_{i,e}$ and $\TT'_{i,e}$ used in \cite[Conjecture~ 6.5]{CLW21} for the following reason. The leading terms (i.e., the $u=0$ summands) in the formulas for the symmetries $\TT'_{i,e}, \TT''_{i,e}$ acting on $B_j$ (for $j\neq i=\tau i$) in Theorem~\ref{thm:BG1}--\ref{thm:BG2} are precisely the formulas for Lusztig's symmetries $T'_{i,e}, T''_{i,e}$ on $F_j$; hence this is most compatible with the view that $F_j$ (not $E_j$) is a leading term for $B_j$.  

The formulas for the automorphisms $\TT'_{i,e}$, $\TT''_{i,e}$ remain valid for quasi-split symmetrizable Kac-Moody type, when $v$ is replaced by $v_i =v^{\frac{i\cdot i}2}$ in the formulas above and also in the $\imath$divided powers. They can be verified in the framework of  $\imath$Hall algebras associated to valued $\imath$quivers (which are to be developed).
\end{remark}

\begin{remark}
The operators $\TT''_{i,e}$, for $i\in\bar{\I}_\btau$, are expected to satisfy the braid relations of the restricted Weyl group $W^{\btau}$ defined in \eqref{eq:Wtau} (extending the suggestion in \cite{KP11} for finite type). For $\tUi$ of rank 2, the only nontrivial braid relation appears in finite type, and hence it holds by \cite{LW21}.  If one can show that $\TT''_{i,e}$ acts on module level in a way compatible with its action on $\tUi$, then the braid relation follows from the rank 2 results. Other than that, for various simply-laced locally finite Kac-Moody types (including all affine types), the braid relations hold thanks again to the computations in \cite{LW21}. %But new ideas are needed to handle the general cases.
\end{remark}

\subsection{Symmetries of $\Ui$}

Set $\bvsd =(\vs_{\diamond,i})_{i\in \I}$, where
\begin{equation}
  \label{eq:bvsd}
\vs_{\diamond,i}=-v^{-2} \;  \text{ if } i= \btau i,
\qquad
\vs_{\diamond,i}=v^{\frac{-c_{i,\btau i}}{2}}\;  \text{ if } i\neq \btau i.
\end{equation}
 We call the parameters $\bvsd$ %in \ref{eq:special split}--\eqref{eq:specialE}
 {\em distinguished}, for the corresponding $\imath$quantum groups $\tUi$ and $\Ui$.

We shall use the index $\bvsd$ to indicate the algebras under consideration are associated to the distinguished parameters $\bvsd$. Let $\La^\imath$ be the $\imath$quiver algebra associated to an $\imath$quiver $(Q, \btau)$. Recall the reduced $\imath$Hall algebra $\rMH$ with a general parameter $\bvs$; cf. \eqref{eqn: reduce1}. If follows that the reduced $\imath$Hall algebra $\rMHd$ with the distinguished parameter $\bvsd$  is  the quotient algebra of $\tMH$ by the ideal generated by
\begin{align}
\label{ideal of MRH}
[\E_i]-1 \; (i\in \I \text{ with } i=\btau i), \quad \text{ and } \quad 
[\E_i]*[\E_{\btau i}]-1\; (i\in \I \text{ with }i\neq \btau i).
\end{align}

\begin{proposition}
[\text{\cite[Proposition 7.1]{LW21}}]
   \label{prop:derived invariant of QSP}
For any sink $i\in Q_0$, the isomorphism $\Gamma_i$ induces an isomorphism of algebras
$\bar{\Gamma}_i: \rMHd \xrightarrow{\cong} \rMHdi$.
\end{proposition}

Recall from Theorem \ref{thm:Ui=iHall} that the isomorphism $\Psi_{Q}: \Ui_{\bvsd |v={\sqq}} \rightarrow \rMHd$ sends
\begin{align*}
B_j \mapsto \frac{-1}{q-1}[S_{j}],\text{ for } j\in\ci;
&\qquad\qquad
B_{j} \mapsto \frac{\sqq}{q-1}[S_{j}],\text{ for } j \in \I\setminus \ci;
 \notag \\
k_j \mapsto [\E_j],\quad \text{ for }  j \in \I\setminus \ci.
\notag
\end{align*}
Similarly, there exists an isomorphism of algebras $\Psi_{Q'} : \Ui_{\bvsd |v={\sqq}} \longrightarrow \rMHdi$ (where the $[S_{j}], [\E_j]$ above are replaced by $[S'_{j}], [\E'_j]$).

Using the same argument of \cite[Proposition 7.2]{LW21}, we obtain the following.
\begin{proposition}
   \label{prop:braidDist}
Let $(Q, \btau)$ be an $\imath$quiver. Then $\TT''_{i,1}: \tUi_{\bvsd} \rightarrow \tUi_{\bvsd}$ in \eqref{eq:defT} induces an automorphism $\TT''_{i,1}$ on $\Ui_{\bvsd}$, for each $i\in\bar{\I}_\btau$.
Moreover, for any sink $i$ in $Q_0$, we have the following commutative diagram of isomorphisms:
\[
\xymatrix{ \Ui_{\bvsd |v={\sqq}}\ar[r]^{\TT''_{i,1}}  \ar[d]^{\Psi_Q} & \Ui_{\bvsd |v={\sqq}}\ar[d]^{\Psi_{Q'}} \\
\cc\rMHd \ar[r]^{\bar{\Gamma}_i} & \cc\rMHdi }
\]
\end{proposition}

The explicit formulas for $\TT''_{i,1}: \Ui_{\bvsd} \rightarrow \Ui_{\bvsd}$ are given as follows:
\begin{enumerate}
\item
$\underline{(i=\btau i)}:$ \;
$\TT''_{i,1}(k_j)= k_j$, and
\begin{align*}
\TT''_{i,1}(B_i) &=B_i,
\\
\TT''_{i,1}(B_j) &= \sum_{r+s=-c_{ij}}(-1)^{r}v^{r}B_{i,\ov{p}}^{(r)}B_jB_{i,\ov{c_{ij}}+\ov{p}}^{(s)}\\
& +\sum_{u\geq1}\sum_{\stackrel{r+s+2u=-c_{ij}}{
\ov{r}=\ov{p}}}(-1)^{r}v^{r}B_{i,\ov{p}}^{(r)}B_jB_{i,\ov{c_{ij}}+\ov{p}}^{(s)}, \qquad \text{ for }j\neq i;
\end{align*}
\item
$\underline{(i\neq \btau i)}:$ \;
$\TT''_{i,1}(k_j)=  k_{\btau i}^{c_{ij}-c_{\tau i,j}} k_j$, and
\begin{align*}
\TT''_{i,1}(B_j) &=
\begin{cases}  -k_{\btau i} B_{\btau i},  & \text{ if }j=i,\\
-B_ik_{\btau i}^{-1}  ,  &\text{ if }j=\btau i,\end{cases}
\end{align*}
and for $j\neq i,\btau i$,
\begin{align*}
\TT''_{i,1}(B_j)
&= \sum^{-\max(c_{ij},c_{\tau i,j})}_{u=0} \; \sum^{-c_{ i,j}-u}_{r=0} \; \sum_{s=0}^{-c_{\tau i,j}-u} (-1)^{r+s} v^{\big(r-s+(-c_{ij}-r-s-u)u\big)} \\
&\qquad\qquad\qquad\qquad  B_i^{(r)} B_{\tau i}^{(-c_{\tau i,j}-u-s)} B_j B_{\tau i}^{(s)} B_i^{(-c_{ij}-r-u)}k_{\tau i}^u.
\end{align*}
\end{enumerate}

Below we write $\Ui=\Ui_{\bvs}$ to indicate its dependence on a parameter $\bvs$.
It is well known that the $\Q(v)$-algebras $\Ui_{\bvs}$ (up to some field extension) are isomorphic for different choices of  parameters $\bvs$  \cite{Let02}; see Lemma~\ref{lem:base change} below. %We shall use the index ${\bvs_{\diamond}}$ to indicate the relevant algebras, e.g., $\Ui_{{\bvsd}}$  in distinguished parameters $\bvsd=(\vs_{\diamond,i})_{i\in \I}$ in \eqref{eq:special split}--\eqref{eq:specialE}.

Consider a field extension of $\Q(v)$
\begin{align}
\label{def:ai}
{\F}= \Q(v)(a_i\mid i\in I),
\qquad
\text{ where }a_i=\sqrt{\frac{\vs_{\diamond,i}}{\vs_i}} \quad (i\in \I).
\end{align}
Denote by
\begin{align}
\label{def:basechange}
{}_{\F}\Ui_{\bvs} =\F \otimes_{\Q(v)} \Ui_{\bvs}
\end{align}
 the $\F$-algebra obtained by a base change. By a direct computation a rescaling automorphism on ${}_{\F}\U$ induces an isomorphism in the lemma below.

\begin{lemma}
[\text{\cite[Lemma 8.6]{LW21}}]
  \label{lem:base change}
There exists an isomorphism of ${\mathbb F}$-algebras
\begin{align*}
\phi_{\bf u}: {}_{\F}\Ui_{\bvs_{\diamond}} & \longrightarrow {}_{\F}\Ui_{\bvs}
\\
B_i \mapsto a_iB_i,\qquad k_j & \mapsto k_j, \quad (\forall i\in \I, j\in\I\backslash \ci).
\end{align*}
\end{lemma}

%As before the notation $\Ui_{\bvs}$ stands for the $\imath$quantum group with general parameters $\bvs$. 
Below we shall denote by $\TT_{\diamond,i}:\Ui_{\bvsd} \rightarrow \Ui_{\bvsd}$ the isomorphism $\TT''_{i,1}: \Ui_{\bvsd} \rightarrow \Ui_{\bvsd}$ obtained in Proposition \ref{prop:braidDist}.
Let $(Q, \btau)$ be an $\imath$quiver.
We now define a braid group action $\TT''_{i,1}$ on ${}_\F\Ui$  from the $\TT_{\diamond,i}$ on $\Ui_{\bvsd}$ via a conjugation by the isomorphism $\phi_{\bf u}$:
\begin{equation}
  \label{eq:TTTT}
\TT''_{i,1} = \phi_{\bf u} \TT_{\diamond,i} \phi_{\bf u}^{-1}.
\end{equation}

\begin{proposition}
  \label{thm:braidUigeneral}
  Let $(Q, \btau)$ be an $\imath$quiver.
Then there is an automorphism $\TT''_{i,1}$ on ${}_\F\Ui_{\bvs}$, for each $i\in \bar{\I}_\btau$. %Moreover there exists a homomorphism
%$\brW \rightarrow \aut({}_\F\Ui)$ such that $t_i\mapsto \TT_i$ in \eqref{eq:TTTT}, for all $i\in \ci$.
\end{proposition}

\begin{remark}
By using the same argument, one can construct automorphisms $\TT''_{i,-1}$, $\TT'_{i,\pm1}$ for the $\imath$quantum group $\Ui_{\bvs}$ with general parameters $\bvs$ by considering {\em different} distinguished parameters. In fact,
for $\TT''_{i,1}$, $\TT'_{i,-1}$, the distinguished parameters are as shown in \eqref{eq:bvsd};
for $\TT''_{i,-1}$, $\TT'_{i,1}$,  the distinguished parameters  are chosen to be
 \begin{align*}
\vs_{\diamond,i}=-1 \;  \text{ if } i= \btau i,
\qquad
\vs_{\diamond,i}=v^{\frac{-c_{i,\btau i}}{2}}\;  \text{ if } i\neq \btau i.
\end{align*}
However, the relation \eqref{eq:psiT} no longer holds in the setting of $\Ui$. In fact, there are some strong constraints on the parameters $\bvs$ (different from the distinguished parameters in \eqref{eq:bvsd}) for the existence of the bar involution $\psi_\imath$ on $\Ui$; see \cite[Proposition 3.7]{CLW18}.
\end{remark}

\subsection{Symmetries of the Drinfeld double $\tU$}

Recall that the Drinfeld double $\tU$ is the universal $\imath$quantum group of diagonal type; see \cite[Lemma 8.3]{LW19}. Below we shall write down explicitly the automorphisms of the Drinfeld double $\tU$ following Theorems~\ref{thm:BG1}--\ref{thm:BG2}, for the convenience of the reader.

There exists a $\Q(v)$-algebra anti-involution $\sigma:\tU\rightarrow \tU$ such that
\begin{align*}
\sigma(E_i)=E_{i},\quad \sigma(F_i)=F_i,\quad \sigma(\tK_i)= \tK'_{i},
\quad \forall i\in \I.
\end{align*}
%(b)
There exists a $\Q$-algebra automorphism $ \psi: \tU\rightarrow \tU$ (called bar involution) such that
\[
\psi(v)=v^{-1}, \quad
\psi(\tK_i)=\tK'_{i}, \quad
\psi (E_i)=E_i,  \quad \psi(F_i)=F_i,\quad
\forall i\in \I.
\]

\begin{proposition}
\label{prop:BG1U}
For $i\in \I$ and $e \in \{\pm 1\}$,
there are automorphisms $\TT''_{i,e}$ on $\tU$  such that
\begin{align*}
\psi \TT''_{i,e} \psi &=\TT''_{i,-e}.
\end{align*}
Moreover, we have
\begin{align*}
\TT''_{i,e}(\tK_j)&= \tK_i^{-c_{ij}} \tK_j,\qquad \TT''_{i,e}(\tK'_j)= (\tK'_i)^{-c_{ij}} \tK'_j,
\\
\TT''_{i,1}(F_i) &=  -\tK_{i}^{-1}E_{ i}, \qquad \TT''_{i,1}(E_i)= -F_i(\tK'_{i})^{-1},
\\
\TT''_{i,-1}(F_i) &=  -(\tK'_{i})^{-1}E_{ i}, \qquad \TT''_{i,-1}(E_i)= -F_i\tK_{i}^{-1},
\\
\TT''_{i,-e}(F_j)
&=  \sum^{-c_{ i,j}}_{r=0} (-1)^{r} v^{ -er} F_i^{(r)}  F_j F_i^{(-c_{ij}-r)},
\\
\TT''_{i,-e}(E_j)
= & \sum^{-c_{ i,j}}_{r=0}  (-1)^{r} v^{er} E_{i}^{(-c_{\tau i,j}-r)} E_j E_{i}^{(r)},\qquad \text{ for } j\neq i.
\end{align*}
\end{proposition}

\begin{proposition}
\label{prop:BG2U}
For $i\in \I$ and $e \in \{\pm 1\}$, there are automorphisms $\TT'_{i,e}$ on $\tU$  such that
\begin{align*}
\sigma\TT_{i,e}' \sigma &= \TT_{i,-e}'',
\qquad
\psi \TT_{i,e}' \psi=\TT_{i,-e}'.
\end{align*}
Moreover, we have
\begin{align*}
\TT'_{i,e}(\tK_j)&= \tK_i^{-c_{ij}} \tK_j,\qquad \TT'_{i,e}(\tK'_j)= (\tK'_i)^{-c_{ij}} \tK'_j,
\\
\TT'_{i,-1}(F_i) &=  -E_i(\tK'_{i})^{-1},  \qquad \TT'_{i,-1}(E_i) =  -\tK_{i}^{-1}F_i,
\\
\TT'_{i,1}(F_i) &=  -E_i\tK_{i}^{-1},  \qquad \TT'_{i,1}(E_i) =  -(\tK'_{i})^{-1}F_i,
\\
\TT'_{i,e}(F_j)
&=  \sum^{-c_{ i,j}}_{r=0} \; (-1)^{r} v^{-er }  F_{i}^{(-c_{ij}-r)} F_j F_{ i}^{(r)},
\\
\TT'_{i,e}(E_j)
&= \sum^{-c_{ i,j}}_{r=0} \; (-1)^{r} v^{er }   E_{ i}^{(r)} E_j E_{i}^{(-c_{ij}-r)},\qquad \text{ for } j\neq i.
\end{align*}
\end{proposition}

\begin{remark}
The actions of $\TT'_{i,e}, \TT''_{i,e}$ on $\tU$ factor through the quotient $\U =\tU \big/ ( \tK_i' \tK_i-1\mid i\in\I )$ to the corresponding automorphisms on $\U$, and the formulas in Propositions~\ref{prop:BG1U}--\ref{prop:BG2U} are then reduced to Lusztig's formulas \cite[\S37.1.3]{Lus93} upon the identification $\tK_i' =\tK_i^{-1}$.
\end{remark}

%%%%%%%%%%%%%%
%%%%%%%%%%%%%%

\appendix

%%%%%%
%%%%%%
\section{Proof of Proposition \ref{prop:build-block}}
\label{app:prop}

This appendix is devoted to a proof of Proposition \ref{prop:build-block}, which concerns about the computation of $[(S'_i)^{\oplus m_1}\oplus (S'_{\btau i})^{\oplus n_1}]*[S'_j]* [(S'_i)^{\oplus m_2} \oplus (S'_{\btau i})^{\oplus n_2}]$.
\subsection{The setup}

By definition, we have
\begin{align*}
&[S_i^{\oplus m_1}\oplus S_{\btau i}^{\oplus n_1}]*[S_{j}]* [S_i^{\oplus m_2} \oplus S_{\btau i}^{\oplus n_2}]
= [S_i^{\oplus m_1}\oplus S_{\btau i}^{\oplus n_1}]*[S_{j}\oplus S_i^{\oplus m_2} \oplus S_{\btau i}^{\oplus n_2}]\\
& =\sqq^{ c_{ij}m_1+c_{\btau i,j}n_1 -m_1m_2-n_1n_2 }\sum_{[L]\in\Iso(\mod(\Lambda^\imath))}\big|\Ext^1_{\Lambda^\imath}( S_i^{\oplus m_1}\oplus S_{\btau i}^{\oplus n_1}, S_{j}\oplus  S_i^{\oplus m_2} \oplus S_{\btau i}^{\oplus n_2})_L\big| \cdot [L].
\end{align*}

For any $[L]\in \Iso(\mod(\Lambda^\imath))$ such that
\begin{align}
\label{eqn:shortexact L}
\big|\Ext^1_{\Lambda^\imath}( S_i^{\oplus m_1}\oplus S_{\btau i}^{\oplus n_1}, S_{j}\oplus  S_i^{\oplus m_2} \oplus S_{\btau i}^{\oplus n_2})_L\big|\neq0,
\end{align}
there exist a unique $[M]\in\Iso(\mod(kQ))$ and $d,e\in\N$ such that $[L]=[M\oplus \E_i^{\oplus d} \oplus \E_{\btau i}^{\oplus e}]$ in $\tMH$.
In this case, $M$ admits the following exact sequence
$$0\longrightarrow S_{j}\oplus S_i^{\oplus (m_2-e)}\oplus S_{\btau i}^{\oplus (n_2-d)} \longrightarrow M\longrightarrow S_i^{\oplus (m_1-d)} \oplus S_{\btau i}^{\oplus(n_1-e)} \longrightarrow0.$$

Fix $[M]\in\Iso(\mod(kQ))$, $0\leq d\leq\min(n_2,m_1)$ and $0\leq e\leq \min(n_1,m_2)$. Then there exists a unique indecomposable $kQ$-module $N$ such that $M\cong N\oplus S_i^{\oplus t_1^M}\oplus S_{\btau i}^{\oplus t_3^M}$ for some unique $t_1^M,t_3^M$.
Denote by
\begin{align*}
\cc_M:=\{[\xi]&\in\Ext^1_{\Lambda^\imath}(  S_i^{\oplus m_1}\oplus S_{\btau i}^{\oplus n_1}, S_{j}\oplus  S_i^{\oplus m_2} \oplus S_{\btau i}^{\oplus n_2})_L\mid L\text{ admits }\\
 &\qquad\qquad\text{a short exact sequence } 0\rightarrow M\rightarrow L\rightarrow \E_i^{\oplus d}\oplus \E_{\btau i}^{\oplus e}\rightarrow0\}.
\end{align*}
In this way, we have
\begin{align}
\label{eq:builtb1}
&[S_i^{\oplus m_1}\oplus S_{\btau i}^{\oplus n_1}]*[S_{j}]* [S_i^{\oplus m_2} \oplus S_{\btau i}^{\oplus n_2}]\\\notag
%=&[S_i^{\oplus m_1}\oplus S_{\btau i}^{\oplus n_1}]*[S_{j}\oplus S_i^{\oplus m_2} \oplus S_{\btau i}^{\oplus n_2}]\\
=&\sum_{e=0}^{\min( n_1,m_2)}\sum_{d=0}^{\min(n_2,m_1 )}\sum_{[M]\in\mathcal{I}'_{m_1+m_2-d-e,n_1+n_2-d-e}} \sqq^{ c_{ij}m_1+c_{\btau i,j}n_1 -m_1m_2-n_1n_2 } |\cc_M| \cdot [M\oplus \E_i^{\oplus d}\oplus \E_{\btau i}^{\oplus e}]\\\notag
=&\sum_{e=0}^{\min( n_1,m_2)}\sum_{d=0}^{\min(n_2,m_1 )}\sum_{[M]\in\mathcal{I}'_{m_1+m_2-d-e,n_1+n_2-d-e}} \sqq^{ c_{ij}m_1+c_{\btau i,j}n_1 -m_1m_2-n_1n_2 +(e-d)(m_1+m_2-n_1-n_2) }  |\cc_M| \\\notag
&\qquad\qquad\qquad \qquad\qquad\qquad   \times [M]*[\E_i]^{ d}* [\E_{\btau i}]^{e}.
\end{align}

\subsection{Computation of $ |\cc_M|$}

We shall compute $ |\cc_M|$.
Let $Q'$ and $Q''$ be the full subquivers of $Q$ formed by the vertices $i,j$ and the vertices $\btau i,j$ respectively. Then we have two restriction functors $\res_{ij}:\mod(kQ)\rightarrow \mod(kQ')$ and $\res_{\btau i,j}: \mod(kQ)\rightarrow \mod(kQ'')$. Set
\begin{equation}
  \label{eq:MN}
M_1:=\res_{ij}(M), \quad
N_1:=\res_{ij}(N), \quad
M_2:=\res_{\btau i,j}(M), \quad
N_2:=\res_{\btau i,j}(N).
\end{equation}
Denote by
\begin{align*}
\cc_1:=&\{[\xi]\in\Ext^1_{\Lambda^\imath}(  S_i^{\oplus m_1}, S_{j}\oplus  S_i^{\oplus m_2} \oplus S_{\btau i}^{\oplus n_2})_{L_1}\mid L_1\text{ admits a short }\\
 &\qquad\qquad\text{ exact sequence } 0\rightarrow M_1\rightarrow L_1\rightarrow \E_i^{\oplus d}\oplus S_i^{\oplus e} \oplus S_{\btau i}^{\oplus (n_2-d)} \rightarrow0\},
 \\
\cc_2:=&\{[\xi]\in\Ext^1_{\Lambda^\imath}(   S_{\btau i}^{\oplus n_1}, S_{j}\oplus  S_i^{\oplus m_2} \oplus S_{\btau i}^{\oplus n_2})_{L_2}\mid L_2\text{ admits a short }\\
 &\qquad\qquad\text{ exact sequence } 0\rightarrow M_2\rightarrow L_2\rightarrow \E_{\btau i}^{\oplus e}\oplus S_{\btau i}^{\oplus d} \oplus S_i^{\oplus (m_2-e)}\rightarrow0\}.
\end{align*}

\begin{lemma}
\label{lem:decom2}
Retain the notations as above.
Then $|\cc_M|=|\cc_1|\cdot|\cc_2|$.
\end{lemma}

\begin{proof}
Applying $\Hom_{\Lambda^\imath}(-,S_{j}\oplus S_i^{\oplus m_2}\oplus S_{\btau i}^{\oplus n_2})$ to the split short exact sequence
\begin{align}
\label{eqn:split1}
0\longrightarrow S_i^{\oplus m_1} \longrightarrow S_i^{\oplus m_1}\oplus S_{\btau i}^{\oplus n_1}\longrightarrow S_{\btau i}^{\oplus n_1}\longrightarrow 0,
\end{align}
we have the following short exact sequence
\begin{align}
\label{eqn:ses1}
0\longrightarrow \Ext^1(S_{\btau i}^{\oplus n_1},S_{j}\oplus S_i^{\oplus m_2}\oplus S_{\btau i}^{\oplus n_2})\stackrel{\beta}{\longrightarrow} \Ext^1( S_i^{\oplus m_1}\oplus S_{\btau i}^{\oplus n_1},S_{j}\oplus S_i^{\oplus m_2}\oplus S_{\btau i}^{\oplus n_2})\\
\stackrel{\alpha}{\longrightarrow}\Ext^1( S_i^{\oplus m_1},S_{j}\oplus S_i^{\oplus m_2}\oplus S_{\btau i}^{\oplus n_2})\longrightarrow 0.\notag
\end{align}
Then $\alpha(\cc_M)=\cc_1$, and thus,
$\cc_M= \bigsqcup_{[\xi]\in \cc_1} (\alpha|_{\cc_M})^{-1}([\xi])$.

Since \eqref{eqn:split1} is split, we have the following exact sequence
\begin{align}
\label{eqn:ses1}
0\longrightarrow \Ext^1(S_i^{\oplus m_1},S_{j}\oplus S_i^{\oplus m_2}\oplus S_{\btau i}^{\oplus n_2})\stackrel{\delta}{\longrightarrow} \Ext^1( S_i^{\oplus m_1}\oplus S_{\btau i}^{\oplus n_1},S_{j}\oplus S_i^{\oplus m_2}\oplus S_{\btau i}^{\oplus n_2})\\
\stackrel{\gamma}{\longrightarrow} \Ext^1( S_{\btau i}^{\oplus n_1},S_{j}\oplus S_i^{\oplus m_2}\oplus S_{\btau i}^{\oplus n_2})\longrightarrow 0\notag
\end{align}
such that $\gamma \circ \beta=\Id$ and $\alpha\circ \delta=\Id$.
For any $[\eta_1],[\eta_2]\in \cc_M$ such that $\alpha([\eta_1])=\alpha([\eta_2])=[\xi]$, if $\gamma([\eta_1])=\gamma([\eta_2])$, then there exists a unique $[\xi']\in\Ext^1(S_i^{\oplus m_1},S_{j}\oplus S_i^{\oplus m_2}\oplus S_{\btau i}^{\oplus n_2})$ such that $\delta([\xi'])=[\eta_1] -[\eta_2]$. It follows that $[\xi']=\alpha\circ\delta([\xi'])= \alpha([\eta_1])-\alpha([\eta_2])=[0]$. So $\delta([\xi'])=0$ and $[\eta_1]=[\eta_2]$.
Therefore, we obtain that  $\gamma|_{(\alpha|_{\cc_M})^{-1} ([\xi])}$ is injective for any $[\xi]\in \cc_1$.

By the above calculations,
\begin{align}
\label{eqn: disjoint}
\cc_M=\bigsqcup_{[\xi]\in \cc_1} \cc_\xi, \text{ where }\cc_\xi =\gamma(\alpha|_{\cc_M})^{-1} ([\xi]).
\end{align}
%where $\cc_\xi =\gamma(\alpha|_{\cc_M})^{-1} ([\xi])$.
One can show that $\cc_\xi=\cc_2$, which is independent of $[\xi]\in \cc_1$. Therefore,
$|\cc_M|=|\cc_1|\cdot|\cc_2|$.
\end{proof}

\subsection{Computation of $ |\cc_1|$ and $ |\cc_2|$}

It remains to compute $ |\cc_1|$ and $ |\cc_2|$.
\begin{lemma}
\label{lem:decom3}
Retain the notations as above. Then
\begin{align*}
|\cc_1|=& q^{-c_{ij}d} \big|\Ext^1_{kQ}(S_i^{\oplus (m_1-d)}, S_{j})_{N_1\oplus S_i^{\oplus (t_1^M+e-m_2)}}\big|\cdot \big|\Ext_{\Lambda^\imath}^1 (S_i^{\oplus m_1},S_{\btau i}^{\oplus n_2})_{\E_i^{\oplus d} \oplus S_i^{\oplus (m_1-d)} \oplus S_{\btau i}^{\oplus (n_2-d)} }\big|,
%=&q^{ad} |\Ext^1_{kQ}(S_i^{\oplus (m_1-d)}, S_{j}\oplus S_i^{\oplus (m_2-e)})_{M_1}|\cdot |\Ext_{\Lambda^\imath}^1 (S_i^{\oplus m_1},S_{\btau i}^{\oplus n_2})_{\E_i^{\oplus d} \oplus S_i^{\oplus (m_1-d)} \oplus S_{\btau i}^{\oplus (n_2-d)} }|.
\\
|\cc_2|=& q^{-c_{\btau i,j}e} \big|\Ext^1_{kQ}(S_{\btau i}^{\oplus (n_1-e)}, S_{j})_{N_2\oplus S_{\btau i}^{\oplus( t_3^M+d-n_2 )}}\big|\cdot \big|\Ext_{\Lambda^\imath}^1 (S_{\btau i}^{\oplus n_1},S_i^{\oplus m_2})_{\E_{\btau i}^{\oplus e} \oplus S_i^{\oplus (m_2-e)} \oplus S_{\btau i}^{\oplus (n_1-e)} }\big|.
\end{align*}
\end{lemma}

\begin{proof}
The 2 formulas are equivalent, and we only prove the first one.
Let $M'_1:=N_1\oplus S_i^{\oplus(t_1^M-m_2)}$, and
\begin{align*}
\cc'_1:=&\{[\xi]\in\Ext^1_{\Lambda^\imath}(  S_i^{\oplus m_1}, S_{j} \oplus S_{\btau i}^{\oplus n_2})_{L'_1}\mid L'_1\text{ admits a short }\\
 &\qquad\qquad\text{ exact sequence } 0\rightarrow M'_1\rightarrow L'_1\rightarrow \E_i^{\oplus d}\oplus S_i^{\oplus e} \oplus S_{\btau i}^{\oplus (n_2-d)} \rightarrow0\}.\\
\end{align*}
Consider the $\imath$quiver algebra with its quiver as in the right figure of \eqref{diag:split KM} (the number of arrows from $i$ to $j$ is $ a=-c_{ij}$), which is denoted by ${}^s\Lambda^\imath$ to avoid confusions.
Then any $\bs_i\Lambda^\imath$-module $L=(L_k,L_\alpha,L_{\varepsilon_k})$ supported at $i,\btau i$ and $j$ with $L_{\alpha}=0$ for any $\alpha:\btau i\rightarrow j$ can be viewed as a ${}^s\Lambda^\imath$-module $G(L)$, that is,
$G(L)_i:=L_i\oplus L_{\btau i}$, $G(L)_j:=L_j$. Let
\begin{align*}
\cc''_1:=&\{[\xi]\in\Ext^1_{^s\Lambda^\imath}(  S_i^{\oplus m_1}, S_{j} \oplus S_{i}^{\oplus n_2})_{L''_1}\mid L''_1\text{ admits a short }\\
 &\qquad\qquad\text{ exact sequence } 0\rightarrow G(M'_1)\rightarrow L''_1\rightarrow \E_i^{\oplus d}\oplus S_i^{\oplus (n_2-d+e)} \rightarrow0\}.\\
\end{align*}
Then $|\cc'_1|=|\cc''_1|$.
By applying \cite[Proposition 3.10]{LW20a} to compute $[S_i^{\oplus m_1}]\diamond[ S_{j} \oplus S_{i}^{\oplus n_2}]$ in $\cs\cd\ch(^s\Lambda^\imath)$, one obtains that
\begin{align*}
|\cc''_1|=& q^{-c_{ij}d} \big|\Ext^1_{kQ^s}(S_i^{\oplus (m_1-d)}, S_{j})_{G(N_1)\oplus S_i^{\oplus (t_1^M+e-m_2)}}\big|\cdot \big|\Ext_{^s\Lambda^\imath}^1 (S_i^{\oplus m_1},S_{ i}^{\oplus n_2})_{\E_i^{\oplus d} \oplus S_i^{\oplus (m_1+n_1-2d)} }\big|.
\end{align*}
Here $Q^s$ is the quiver $\xymatrix{i\ar[rr]|-{-c_{ij}} && j}$.
Clearly, we have
\begin{align*}
\big|\Ext^1_{kQ^s}(S_i^{\oplus (m_1-d)}, S_{j})_{G(N_1)\oplus S_i^{\oplus (t_1^M+e-m_2)}}\big|
&= \big|\Ext^1_{kQ}(S_i^{\oplus (m_1-d)}, S_{j})_{N_1\oplus S_i^{\oplus (t_1^M+e-m_2)}}\big|,
\\
\big|\Ext_{^s\Lambda^\imath}^1 (S_i^{\oplus m_1},S_{ i}^{\oplus n_2})_{\E_i^{\oplus d} \oplus S_i^{\oplus (m_1+n_1-2d)} }\big|
&=\big|\Ext_{\Lambda^\imath}^1 (S_i^{\oplus m_1},S_{\btau i}^{\oplus n_2})_{\E_i^{\oplus d} \oplus S_i^{\oplus (m_1-d)} \oplus S_{\btau i}^{\oplus (n_2-d)} }\big|,
\end{align*}
and then the desired formula follows.
\end{proof}

Recall the notation $N_1, N_2$ from \eqref{eq:MN}.

\begin{lemma}
\label{lem:decom1}
Retain the notations as above. Then
\begin{align*}
&\big|\Ext^1_{kQ}(S_i^{\oplus (m_1-d)}\oplus S_{\btau i}^{\oplus (n_1-e)}, S_{j})_{N\oplus S_i^{\oplus (t_1^M+e-m_2)} \oplus S_{\btau i}^{\oplus(t_3^M+d-n_2)}} \big|
\\
=&\big|\Ext^1_{kQ}(S_i^{\oplus (m_1-d)}, S_{j})_{N_1\oplus S_i^{\oplus (t_1^M+e-m_2)}}\big|\cdot \big|\Ext^1_{kQ}(S_{\btau i}^{\oplus (n_1-e)}, S_{j})_{N_2\oplus S_{\btau i}^{\oplus( t_3^M+d-n_2 )}}\big|.
\end{align*}
\end{lemma}

\begin{proof}
Note that $N$, $N_1$ and $N_2$ are indecomposable. Then the orders of their automorphism groups are $q-1$.
A direct computation shows that
\begin{align*}
F_{S_i^{\oplus (m_1-d)}\oplus S_{\btau i}^{\oplus (m_1-e)},S_j}^{N\oplus S_i^{\oplus (t_1^M+e-m_2)} \oplus S_{\btau i}^{\oplus (t_3^M+d-n_2)}}=1=F_{S_i^{\oplus (m_1-d)}, S_{j}}^{N_1\oplus S_i^{\oplus (t_1^M+e-m_2)}}=F_{S_{\btau i}^{\oplus (m_1-e)}, S_{j}}^{N_2\oplus S_{\btau i}^{\oplus (t_3^M+d-n_2) }},
\end{align*}
if they are nonzero.

We have
\begin{align*}
(q-1)|&\Aut(N\oplus S_i^{\oplus (t_1^M+e-m_2)} \oplus S_{\btau i}^{\oplus (t_3^M+d-n_2)})|
\\
&=|\Aut(N_1\oplus S_i^{\oplus (t_1^M+e-m_2)})|\cdot |\Aut(N_2\oplus S_{\btau i}^{\oplus (t_3^M+d-n_2) })|.
\end{align*}
Using the Riedtman-Peng formula, we have
\begin{align*}
\big|\Ext^1_{kQ} & (S_i^{\oplus (m_1-d)}\oplus S_{\btau i}^{\oplus (m_1-e)}, S_{j})_{N\oplus S_i^{\oplus (t_1^M+e-m_2)} \oplus S_{\btau i}^{\oplus (t_3^M+d-n_2)}} \big|
\\
&= \frac{\big|\Aut(S_i^{\oplus (m_1-d)}\oplus S_{\btau i}^{\oplus (m_1-e)})\big|\cdot |\Aut(S_j)|}{\big|\Aut(N\oplus S_i^{\oplus (t_1^M+e-m_2)} \oplus S_{\btau i}^{\oplus (t_3^M+d-n_2)})\big|}
\\
&= \frac{\big|\Aut(S_i^{\oplus (m_1-d)})\big|\cdot \big|\Aut(S_{\btau i}^{\oplus (m_1-e)})\big|\cdot |\Aut(S_j)|^2}{\big|\Aut(N_1\oplus S_i^{\oplus(t_1^M+e-m_2)})\big|\cdot \big|\Aut(N_2\oplus S_{\btau i}^{\oplus(t_3^M+d-n_2)})\big|}
\\
&= \big|\Ext^1_{kQ}(S_i^{\oplus (m_1-d)}, S_{j})_{N_1\oplus S_i^{\oplus (t_1^M+e-m_2)}}\big|\cdot \big|\Ext^1_{kQ}(S_{\btau i}^{\oplus (m_1-e)}, S_{j})_{N_2\oplus S_{\btau i}^{\oplus (t_3^M+d-n_2) }}\big|.
\end{align*}
The lemma is proved.
\end{proof}

\subsection{The proof}

Now we can complete the proof of Proposition \ref{prop:build-block}.
By Lemma \ref{lem:decom1} we have
\begin{align*}
&\big|\Ext^1_{kQ}(S_i^{\oplus (m_1-d)}\oplus S_{\btau i}^{\oplus (n_1-e)}, S_{j}\oplus S_i^{\oplus (m_2-e)}\oplus S_{\btau i}^{\oplus (n_2 -d) })_M  \big|\\
=&\big|\Ext^1_{kQ}(S_i^{\oplus (m_1-d)}\oplus S_{\btau i}^{\oplus (n_1-e)}, S_{j})_{N\oplus S_i^{\oplus (t_1^M+e-m_2)} \oplus S_{\btau i}^{\oplus(t_3^M+d-n_2)}} \big|\\
=&\big|\Ext^1_{kQ}(S_i^{\oplus (m_1-d)}, S_{j})_{N_1\oplus S_i^{\oplus (t_1^M+e-m_2)}}\big|\cdot \big|\Ext^1_{kQ}(S_{\btau i}^{\oplus (n_1-e)}, S_{j})_{N_2\oplus S_{\btau i}^{\oplus( t_3^M+d-n_2 )}}\big|.
\end{align*}
Together with Lemma \ref{lem:decom2} and Lemma \ref{lem:decom3}, we have
\begin{align}
\label{eq:CM}
|\cc_M|&= q^{-c_{ij}d-c_{\btau i,j}e}\big|\Ext^1_{kQ}(S_i^{\oplus (m_1-d)}\oplus S_{\btau i}^{\oplus (n_1-e)}, S_{j}\oplus S_i^{\oplus (m_2-e)}\oplus S_{\btau i}^{\oplus (n_2-d) })_M  \big|
 \\
&\cdot\big|\Ext_{\Lambda^\imath}^1 (S_i^{\oplus m_1},S_{\btau i}^{\oplus n_2})_{\E_i^{\oplus d} \oplus S_i^{\oplus (m_1-d)} \oplus S_{\btau i}^{\oplus (n_2-d)} }\big|\cdot \big|\Ext_{\Lambda^\imath}^1 (S_{\btau i}^{\oplus n_1},S_i^{\oplus m_2})_{\E_{\btau i}^{\oplus e} \oplus S_i^{\oplus (m_2-e)} \oplus S_{\btau i}^{\oplus (n_1-e)} }\big|.
\notag
\end{align}

By a standard computation (see e.g. \cite{Rin90}), we have
\begin{align*}
&F_{S_i^{\oplus (m_1-d)}\oplus S_{\btau i}^{\oplus (n_1-e)}, S_{j}\oplus S_i^{\oplus (m_2-e)}\oplus S_{\btau i}^{\oplus (n_2-d) }}^{ M}\\
=& \sqq^{(t_1^M-(m_2-e))(m_2-e)} \qbinom{t_1^M}{m_2-e}_\sqq \cdot  \sqq^{(t_3^M-(n_2-d))(n_2-d)} \qbinom{t_3^M}{n_2-d}_\sqq.
\end{align*}
Using \eqref{Ried-P}, one can obtain that
\begin{align*}
&\big|\Ext^1_{kQ}(S_i^{\oplus (m_1-d)}\oplus S_{\btau i}^{\oplus (n_1-e)}, S_{j}\oplus S_i^{\oplus (m_2-e)}\oplus S_{\btau i}^{\oplus (n_2-d) })_M  \big|\\
%=& | \Ext^1_{\Lambda^\imath}( S_{\btau i}^{\oplus (r-x-d)} \oplus S_i^{\oplus(s-x-e)},  S_{j}\oplus S_i^{\oplus (a-u-s-y-d)}\oplus S_{\btau i}^{\oplus (a-u-r-y-e)})_{M}  |\\
=& \prod_{i=0}^{m_1-d-1}(q^{m_1-d}-q^i)\prod_{i=0}^{n_1-e-1}(q^{n_1-e}-q^i)  \prod_{i=0}^{m_2-e-1}(q^{m_2-e}-q^i)\prod_{i=0}^{n_2-d-1}(q^{n_2-d}-q^i)\times\\
& \sqq^{ 2(m_1-d)(m_2-e)+ 2(n_1-e)(n_2-d)}\sqq^{(t_1^M-(m_2-e))(m_2-e)+(t_3^M-(n_2-d))(n_2-d)} \times\\
&\qbinom{t_1^M}{m_2-e}_\sqq \qbinom{t_3^M}{n_2-d}_\sqq \frac{(q-1)}{|\aut(M)|}\\
=& \sqq^{(m_1-d)^2+{ m_1-d \choose 2} + (n_1-e)^2 + { n_1-e \choose 2}  + (m_2-e)^2 +{ m_2-e \choose 2} +(n_2-d)^2 + {n_2-d \choose 2} } \times\\
&\sqq^{ 2(m_1-d)(m_2-e)+ 2(n_1-e)(n_2-d)}\sqq^{(t_1^M-(m_2-e))(m_2-e)+(t_3^M-(n_2-d))(n_2-d)} (\sqq-\sqq^{-1})^{m_1+m_2+n_1+n_2-2d-2e} \times\\
&[m_1-d]_\sqq^! [n_1-e]_\sqq^! [m_2-e]_\sqq^! [n_2-d]_\sqq^!\qbinom{t_1^M}{m_2-e}_\sqq \qbinom{t_3^M}{n_2-d}_\sqq \frac{(q-1)}{|\aut(M)|}.
\end{align*}
Furthermore, we have
\begin{align}
\label{eq:rank d}
\big|\Ext_{\Lambda^\imath}^1 (S_i^{\oplus m_1},S_{\btau i}^{\oplus n_2})_{\E_i^{\oplus d} \oplus S_i^{\oplus (m_1-d)} \oplus S_{\btau i}^{\oplus (n_2-d)} }\big|
&= \frac{\prod_{i=0}^{d-1} (q^{m_1}-q^{i})   \prod_{i=0}^{d-1} (q^{n_2}-q^{i})}{\prod_{i=0}^{d-1}(q^d-q^i)}\\\notag
&= \sqq^{d(m_1+n_2-d)+ {d \choose 2} } (\sqq -\sqq^{-1})^d \qbinom{m_1}{d}_\sqq \qbinom{n_2}{d}_\sqq[d]_\sqq^!,
%\frac{ [m_1]_\sqq [m_1-1]_\sqq \ldots [m_1-d+1]_\sqq [n_2]_\sqq [n_2-1]_\sqq \ldots [n_2-d+1]_\sqq}{ [d]_\sqq^! };
\end{align}
and
\begin{align}
\notag
&\big|\Ext_{\Lambda^\imath}^1 (S_{\btau i}^{\oplus n_1},S_i^{\oplus m_2})_{\E_{\btau i}^{\oplus e} \oplus S_i^{\oplus (m_2-e)} \oplus S_{\btau i}^{\oplus (n_1-e)} }\big|=\sqq^{e(n_1+m_2-e)+ {e \choose 2} } (\sqq -\sqq^{-1})^e\qbinom{m_2}{e}_\sqq  \qbinom{n_1}{e}_\sqq[e]_\sqq^!. %\frac{ [n_1]_\sqq [n_1-1]_\sqq \ldots [n_1-e+1]_\sqq [m_2]_\sqq [m_2-1]_\sqq \ldots [m_2-e+1]_\sqq}{ [e]_\sqq^! }.
\end{align}

Plugging into \eqref{eq:CM}, we have
\begin{align*}
|\cc_M|=&\sqq^{-2c_{ij}d-2c_{\btau i,j}e} \sqq^{(m_1-d)^2+{ m_1-d \choose 2} + (n_1-e)^2 + { n_1-e \choose 2}  + (m_2-e)^2 +{ m_2-e \choose 2} +(n_2-d)^2 + {n_2-d \choose 2} } \\
&\quad \times \sqq^{d(m_1+n_2-d)+ {d \choose 2} +e(n_1+m_2-e)+ {e \choose 2}}
\sqq^{ 2(m_1-d)(m_2-e)+ 2(n_1-e)(n_2-d)}
\\
& \quad \times  \sqq^{(t_1^M-(m_2-e))(m_2-e)+(t_3^M-(n_2-d))(n_2-d)}
 (\sqq-\sqq^{-1})^{m_1+m_2+n_1+n_2-d-e}
\\
&\quad \times \frac{[m_1]_\sqq^![m_2]_\sqq^![n_1]_\sqq^![n_2]_\sqq^!}{[d]_\sqq^![e]_\sqq^!}\qbinom{t_1^M}{m_2-e}_\sqq \qbinom{t_3^M}{n_2-d}_\sqq \frac{(q-1)}{|\aut(M)|}
\\
=& \sqq^{-2c_{ij}d-2c_{\btau i,j}e} \sqq^{ 2(m_1-d)(m_2-e)+ 2(n_1-e)(n_2-d)} (\sqq-\sqq^{-1})^{m_1+m_2+n_1+n_2-d-e}\\
& \quad \times
\sqq^{p'(t_3^M,d,m_1,n_2)} \frac{[m_1]_\sqq^![n_2]_\sqq^!}{[d]_\sqq^!} \qbinom{t_3^M}{n_2-d}_\sqq \sqq^{p'(t_1^M,e,n_1,m_2)} \frac{[m_2]_\sqq^![n_1]_\sqq^!}{[e]_\sqq^!} \qbinom{t_1^M}{m_2-e}_\sqq\frac{(q-1)}{|\aut(M)|}.
\end{align*}
The desired formula \eqref{eq:built1} follows from the above formula and \eqref{eq:builtb1}. This completes the proof of Proposition~ \ref{prop:build-block}.

%\begin{proposition}
%   \label{prop:qHall=Ui}
%Let $(Q, \btau)$ be a Dynkin $\imath$quiver. Then there exists a $\Q({\sqq})$-algebra homomorphism
%\begin{align}
 %  \label{eqn:psi morphism}
%\widetilde{\Psi}: \tUi_{|v={\sqq}} &\longrightarrow \tMH,
%\end{align}
%which sends
%\begin{align}
%B_j \mapsto \frac{-1}{q-1}[S_{j}],\text{ if } j\in\ci,
%&\qquad\qquad
%\tk_i \mapsto - q^{-1}[\E_i], \text{ if }i=\btau i;
%  \label{eq:split}
%\\
%B_{j} \mapsto \frac{{\sqq}}{q-1}[S_{j}],\text{ if }j\notin \ci,
%\qquad\qquad
%\tk_i \mapsto [\E_i],\quad \text{ if }i\neq \btau i.
%  \label{eq:extra}
%\end{align}
%\end{proposition}

%\begin{lemma}[see e.g. \cite{M06}]
%\label{lem: number rank r}
%We have
%$$\sharp\{A\in M_{m\times n}(\F_q)\mid \rank A=r\}=\prod_{i=0}^{r-1}\frac{(q^m-q^i)(q^n-q^i)}{q^r-q^i}.$$
%\prod_{i=1}^r \frac{q^{\min\{m,n\}-i+1}-1}{q^i-1}\prod_{i=1}^r (q^{\max\{m,n\}}-q^{i-1}).$$
%\end{lemma}

%%%%\iffalse
%%%%
%%%%%%
%%%%%%
\section{Proof of the formula \eqref{eqn:braidsplitodd}}
\label{app:braidsplitodd}

In this appendix, we provide the details for the proof of the formula \eqref{eqn:braidsplitodd}.

\subsection{Computation of $[S'_i]_{\odd}^{(r)}*[S'_j] *[S'_i]_{\ov{1}+\ov{a}}^{(s)}$}
\label{subsec:braidsplitodd}

Let us first compute $[S'_i]_{\odd}^{(r)}*[S'_j] *[S'_i]_{\ov{1}+\ov{a}}^{(s)}$, depending on the parity of $r$.

\subsubsection{$r$ is even}
%{\bf Case (I)}. Assume $r$ is even.
For any $s\geq0$ such that $r+s+2t=a$ with $t\geq0$, we have by Lemma~\ref{lem:iDPev} and Lemma~\ref{lem:SSSM}
\begin{align*}
[S'_i]_{\odd}^{(r)}*[S'_j] *[S'_i]_{\ov{1}+\ov{a}}^{(s)}
&= \sum_{k=0}^{\frac{r}{2}}\frac{\sqq^{k(k+1)-\binom{r-2k}{2}} \cdot (\sqq-\sqq^{-1})^k }{[r-2k]_\sqq^{!}[2k]_\sqq^{!!}} [(r-2k)S'_i]*[\E'_i]^k*[S'_j]
\\
&\qquad * \sum_{m=0}^{\lfloor\frac{s}{2}\rfloor} \frac{\sqq^{m(m+1)-\binom{s-2m}{2}}\cdot (\sqq-\sqq^{-1})^{m}}{[s-2m]_\sqq^{!} [2m]_\sqq^{!!}} [(s-2m)S'_i]*[\E'_i]^m
\\
&= \sum_{k=0}^{\frac{r}{2}}\sum_{m=0}^{\lfloor\frac{s}{2}\rfloor} \frac{\sqq^{k(k+1)+m(m+1)-\binom{r-2k}{2} -\binom{s-2m}{2}}\cdot (\sqq-\sqq^{-1})^{k+m} }{[r-2k]_\sqq^![s-2m]_\sqq^![2k]_\sqq^{!!}[2m]_\sqq^{!!}}
\\
& \qquad\qquad \times [(r-2k)S'_i]*[S'_j]*[(s-2m)S'_i]*[\E'_i]^{k+m}
\\
&= \sum_{k=0}^{\frac{r}{2}}\sum_{m=0}^{\lfloor\frac{s}{2}\rfloor} \sum_{n=0}^{\min\{r-2k,s-2m\}} \sum_{[M]\in\mathcal{I}_{r+s-2k-2m-2n}}
\frac{\sqq^{k(k+1)+m(m+1)-\binom{r-2k}{2} -\binom{s-2m}{2}} }{[r-2k]_\sqq^![s-2m]_\sqq^![2k]_\sqq^{!!}[2m]_\sqq^{!!}}\\
&\qquad\qquad \times \sqq^{p(a,n,r-2k,s -2m)}
(\sqq -\sqq^{-1})^{r+s -k -m -n+1 } \frac{[r-2k]_\sqq^{!} [s -2m]_\sqq^{!}}{[n]_\sqq^{!}}
\\ &
\qquad\qquad \times \qbinom{u_M}{s -2m-n}_\sqq \frac{[M]}{|\aut(M)|}*[\E'_i]^{n+k+m}.
\end{align*}
This can be simplified to be
\begin{align*}
[S'_i]_{\odd}^{(r)}*[S'_j] *[S'_i]_{\odd}^{(s)}
&= \sum_{k=0}^{\frac{r}{2}}\sum_{m=0}^{\lfloor\frac{s}{2}\rfloor} \sum_{n=0}^{r-2k} \sum_{[M]\in\mathcal{I}_{r+s-2k-2m-2n}}
\frac{\sqq^{\cz(a,r,s,k,m,n)+2k+2m}
(\sqq-\sqq^{-1})^{r+s-k-m-n+1} }{[n]_\sqq^{!} [2k]_\sqq^{!!}[2m]_\sqq^{!!}}\\
&
\qquad\qquad\qquad\qquad \times \qbinom{u_M}{s -2m-n}_\sqq \frac{[M] *[\E'_i]^{n+k+m}}{|\aut(M)|}.
\end{align*}

\subsubsection{$r$ is even}

%{\bf Case (II)}. Assume $r$ is odd.
For any $s\geq0$ such that $r+s+2t=a$ with $t\geq0$, we have by Lemma~\ref{lem:iDPev} and Lemma~\ref{lem:SSSM}
\begin{align*}
[S'_i]_{\odd}^{(r)}*[S'_j] *[S'_i]_{\ov{1}+\ov{a}}^{(s)}
=&\sum_{k=0}^{\lfloor\frac{r}{2}\rfloor}\frac{\sqq^{k(k-1)-\binom{r-2k}{2}} \cdot (\sqq-\sqq^{-1})^k }{[r-2k]_\sqq^![2k]_\sqq^{!!}} [(r-2k)S'_i]*[\E'_i]^k*[S'_j]\\
& \qquad * \sum_{m=0}^{\lfloor\frac{s}{2}\rfloor} \frac{\sqq^{m(m-1)-\binom{s-2m}{2}}\cdot (\sqq-\sqq^{-1})^{m}}{[s-2m]_\sqq^{!}[2m]_\sqq^{!!}} [(s-2m)S'_i]*[\E'_i]^m\\
=&\sum_{k=0}^{\lfloor\frac{r}{2}\rfloor}\sum_{m=0}^{\lfloor\frac{s}{2}\rfloor} \sum_{n=0}^{\min\{r-2k,s-2m\}} \sum_{[M]\in\mathcal{I}_{r+s-2k-2m-2n}}
 \frac{\sqq^{k(k-1)+m(m-1)-\binom{r-2k}{2} -\binom{s-2m}{2}} }{[r-2k]_\sqq^![s-2m]_\sqq^![2k]_\sqq^{!!}[2m]_\sqq^{!!}}\\
&\qquad \times \sqq^{p(a,n,r-2k, s -2m)}
(\sqq -\sqq^{-1})^{r+s -k -m -n+1 } \frac{[r-2k]_\sqq^{!} [s -2m]_\sqq^{!}}{[n]_\sqq^{!}}
\\
& \qquad \times \qbinom{u_M}{s -2m-n}_\sqq \frac{[M]}{|\aut(M)|}*[\E'_i]^{n+k+m}.
\end{align*}
This can be simplified to be
\begin{align*}
&[S'_i]_{\odd}^{(r)}*[S'_j] *[S'_i]_{\ov{1}+\ov{a}}^{(s)}
= \sum_{k=0}^{\lfloor\frac{r}{2}\rfloor}\sum_{m=0}^{\lfloor\frac{s}{2}\rfloor} \sum_{n=0}^{r-2k} \sum_{[M]\in\mathcal{I}_{r+s-2k-2m-2n}}
  \\
& \qquad\qquad \frac{\sqq^{\cz(a,r,s,k,m,n)}
 (\sqq-\sqq^{-1})^{r+s -k -m -n+1} }{[n]_\sqq^{!} [2k]_\sqq^{!!}[2m]_\sqq^{!!}}
\qbinom{u_M}{s -2m-n}_\sqq \frac{[M] *[\E'_i]^{n+k+m}}{|\aut(M)|}.
\end{align*}

\subsection{Reduction for \eqref{eqn:braidsplitodd}}

Summing up the above two cases, we obtain
\begin{align}
\text{RHS}\eqref{eqn:braidsplitodd}
  \label{eq:long2}
&= \sum_{r=0,2\mid r}^{a} \sum_{k=0}^{\frac{r}{2}}\sum_{m=0}^{\lfloor\frac{a-r}{2}\rfloor} \sum_{n=0}^{r-2k} \sum_{[M]\in\mathcal{I}_{a-2k-2m-2n}}
 (-1)^a  (\sqq-\sqq^{-1})^{ -k -m -n+1}
 \\
& \qquad\qquad \times  \frac{\sqq^{r+\cz(a,r,a-r,k,m,n)+2k+2m-a} }{[n]_\sqq^{!} [2k]_\sqq^{!!}[2m]_\sqq^{!!}}
\qbinom{u_M}{a-r -2m-n}_\sqq \frac{[M] *[\E'_i]^{n+k+m}}{|\aut(M)|}
 \notag\\
&\quad -\sum_{r=0,2\nmid r}^{a}\sum_{k=0}^{\lfloor\frac{r}{2}\rfloor}\sum_{m=0}^{\lfloor\frac{a-r}{2}\rfloor} \sum_{n=0}^{r-2k} \sum_{[M]\in\mathcal{I}_{a-2k-2m-2n}}
(-1)^a  (\sqq-\sqq^{-1})^{ -k -m -n+1}
 \notag\\
& \qquad\qquad \times \frac{\sqq^{r+\cz(a,r,a-r,k,m,n)-a} }{[n]_\sqq^{!} [2k]_\sqq^{!!}[2m]_\sqq^{!!}}
\qbinom{u_M}{a-r -2m-n}_\sqq \frac{[M] *[\E'_i]^{n+k+m}}{|\aut(M)|}
 \notag\\
&\quad -\sum_{t\geq1}\sum_{r=0,2\nmid r}^{a-2t}\sum_{k=0}^{\lfloor\frac{r}{2}\rfloor}\sum_{m=0}^{\lfloor\frac{a-r}{2}\rfloor-t} \sum_{n=0}^{r-2k} \sum_{[M]\in\mathcal{I}_{a-2t-2k-2m-2n}}
(-1)^a  (\sqq-\sqq^{-1})^{ -k -m -n+1}
 \notag\\
& \quad \times \frac{\sqq^{r+ \cz(a,r,a-2t-r,k,m,n)  -a+2t} }{[n]_\sqq^{!} [2k]_\sqq^{!!}[2m]_\sqq^{!!}}
\qbinom{u_M}{a-2t-r -2m-n}_\sqq \frac{[M] *[\E'_i]^{n+k+m+t}}{|\aut(M)|}.
\notag
\end{align}

Fix
\[
d =t+k+m+n.
\]
 %and fix $[M]   \in\mathcal{I}_{a-2d}$.
 We have reduced the proof of \eqref{eqn:braidsplitodd} to  proving the coefficient of $\frac{[M] *[\E'_i]^{d}}{|\aut(M)|}$ of the RHS of \eqref{eq:long2} is 0 for any given $[M] \in \mathcal I_{a-2d}$ if not both $d$ and $u_M$ are zeros.

 %We have $n=d-k-m-t \geq0$.
Set $u=u_M$, and set
\begin{align}
\label{def:T'}
A'(a,d,u)
&= \sum_{t\geq0}\sum_{r=0,2\nmid r}^{a-2t}\sum_{k=0}^{\lfloor\frac{r}{2}\rfloor}\sum_{m=0}^{\lfloor\frac{a-r}{2}\rfloor-t} \delta\{0\le n \le r-2k\}
\\\notag
&\qquad \times \frac{\sqq^{r+\cz(a,r,a-2t-r,k,m,n)-a+2t}
(\sqq-\sqq^{-1})^{-k-m-n+1} }{[n]_\sqq^{!} [2k]_\sqq^{!!}[2m]_\sqq^{!!}}
\qbinom{u}{a-2t-r -2m-n}_\sqq
\\\notag
&\quad -\sum_{r=0,2\mid r}^{a} \sum_{k=0}^{\frac{r}{2}}\sum_{m=0}^{\lfloor\frac{a-r}{2}\rfloor}
  \delta\{0\le n \le r-2k\}
  \\
  & \qquad \times \frac{\sqq^{r+\cz(a,r,a-r,k,m,n)+2k+2m-a}
 (\sqq-\sqq^{-1})^{ -k -m -n+1} }{[n]_\sqq^{!} [2k]_\sqq^{!!}[2m]_\sqq^{!!}}
\qbinom{u}{a -r-2m-n}_\sqq.
 \notag
\end{align}
%for any $a\geq0$, $0\leq d\leq \frac{a}{2}$, $0\leq u\leq a-2d$, %$d$ and $u$ not both zero,
See \eqref{eq:z} for $\cz(\cdot, \cdot, \cdot, \cdot, \cdot, \cdot)$ and also see \eqref{eq:p} for $p(\cdot, \cdot, \cdot, \cdot)$. %Here $\delta\{X\}=1$ if the statement $X$ holds and $=0$ if $X$ is false.

Then the coefficient of $\frac{[M] *[\E'_i]^{d}}{|\aut(M)|}$ of the RHS of \eqref{eq:long2} is $(-1)^a A'(a,d,u_M)$. Summarizing, we have established the following (which is the counterpart of Proposition \ref{prop:GammaA}).

\begin{proposition}
 \label{prop:GammaB}
The formula \eqref{eqn:braidsplitodd} is equivalent to the following identity
\begin{equation}
\label{eq:B=0}
A'(a,d,u) =0,
\end{equation}
for non-negative integers $a,d,u$ subject to the constraints:
\begin{equation}
0\leq d\leq \frac{a}{2}, \quad 0\leq u \leq a-2d, \quad  d \text{ and $u$ not both zero}.
\end{equation}
\end{proposition}

\subsection{Reduction for the identity \eqref{eq:B=0} }
\label{subsec:A'=0}

We shall denote the 2 summands in $A' =A'(a,d,u)$ in \eqref{def:T'} as $A'_0, A'_1$, and thus
\[
A' =A'_0 -A'_1.
\]
Compare with \eqref{eq:defA}.
We also denote
\begin{align}
w &=r-2k-n. %,  \qquad a=a.
\end{align}
Set
\begin{align}
d= k+m+n +t
\end{align}
in the $A'_0$ side, and $d =k+m+n$ in the $A'_1$ side.
Using the same argument as in \S\ref{subsec:A=0}, for $d>0$, we see that the identity \eqref{eq:B=0} for $d>0$ is equivalent to the following identity
\begin{align}
 \label{eq:coeff2b}
\sum_{\stackrel{w+n \blue{\text{ odd}}}{t+k+m+n =d}}
\frac{\sqq^{t^2 -2dt +t  +2nt +\binom{n+1}{2} -2km -2m}}{[n]^!_\sqq [2k]^{!!}_\sqq [2m]^{!!}_\sqq} (\sqq -\sqq^{-1})^{t}
- \sum_{\stackrel{w+n \blue{\text{ even}}}{k+m+n =d}}
\frac{\sqq^{\binom{n+1}{2} -2km +2k}}{[n]^!_\sqq [2k]^{!!}_\sqq [2m]^{!!}_\sqq} =0.
\end{align}
The identity \eqref{eq:coeff2b} is clearly equivalent to the identity \eqref{eq:coeff} (by switching the parity of $w$), which was established in Section~\ref{sec:identities}.

The identity \eqref{eq:B=0} for $d=0$ holds exactly in the same way as for \eqref{eq:d=0} (up to an irrelevant overall sign).
Therefore, the identity \eqref{eq:B=0} is fully established, and then the formula \eqref{eqn:braidsplitodd}
follows by Proposition~\ref{prop:GammaB}.

%%%%\fi
%%%%

%%%%%%%
%%%%%%%


\begin{thebibliography}{CLW21c}

\bibitem{BW18} H. Bao and W. Wang,
{\em  A new approach to Kazhdan-Lusztig theory  of type $B$ via quantum symmetric pairs}, Ast\'erisque {\bf 402}, 2018, vii+134pp, \href{https://arxiv.org/abs/1310.0103}{arXiv:1310.0103v3}

\bibitem{BW18b} H. Bao and W. Wang,
{\em Canonical bases arising from quantum symmetric pairs}, Invent. Math. {\bf 213} (2018), 1099--1177.

\bibitem{BK20} P. Baseilhac and S. Kolb,
{\em Braid group action and root vectors for the $q$-Onsager algebra}, Transform. Groups {\bf 25} (2020),  363--389.
%\href{https://arxiv.org/abs/1706.08747}{arXiv:1706.08747}

\bibitem{BeW18} C. Berman and W. Wang, {\em Formulae of $\imath$-divided powers in ${\mathbf U}_q(\mathfrak{sl}_2)$}, J. Pure Appl. Algebra {\bf 222} (2018), 2667--2702.

%\bibitem[BGP73]{BGP73} I.N. Bernstein, I.M. Gelfand, and V.A. Ponomarev,
%{\em Coxeter functors and Gabriel's theorem}, Uspehi. Mat. Nauk {\bf 28}  (1973), no. 2 (170), 19--33.

\bibitem{Br13} T. Bridgeland,
{\em Quantum groups via Hall algebras of complexes}, Ann. Math. {\bf 177} (2013), 739--759.

\bibitem{Ch07} L. Chekhov,
{\em Teichm\"uller theory of bordered surfaces}, SIGMA Symmetry Integrability Geom. Methods Appl. {\bf 3} (2007),
Paper 066, 37 pp.

\bibitem{CLW18} X. Chen, M. Lu and W. Wang,
{\em A Serre presentation for the $\imath$quantum groups}, Transform. Groups {\bf26} (2021), 827--857,  \href{http://arxiv.org/abs/1810.12475}{arxiv:1810.12475}

\bibitem{CLW21} X. Chen, M. Lu and W. Wang,
{\em Serre-Lusztig relations for $\imath$quantum groups}, Commun. Math. Phys. {\bf 382} (2021), 1015--1059,
\href{http://arxiv.org/abs/2001.03818}{arxiv:2001.03818}


\bibitem{CLW21c} X. Chen, M. Lu and W. Wang,
{\em Serre-Lusztig relations for $\imath$quantum groups III}, J. Pure Appl. Algebra (to appear),  \href{http://arxiv.org/abs/2106.06888}{arxiv:2106.06888}



%\bibitem[DDPW08]{DDPW}
%B.~Deng, J.~Du, B.~Parshall and J.~Wang,
%{\em Finite dimensional algebras and quantum groups},
%Mathematical Surveys and Monographs {\bf 150}.
%American Mathematical Society, Providence, RI, 2008.

\bibitem{DR74} V. Dlab and C.M. Ringel,
{\em Representations of graphs and algebras}, Carleton Mathematical Lecture
Notes  {\bf 8}, Carleton University, Ottawa, 1974, iii+86 pp.

\bibitem{Dob20} L. Dobson,
{\em Braid group actions for quantum symmetric pairs of type AIII/AIV}, J. Algebra {\bf 564} (2020), 151--198.
%\href{https://arxiv.org/abs/1909.11215}{arxiv:1909.11215}


\bibitem{GLS17} C. Geiss, B. Leclerc and J. Schr\"{o}er,
{\em Quivers with relations for symmetrizable Cartan matrices I: Foundations}, Invent. Math. {\bf 209} (2017), 61--158.

%\bibitem[Gor13]{Gor13} M. Gorsky,
%{\em Semi-derived Hall algebras and tilting invariance of Bridgeland-Hall algebras}, \href{https://arxiv.org/abs/1303.5879}{arXiv:1303.5879v2}

\bibitem{Gor18} M. Gorsky,
{\em Semi-derived and derived Hall algebras for stable categories}, IMRN, Vol. {\bf 2018}, No.~1,  138--159. %, \href{https://arxiv.org/abs/1409.6798}{arXiv:1409.6798}

\bibitem{Gr95} J.A. Green,
{\em Hall algebras, hereditary algebras and quantum groups}, Invent. Math. {\bf 120} (1995), 361--377.

\bibitem{KR90} A.N. Kirillov and N. Reshetikhin,
{\em q-Weyl group and a multiplicative formula for universal R-matrices}, Commun. Math. Phys. {\bf 134} (1990), 421--431.

\bibitem{Ko14} S. Kolb,
{\em Quantum symmetric Kac-Moody pairs}, Adv. Math. {\bf 267} (2014), 395--469.



\bibitem{KP11} S. Kolb and J. Pellegrini,
{\em Braid group actions on coideal subalgebras of quantized enveloping algebras}, J. Algebra {\bf 336} (2011), 395--416.

\bibitem{Let99} G. Letzter,
{\em Symmetric pairs for quantized enveloping algebras}, J. Algebra {\bf 220} (1999), 729--767.

\bibitem{Let02}
G. Letzter,
{\em Coideal subalgebras and quantum symmetric pairs},
New directions in Hopf algebras (Cambridge), MSRI publications, {\bf 43}, Cambridge Univ. Press, 2002, pp. 117--166.

\bibitem{LS90} S. Levendorskii and I. Soibelman,
{\em Some applications of quantum Weyl groups}, J. Geom. and Phys. {\bf 7} (1990), 241--254.

\bibitem{Li12} F. Li,
{\em Modulation and natural valued quiver of an algebra}, Pacific J. Math. {\bf 256} (2012), 105--128.

\bibitem{Lu19} M. Lu,
{\em Appendix A to  \cite{LW19}}, \href{http://arxiv.org/abs/1901.11446}{arXiv:1901.11446}

\bibitem{LP21} M. Lu and L. Peng,
{\em Semi-derived Ringel-Hall algebras and Drinfeld double}, Adv. Math. {\bf 383} (2021), 107668. %,  \href{https://arxiv.org/abs/1608.03106}{arXiv:1608.03106}





\bibitem{LW20a} M. Lu and W. Wang,
{\em Hall algebras and quantum symmetric pairs of Kac-Moody type}, \href{http://arxiv.org/abs/2006.06904}{ arXiv:2006.06904}

\bibitem{LW20b} M. Lu and W. Wang,
{\em A Drinfeld type presentation of affine $\imath$quantum groups I: split ADE type}, Adv. Math. {\bf393} (2021), 108111,  \href{https://arxiv.org/abs/2009.04542}{arXiv:2009.04542}

\bibitem{LW21} M. Lu and W. Wang, {\em Hall algebras and quantum symmetric pairs II: reflection functors},  Commun. Math. Phys.  {\bf 381} (2021), 799--855, \href{http://arxiv.org/abs/1904.01621}{arXiv:1904.01621}

\bibitem{LW19} M. Lu and W. Wang,
{\em Hall algebras and quantum symmetric pairs I: foundations},  Proc. Lond. Math. Soc. (3) {\bf124} (2022), 1--82,  \href{http://arxiv.org/abs/1901.11446}{arXiv:1901.11446}

\bibitem{Lus90} G. Lusztig,
{\em Canonical bases arising from quantized enveloping algebras}, J. Amer. Math. Soc. {\bf 3} (1990),  447--498.

\bibitem{Lus90b} G.~ Lusztig,
{\em Quantum groups at roots of $1$}, Geom. Dedicata {\bf 35} (1990), 89--114. 

\bibitem{Lus93} G. Lusztig, Introduction to Quantum Groups, Birkh\"{a}user, Boston, 1993.

\bibitem{Lus03}
G. Lusztig,
{\em Hecke algebras with unequal parameters},
CRM Monograph Series {\bf 18}, Amer. Math. Soc., Providence, RI, 2003, \href{https://arxiv.org/abs/math/0208154}{arXiv:0208154v2}

\bibitem{MR08} A. Molev and E. Ragoucy,
{\em Symmetries and invariants of twisted quantum algebras and associated Poisson algebras}, Rev. Math. Phys. {\bf 20} (2008), 173--198.

\bibitem{Rin90} C.M. Ringel,
{\em Hall algebras and quantum groups}, Invent. Math. {\bf 101} (1990), 583--591.

\bibitem{Rin96} C.M. Ringel,
{\em PBW-bases of quantum groups}, J. Reine Angrew. Math. {\bf 470} (1996), 51--88.

\bibitem{SV99} B. Sevenhant and M. Van den Bergh,
{\em On the double of the Hall algebra of a quiver}, J. Algebra {\bf 221} (1999), 135--160.

\bibitem{T06} B. T\"oen,
{\em Derived Hall algebras}, Duke Math. J. {\bf 135} (2006), 587--615.

\bibitem{X97} J. Xiao,
{\em Drinfeld double and Ringel-Green theory of Hall algebras}, J. Algebra {\bf 190} (1997), 100--144.

\bibitem{XX08}  J. Xiao and F. Xu,
{\em Hall algebras associated to triangulated categories}, Duke Math. J. {\bf 143} (2008), 357--373.

\bibitem{XY01} J. Xiao and S. Yang,
{\em BGP-reflection functors and Lusztig's symmetries: A Ringel-Hall approach to quantum groups}, J.
Algebra {\bf 241} (2001), 204--246.

%\bibitem[Ter17]{Ter17} P. Terwilliger, {\em The Lusztig automorphism of the $q$-Onsager algebra},  J. Algebra {\bf 506} (2018), 56--75. %\href{http://arxiv.org/abs/1706.05546}{arXiv:1706.05546}

\end{thebibliography}
\end{document}